# RELATIVE PERFORMANCE OF EXPECTED AND OBSERVED FISHER INFORMATION IN COVARIANCE ESTIMATION FOR MAXIMUM LIKELIHOOD ESTIMATES

by
Xumeng Cao

A dissertation submitted to Johns Hopkins University in conformity with the requirements for the degree of Doctor of Philosophy

# Abstract

Maximum likelihood estimation is a popular method in statistical inference. As a way of assessing the accuracy of the maximum likelihood estimate (MLE), the calculation of the covariance matrix of the MLE is of great interest in practice. Standard statistical theory shows that the normalized MLE is asymptotically normally distributed with covariance matrix being the inverse of the Fisher information matrix (FIM) at the unknown parameter. Two commonly used estimates for the covariance of the MLE are the inverse of the observed FIM (the same as the inverse Hessian of the negative log-likelihood) and the inverse of the expected FIM (the same as the inverse FIM). Both of the observed and expected FIM are evaluated at the MLE from the sample data. In this dissertation, we demonstrate that, under reasonable conditions similar to standard MLE conditions, the inverse expected FIM outperforms the inverse observed FIM under a mean squared error criterion. Specifically, in an asymptotic sense, the inverse expected FIM (evaluated at the MLE) has no greater mean squared error with respect to the true covariance matrix than the inverse observed FIM (evaluated at the MLE) at the element level. This result is different from widely accepted results showing preference for the observed FIM. In this dissertation, we present theoretical derivations that lead to the conclusion above. We also present numerical studies on three distinct problems to support the theoretical result.

This dissertation also includes two appendices on topics of relevance to stochastic systems. The first appendix discusses optimal perturbation distributions for the simultaneous perturbation stochastic approximation (SPSA) algorithm. The second appendix considers Monte Carlo methods for computing FIMs when closed forms are not attainable.

First Reader and Advisor: James C. Spall

Second Reader: Daniel Q. Naiman

# Acknowledgements

Foremost, I want to express my deepest gratitude to my advisor Prof. James Spall for his continuous supervision and support throughout my Ph.D. study and research at Hopkins. Prof. Spall's immense knowledge and professionalism, as well as his great personality, has guided me through this journey. He unconditionally shares his experience to provide intellectual supervision, which has had significant impact on my research achievement. He motivates me to explore with independent thinking and development. I'm grateful for his belief in me to tackle the challenges along the way. Without his motivation and encouragement, I would not have achieved all this. I'm always impressed by his discipline of rigorous conduct in scientific research, which I believe will continue to influence me throughout my career and life.

Besides, I want to greatly thank the second reader Prof. Daniel Naiman. His broad knowledge and experience in statistics has had been a key input to my research success. He and I had a conversation after my GBO (Graduate Board Oral) exam. He kindly shared his understanding of the topic, which broadened my thinking and approach in research. The book (McCullagh, 1987) that he shared has been a major reference to my research development. It gave me great inspiration in developing the approach to solve the main problem of my dissertation.

My sincere appreciation also goes to the rest of the dissertation committee: Prof. Charles Rohde, and Prof. Fred Torcaso. I want to thank each one of them for the insightful discussions with me during my research. Their willingness to share and inspire has been a great support on the completion of my dissertation. Besides, I truly appreciate their time and effort on reviewing my work and attending my defense.

A special thank you should go to Applied Physics Lab (APL) of Johns Hopkins University. I benefited from research funds from APL for several semesters, which has been a great support for my Ph.D. study. I also would like to thank the Duncan Fund offered by the Department of Applied Math and Statistics, JHU. This fund has financially supported me multiple times to attend various conferences. It creates great opportunity for me to share my work with other people in the community as well as to get exposure to frontier discoveries and developments in the field.

I also want to thank my beloved parents. They always believe in me, even when I don't. I would not have made it without their unconditional support and encouragement since the day I was born.

Last but not least, I want to thank my boyfriend, William. He is always there for me with love and support. He cheers for me during good times and backs me up during bad ones. I'm grateful to have him in my life.

# Table of contents

# List of Tables

# List of Figures

# Chapter 1

# Introduction

In this introduction chapter, we start with the motivation that drives our interest in the topic of Fisher information, which is followed by the literature review, where we summarize relevant work done by others. In Section 1.3, we propose our approach to solve the problem of interest. A sketch of our scheme is summarized at a high level. In the last section, some large sample results are discussed as a background for further analysis in Chapter 2.

## 1.1 Motivation

Maximum likelihood (ML) estimation is a standard approach for parameter estimation and statistical inference in the modeling of stochastic systems. Given a set of sample observations and the proposed underlying statistical model, the method of ML selects values of the model parameters that produce a distribution that gives the observed data the greatest probability. ML estimation enjoys great popularity in practice because it has many optimal properties such as asymptotic normality, functional invariance, and convergence to the true parameter in a certain probability sense. Not all of these properties are shared with other parameter estimation methods such as least-squares estimation (LSE).

Because of the nice properties it possesses, ML estimation is commonly used across a wide range of statistical models are fitted to real-life situations. To name a few examples, the generalized linear model (GLM), which is a generalization of ordinary linear regression that allows for response variables that have other than a normal distribution, is extensively used in various industries such as clinical trials, customer relationship marketing, and quantitative finance. The ML method is the standard approach used in practice to estimate parameters associated with GLMs, see Nelder and Wedderburn (1972). ML estimation can also be applied to hypothesis testing (Huelsenbeck and Crandall, 1997). The construction of a likelihood ratio test statistic is based on the idea of ML under the null hypothesis and the alternative hypothesis. System identification, where statistical methods are used in control engineering to build mathematical models of dynamical systems, is another area where ML is commonly seen. Particularly, system parameters are estimated using ML (Prado, 1979; Johari, et. al., 1965; and Ljung, 1999).

ML estimation produces point estimates based on sample data. Like many other point estimation methods, e.g., LSE, people are also interested in the accuracy of maximum likelihood estimates (MLEs). As a way of assessing the accuracy of MLEs, calculations of the associated confidence intervals are of great interest in statistical inference on unknown parameters (e.g., Ljung, 1999, pp. 215–218). At the same level of confidence, tighter confidence intervals indicate better accuracy of the corresponding MLEs and vice versa. To construct confidence intervals for MLEs, typically, one needs to know the distribution and the covariance matrix of the MLE. In fact, under regularity conditions, MLEs are asymptotically normally distributed. Given this asymptotic distribution, the problem of constructing confidence intervals is essentially reduced to finding the covariance matrix of MLEs.

Before we elaborate on the above statement, let us define the relevant notations. Let $\boldsymbol{X} = [\boldsymbol{X}_1, \boldsymbol{X}_2, \ldots, \boldsymbol{X}_n]$ be a sequence of $n$ independent but not necessarily identically distributed (i.n.i.d.) random vectors (variables) where each $\boldsymbol{X}_i$ may contain discrete or continuous components. The probability density/mass function of $\boldsymbol{X}_i$, say $p_i(\boldsymbol{x}_i, \boldsymbol{\theta})$, depends on a $p \times 1$ vector of unknown parameters $\boldsymbol{\theta} = [t_1, t_2, \ldots, t_p]^T$, where $\boldsymbol{\theta} \in \Theta$ and $\Theta$ is a $p$-dimensional parameter space. Let $\hat{\boldsymbol{\theta}}_n$ be an MLE for $\boldsymbol{\theta}$ based on $\boldsymbol{X}$ and the true value of $\boldsymbol{\theta}$ in the underlying distribution be $\boldsymbol{\theta}^*$. We use the notation $t_i$ to denote the $i$th component of $\boldsymbol{\theta}$ because we reserve $\hat{\boldsymbol{\theta}}_n$ for MLEs derived from a sample of size $n$. The joint probability density/mass function of $\boldsymbol{X}$ is $p(\boldsymbol{x}, \boldsymbol{\theta}) \equiv \prod_{i=1}^{n} p_i(\boldsymbol{x}_i, \boldsymbol{\theta})$. If we denote the negative log-likelihood function as $l(\boldsymbol{\theta}, \boldsymbol{x}) = -\log p(\boldsymbol{x}, \boldsymbol{\theta})$, the $p \times p$ Fisher information matrix (FIM) $\boldsymbol{F}_n(\boldsymbol{\theta})$ is defined as

$$\boldsymbol{F}_n(\boldsymbol{\theta}) \equiv E\left(\frac{\partial \log p(\boldsymbol{x}, \boldsymbol{\theta})}{\partial \boldsymbol{\theta}} \times \frac{\partial \log p(\boldsymbol{x}, \boldsymbol{\theta})}{\partial \boldsymbol{\theta}^T}\right)$$

$$= E\left(\frac{\partial\left(-l(\boldsymbol{\theta}, \boldsymbol{x})\right)}{\partial\boldsymbol{\theta}} \times \frac{\partial\left(-l(\boldsymbol{\theta}, \boldsymbol{x})\right)}{\partial\boldsymbol{\theta}^T}\right)$$

$$= E\left(\frac{\partial l(\boldsymbol{\theta}, \boldsymbol{x})}{\partial\boldsymbol{\theta}} \times \frac{\partial l(\boldsymbol{\theta}, \boldsymbol{x})}{\partial\boldsymbol{\theta}^T}\right), \tag{1.1}$$

where $\boldsymbol{\theta}^T$ is the transpose of $\boldsymbol{\theta}$, all expectations are taken with respect to data $\boldsymbol{X}$ and are conditional on the true parameter $\boldsymbol{\theta}^*$. The $p \times p$ Hessian matrix of $l(\boldsymbol{\theta}, \boldsymbol{x})$, $\boldsymbol{H}_n(\boldsymbol{\theta})$, is defined as the second derivative of $l(\boldsymbol{\theta}, \boldsymbol{x})$ with respect to $\boldsymbol{\theta}$:

$$\boldsymbol{H}_n(\boldsymbol{\theta}) \equiv \frac{\partial^2 l(\boldsymbol{\theta}, \boldsymbol{x})}{\partial\boldsymbol{\theta}\partial\boldsymbol{\theta}^T}.$$

Computation of $\boldsymbol{F}_n(\boldsymbol{\theta})$ according to its definition in (1.1) is often formidable because it involves direct calculation of expectation of an outer product form. Under some regularity conditions where the interchange of differentiation and integral is valid (more details are discussed in Chapter 2), $\boldsymbol{F}_n(\boldsymbol{\theta})$ has the following form equivalent to (1.1):

$$\boldsymbol{F}_n(\boldsymbol{\theta}) = E(\boldsymbol{H}_n(\boldsymbol{\theta})), \tag{1.2}$$

where the expectation is taken with respect to $\boldsymbol{X}$ and is conditional on the true parameter $\boldsymbol{\theta}^*$. Expression (1.2) provides an alternative of computing $\boldsymbol{F}_n(\boldsymbol{\theta})$, which is often more computationally friendly than the definition in (1.1).

Standard statistical theory shows that the normalized $\hat{\boldsymbol{\theta}}_n$ from either i.i.d. or i.n.i.d samples is asymptotically Gaussian under some reasonable conditions (Ljung, 1999, pp. 215–218 and Spall, 2003, Sect. 13.3). That is, under modest assumptions (more details discussed in Section 1.4.1),

$$\sqrt{n}(\hat{\boldsymbol{\theta}}_n - \boldsymbol{\theta}^*) \xrightarrow{\text{dist}} N\left(0,\ \bar{\boldsymbol{F}}(\boldsymbol{\theta}^*)^{-1}\right), \tag{1.3}$$

where "$\xrightarrow{\text{dist}}$" denotes convergence in distribution, and $\bar{\boldsymbol{F}}(\boldsymbol{\theta}^*) \equiv \lim_{n\to\infty} \boldsymbol{F}_n(\boldsymbol{\theta}^*)/n$. The superscript "−1" in (1.3) denotes matrix inverse. The asymptotic normality in (1.3) for i.n.i.d samples is of particular interest in our discussion below.

Given the asymptotic normality of MLEs, the problem of constructing confidence intervals reduces largely to the problem of determining the covariance matrix of MLEs, which is the main focus of this dissertation. In fact, other than the essential role in computing confidence intervals, the estimation of the covariance matrix of MLEs is also crucial in other applications. For example, in Nie and Yang (2005), the covariance matrix of MLEs is used in the discussion of the consistency of MLEs. Another example lies in the standard *t*-test, which is used to assess the significance of a parameter. Estimation of the covariance of the MLE is needed in computing the test statistic and the associated *P*-value, which is the probability of obtaining a test statistic at least as extreme as the one that was actually observed, assuming the null hypothesis is true.

In practical applications, one of two matrices is commonly used to approximate the covariance matrix of MLE: $\bar{\boldsymbol{F}}_n(\hat{\boldsymbol{\theta}}_n)^{-1}$ or $\bar{\boldsymbol{H}}_n(\hat{\boldsymbol{\theta}}_n)^{-1}$, where $\bar{\boldsymbol{F}}_n(\boldsymbol{\theta}) \equiv \boldsymbol{F}_n(\boldsymbol{\theta})/n$ and $\bar{\boldsymbol{H}}_n(\boldsymbol{\theta}) \equiv \boldsymbol{H}_n(\boldsymbol{\theta})/n$. Both of the two estimates are evaluated at the MLE $\hat{\boldsymbol{\theta}}_n$. The derivation of these two estimates is not surprising given the covariance term in the right hand side of (1.3) and the relation in (1.2). However, there is not yet a solid theoretical validation for the better choice between $\bar{\boldsymbol{F}}_n(\hat{\boldsymbol{\theta}}_n)^{-1}$ and $\bar{\boldsymbol{H}}_n(\hat{\boldsymbol{\theta}}_n)^{-1}$. In fact, people in practice tend to choose one or the other, depending on which one is easier to obtain for their problems. For instance, in Rice (1995, p. 269), $\bar{\boldsymbol{F}}_n(\hat{\boldsymbol{\theta}}_n)^{-1}$ is used to estimate the variance of the MLE based on i.i.d. Poisson distribution, where the closed form of $\bar{\boldsymbol{F}}_n(\boldsymbol{\theta})$ is easy to compute. Abt and Welch (1998) uses $\bar{\boldsymbol{F}}_n(\hat{\boldsymbol{\theta}}_n)^{-1}$ to estimate the covariance matrix of MLE in Gaussian stochastic processes. Escobar and Meeker (2001) discusses asymptotic equivalent performance of $\bar{\boldsymbol{F}}_n(\hat{\boldsymbol{\theta}}_n)^{-1}$ as an estimate of covariance of MLE for censored data (i.e. partially known observations) from location-scale families, which is a family of univariate probability distributions parameterized by a location parameter (e.g. mean of a normal distribution) and a non-negative scale parameter (e.g. variance of a normal distribution). On the other hand, Cavanaugh and Shumway (1996) mentions that in the setting of state-space models, where the structure of the Gaussian log-likelihood often makes $\bar{\boldsymbol{F}}_n(\hat{\boldsymbol{\theta}}_n)^{-1}$ difficult to compute, people prefer to use $\bar{\boldsymbol{H}}_n(\hat{\boldsymbol{\theta}}_n)^{-1}$ as an approximation of the covariance of MLE. In Prescott and Walden (1983), $\bar{\boldsymbol{H}}_n(\hat{\boldsymbol{\theta}}_n)^{-1}$ is used as an estimation of the covariance matrix of MLE from generalized extreme-value distributions for complete, left censored, right censored or doubly censored samples.

Other than the two estimates discussed above, there are other estimation methods in practical applications as well. For example, Jiang (2005) proposed an estimate of the covariance matrix that consists partially of $\bar{\boldsymbol{F}}_n(\hat{\boldsymbol{\theta}}_n)^{-1}$ and partially of $\bar{\boldsymbol{H}}_n(\hat{\boldsymbol{\theta}}_n)^{-1}$ for mixed linear models with non-normal data, where mixed models contain both fixed effects and random effects. In this case, a combination of $\bar{\boldsymbol{F}}_n(\hat{\boldsymbol{\theta}}_n)^{-1}$ and $\bar{\boldsymbol{H}}_n(\hat{\boldsymbol{\theta}}_n)^{-1}$ is used because the closed form of $\bar{\boldsymbol{F}}_n(\hat{\boldsymbol{\theta}}_n)^{-1}$ is not attainable and $\bar{\boldsymbol{H}}_n(\hat{\boldsymbol{\theta}}_n)^{-1}$ is inconsistent in the sense that it does not converge to the true covariance matrix in probability. Alternative estimation methods can be found in Royall (1986), Reeds (1978), and others. However, our discussion below mainly focuses on the relative performance of $\bar{\boldsymbol{F}}_n(\hat{\boldsymbol{\theta}}_n)^{-1}$ and $\bar{\boldsymbol{H}}_n(\hat{\boldsymbol{\theta}}_n)^{-1}$. Potential extension to other general estimation methods is left for future work.

Given the importance of covariance matrix estimation of MLE and the fact that no theoretical conclusion has been established for the best estimate, the aim of this work is to provide theoretical development for choosing a good estimate of the covariance matrix of MLE. In particular, we explore the properties of $\bar{\boldsymbol{F}}_n(\hat{\boldsymbol{\theta}}_n)^{-1}$ and $\bar{\boldsymbol{H}}_n(\hat{\boldsymbol{\theta}}_n)^{-1}$ and compare their performance in estimating the covariance matrix of a normalized $\hat{\boldsymbol{\theta}}_n$.

## 1.2 Literature review

There has been great interest and discussion in both observed and expected FIM, $\bar{\boldsymbol{H}}_n(\hat{\boldsymbol{\theta}}_n)$ and $\bar{\boldsymbol{F}}_n(\hat{\boldsymbol{\theta}}_n)$, in the literature. In this section, we review some of the work that is relevant to our discussion.

Efron and Hinkley (1978) appears to be the most-cited paper relative to comparing $\bar{\boldsymbol{H}}_n(\hat{\boldsymbol{\theta}}_n)$ and $\bar{\boldsymbol{F}}_n(\hat{\boldsymbol{\theta}}_n)$. Efron and Hinkley demonstrate that for scalar-parameter translation families with an appropriate ancillary statistic $a$ (more explanation below), the conditional variance of normalized $\hat{\theta}_n$ is better approximated by $\bar{H}_n(\hat{\theta}_n)^{-1}$ than by $\bar{F}_n(\hat{\theta}_n)^{-1}$. Specifically, the following ratio decays to zero in a stochastic sense:

$$\frac{\mathrm{var}\left(\sqrt{n}(\hat{\theta}_n - \theta^*)|a\right) - \bar{H}_n(\hat{\theta}_n)^{-1}}{\mathrm{var}\left(\sqrt{n}(\hat{\theta}_n - \theta^*)|a\right) - \bar{F}_n(\hat{\theta}_n)^{-1}}, \tag{1.4}$$

where var(·) denotes variance. Roughly speaking, if $n$ is large enough, the magnitude of error produced by $\bar{H}_n(\hat{\theta}_n)^{-1}$ is less than that produced by $\bar{F}_n(\hat{\theta}_n)^{-1}$ in some stochastic sense, i.e.,

$$\left|\mathrm{var}\left(\sqrt{n}(\hat{\theta}_n - \theta^*)|a\right) - \bar{H}_n(\hat{\theta}_n)^{-1}\right| < \left|\mathrm{var}\left(\sqrt{n}(\hat{\theta}_n - \theta^*)|a\right) - \bar{F}_n(\hat{\theta}_n)^{-1}\right|. \tag{1.5}$$

The ancillary statistic, $a$, in (1.4) and (1.5) is a statistic whose distribution does not depend on $\theta$ but which affects the precision of $\hat{\theta}_n$ as an estimate of $\theta$. An example of an ancillary statistic is given in Cox (1958), which we now summarize. An experiment is conducted to measure a constant $\theta$. Independent unbiased measurements $y$ of $\theta$ can be made with either of two instruments, both of which measure with normal error: instrument $k$ produces independent error that follows a $N(0, \sigma_k^2)$ distribution ($k = 1, 2$), where $\sigma_1^2$ and $\sigma_2^2$ are known and unequal. When a measurement $y$ is obtained, a record is also kept of the instrument used. In this case, the ancillary statistic is defined as the label for the instrument used for a particular observation, i.e., $a_j = k$ if $y_j$ is obtained using instrument $k$. More discussion can be found in Sundberg (2003).

There were several short papers that commented on Efron and Hinkley (1978) that appeared in the same issue of the journal containing Efron and Hinkley (1978). For example, Barndorff-Nielsen (1978) discusses ancillarity properties of $\bar{H}_n(\hat{\theta}_n)$ in a more general sense. He stated that part of Efron and Hinkley's (1978) paper perpetuates the impression that $\bar{H}_n(\hat{\theta}_n)$ is, in general, an approximate ancillary statistic (see remarks immediately after formulae (1.5) and (1.6) in Efron and Hinkely (1978)). He pointed out that this impression is not true. He also argues with an example that the possible ancillarity properties of $\bar{H}_n(\hat{\theta}_n)$ depend on the parameterization chosen. An ancillary statistic under one parameterization may not be ancillary if the model is reparameterized. Besides, in Efron and Hinkley (1978), a number of examples are demonstrated in which $\bar{H}_n(\hat{\theta}_n)$ is preferable to $\bar{F}_n(\hat{\theta}_n)$. To argue that this is not always the case, James (1978) deliberately modified an example of Cox (1958) where $\bar{F}_n(\hat{\theta}_n)$ is superior to $\bar{H}_n(\hat{\theta}_n)$ in estimating the variance of an error term. Likewise, Sprott (1978) also provides an example where $\bar{F}_n(\hat{\theta}_n)$ is more accurate than $\bar{H}_n(\hat{\theta}_n)$.

Efron and Hinkley (1978) appeared at the forefront of the wave of interest in conditional inference and asymptotics for parametric models. The paper was motivated by Fisher's (1934) statement that the information loss for MLE in location-scale families can be recovered completely by basing inference on the conditional distribution of the MLE $\hat{\theta}_n$ given an exact ancillary statistic $a$ for which ($\hat{\theta}_n$, $a$) is sufficient (DiCiccio, 2008). A statistic is sufficient if no other statistic that can be calculated from the same sample provides any additional information as to the value of the parameter.

There has been much subsequent work that follows Efron and Hinkley (1978). Most of such work has focused on developing approximate ancillaries, instead of exact ancillaries, and on approximating the conditional variances of the MLE. For instance, Cox (1980) introduced the concept of local ancillary and discussed second-order local ancillaries for scalar-parameter models. Ryall (1981) extended Cox's (1980) result of second-order local ancillary to the vector parameter case and developed a third-order local ancillary for scalar parameter models. Skovgaard (1985) studied general vector parameter models and developed a second-order local ancillary analogous to the ancillarity in Efron and Hinkley (1978). Barndorff-Nielsen (1980) and Amari (1982) discussed various approximate ancillaries in the context of curved exponential families. More subsequent work based on Efron and Hinkley (1978) can be found in Pedersen (1981), Grambsch (1983), McCullagh (1984), and Sweeting (1992), etc.

However, the reliance on an ancillary statistic imposes a major practical limitation in real-life applications. The ancillary statistic is often hard to specify in practice. That is, it is either difficult to define or is not unique in many practical problems. And it is even harder to find a pair ($\hat{\theta}_n$, $a$) that is sufficient. In DiCiccio (2008), a comment on the conclusion of Efron and Hinkley (1978) appears as: "*One obstacle to extending the results for translation families to more general scalar parameter models is that typically no exact ancillary statistic a exists such that (*$\hat{\theta}_n$*, a) is sufficient*". Thus, the conditional variance approach might not be as applicable in practical problems. Besides, theoretical conclusions in Efron and Hinkley (1978) only hold for one-parameter translation families, which is another constraint on general application.

Despite the practical limitations discussed above, the main message of Efron and Hinkley (1978), that the variance estimates for MLE should be constructed from observed information, is still widely accepted; see, for example, McLachlan and Peel (2000), Agresti (2002), and Lawless (2002). However, it has been found in the literature that some papers use the conclusions of Efron and Hinkley (1978) without strictly following the underlying assumptions. For example, in Caceres et. al. (1999), Efron and Hinkley (1978) is cited to validate the use of observed information in variance estimation of MLE for confidence interval construction. But no discussion on ancillary statistics is seen throughout the paper. Similar reference of Efron and Hinkley's result can also be found in Hosking and Wallis (1987), Kass (1987), and Raftery (1996). The presence of such references reflects the fact that the theoretical foundation for a good covariance estimate for MLE is of great interest and value in the literature. However, there is no solid theoretical development in this area yet. This fact further motivates the pursuit of a theoretical analysis for covariance estimation of MLE in this dissertation.

Unlike Efron and Hinkley (1978), Lindsay and Li (1997) avoided the concept of ancillarity. They showed that for $p$-dimensional parameter models, if an error of magnitude $O(n^{-3/2})$ is ignored, $\bar{\boldsymbol{H}}_n(\hat{\boldsymbol{\theta}}_n)^{-1}$ is the optimal estimator of the realized squared error among all asymptotically linear estimators (see Hampel (1974) and Bickel, Klaassen, Ritov, and Wellner (1993, p.19)). That is, for all $i, j = 1, 2, \ldots, p$, $\bar{\boldsymbol{H}}_n(\hat{\boldsymbol{\theta}}_n)^{-1}$ solves the optimization problem:

$$\min_{T(\boldsymbol{X})} E\left[\left(\left(n(\hat{\boldsymbol{\theta}}_n - \boldsymbol{\theta}^*)(\hat{\boldsymbol{\theta}}_n - \boldsymbol{\theta}^*)^T\right)_{i,j} - \left(T(\boldsymbol{X})\right)_{i,j}\right)^2\right], \tag{1.6}$$

where $(\cdot)_{i,j}$ denotes the $(i, j)$th entry of a matrix and $T(\boldsymbol{X})$ is any statistic chosen from a class of asymptotically linear estimators based on the sample data $\boldsymbol{X}$. Here asymptotically linear estimators are defined as linear combinations of functions of each observation plus a term that converges to zero asymptotically. This class of estimators includes $\bar{\boldsymbol{F}}_n(\hat{\boldsymbol{\theta}}_n)^{-1}$, $\bar{\boldsymbol{H}}_n(\hat{\boldsymbol{\theta}}_n)^{-1}$, etc.

The construction of (1.6) indicates that Lindsay and Li's work does not directly estimate the variance of MLE. Instead, the estimation target is the realized squared error rather than the covariance matrix of normalized $\hat{\boldsymbol{\theta}}_n$, where the two differ by an operation of expectation. Specifically, the expectation of the realized squared error is the covariance matrix. Lindsay and Li's work does not directly solve our problem of interest, which is on covariance matrix estimation of MLE. However, the paper has great value in stimulating the approach that follows.

Compared to (1.6), Cao and Spall (2009, 2010) proposed an alternative to determining the best approximation to the variance of $\hat{\theta}_n$ when $\theta$ is a scalar. Specifically, the optimization problem is revised with the adjustment of the estimation target:

$$\min_{T(\boldsymbol{X})} E\left[\left(n\,\text{var}(\hat{\theta}_n) - T(\boldsymbol{X})\right)^2\right], \tag{1.7}$$

where $T(\boldsymbol{X})$ denotes an estimate of the variance of normalized $\hat{\theta}_n$ based on sample data $\boldsymbol{X}$. In Cao and Spall (2009, 2010), $T(\boldsymbol{X})$ is constrained to two candidates: $\bar{F}_n(\hat{\theta}_n)^{-1}$ or $\bar{H}_n(\hat{\theta}_n)^{-1}$. This idea of minimizing the mean squared error of estimation was discussed in Sandved (1968) in the context of approximating a measure of accuracy for a parameter estimate. In Cao and Spall (2009), it is shown that for scalar $\theta$, $\bar{F}_n(\hat{\theta}_n)^{-1}$ is a better estimator of $n\text{var}(\hat{\theta}_n)$ than $\bar{H}_n(\hat{\theta}_n)^{-1}$ under criterion (1.7) with some reasonable conditions. In this paper, we generalize the above scalar result to multivariate $\boldsymbol{\theta}$.

The comparison of $\bar{\boldsymbol{F}}_n(\hat{\boldsymbol{\theta}}_n)^{-1}$ and $\bar{\boldsymbol{H}}_n(\hat{\boldsymbol{\theta}}_n)^{-1}$ has also been done in other aspects. For example, in a score test, the numerator of the test statistic is the squared score function, which is the first derivative of the log-likelihood function with respect to the parameter of interest. The denominator of the test statistic can be either $\bar{\boldsymbol{F}}_n(\hat{\boldsymbol{\theta}}_n)^{-1}$ or $\bar{\boldsymbol{H}}_n(\hat{\boldsymbol{\theta}}_n)^{-1}$ . In practice, $\bar{\boldsymbol{F}}_n(\hat{\boldsymbol{\theta}}_n)^{-1}$ is preferred to $\bar{\boldsymbol{H}}_n(\hat{\boldsymbol{\theta}}_n)^{-1}$ since the latter may result in a negative test statistic; see Morgan *et al* (2007), Verbeke *et al* (2007), and Freedman (2007). The relative merit of $\bar{\boldsymbol{F}}_n(\hat{\boldsymbol{\theta}}_n)^{-1}$ and $\bar{\boldsymbol{H}}_n(\hat{\boldsymbol{\theta}}_n)^{-1}$ is also discussed in the context of iterative calculation of MLE, where Newton's method or scoring method can be used for situations in which closed form of MLE is not attainable; see Fisher (1925), Green (1984), and Garwood (1941). Another area where $\bar{\boldsymbol{F}}_n(\hat{\boldsymbol{\theta}}_n)^{-1}$ and $\bar{\boldsymbol{H}}_n(\hat{\boldsymbol{\theta}}_n)^{-1}$ is compared is the construction of confidence regions, see Royal (1986) and Rust, et. al. (2011).

### 1.3 New approach

In this section, we first lay out the problem settings discussed in this work and then briefly introduce the approach we take to achieve the theoretical conclusion.

To keep our context as general as possible, we consider sequences of i.n.i.d. random vectors, which is often of more practical interest than i.i.d samples, throughout our discussion. The parameter considered is multivariate to accommodate for general practical situations.

The main goal of this work is to compare the performance of $\bar{\boldsymbol{F}}_n(\hat{\boldsymbol{\theta}}_n)^{-1}$ and $\bar{\boldsymbol{H}}_n(\hat{\boldsymbol{\theta}}_n)^{-1}$ in estimating the scaled covariance matrix of MLE, which is denoted by $n\,\text{cov}(\hat{\boldsymbol{\theta}}_n)$. We follow the idea used in Lindsay and Li (1997) and Cao and Spall (2009). We want to solve the following optimization problem:

$$\min_{T(\boldsymbol{X})} E\left[\left(\left(n\operatorname{cov}(\hat{\boldsymbol{\theta}}_n)\right)_{i,j} - \left(T(\boldsymbol{X})\right)_{i,j}\right)^2\right], \tag{1.8}$$

Specifically, our current discussion focuses on $\boldsymbol{T}(\mathbf{X})$ being either $\bar{\boldsymbol{F}}_n(\hat{\boldsymbol{\theta}}_n)^{-1}$ or $\bar{\boldsymbol{H}}_n(\hat{\boldsymbol{\theta}}_n)^{-1}$. Generalization to other estimation candidate $\boldsymbol{T}(\mathbf{X})$ may be considered in future work.

In essence, we compare the performance of $\bar{\boldsymbol{F}}_n(\hat{\boldsymbol{\theta}}_n)^{-1}$ and $\bar{\boldsymbol{H}}_n(\hat{\boldsymbol{\theta}}_n)^{-1}$ at the individual entry level. If we can show that one is better than the other for every matrix entry, then we have found the better of the two in estimating $n\operatorname{cov}(\hat{\boldsymbol{\theta}}_n)$.

## 1.4 Background

Standard results have been established for large sample properties for i.i.d. samples including the central limit theorem (CLT) for the raw data, and the weak law of large numbers (WLLN). In reality, however, observations are frequently not generated from i.i.d samples. In this section, we discuss the CLT, and the WLLN for i.n.i.d. samples. Specifically, we present sufficient conditions that lead to these properties. These conditions will be used in the theoretical development in Chapter 2. All limits below are as $n \to \infty$.

### 1.4.1. *The central limit theorem*

The CLT states that under certain conditions, the distribution of a normalized sample mean of a sequence approaches a normal distribution, i.e.,

$$n^{-1/2}\sum_{i=1}^{n}\left(\eta_i - E(\eta_i)\right)\xrightarrow{\text{dist}} N\left(0, n^{-1}\sum_{i=1}^{n}\sigma_i^2\right),$$

where $\{\eta_1, \eta_2, \ldots, \eta_n\}$ is a sequence of i.n.i.d random variables with corresponding variances $\{\sigma_1^2, \sigma_2^2, \ldots, \sigma_n^2\}$. Various studies of conditions under which the above asymptotic distribution holds have been made by Chebyshev (1980), Feller (1935), Levy (1935), Lindberg (1922), Lyapunov (1900, 1901), Markov (1900), and others.

For a random sample $\{\eta_1, \eta_2, \ldots, \eta_n\}$, the following well-known Lindberg-Feller condition guarantees the CLT result:

**A1**. $\eta_1, \eta_2, \ldots, \eta_n$ is a sequence of independent samples;

**A2**. For every $\varepsilon > 0$, $\lim_{n\to\infty} \sum_{i=1}^{n} E\left(\left(\eta_i - E(\eta_i)\right)^2 \cdot \mathbf{1}_{\{|\eta_i - E(\eta_i)| > \varepsilon s_n\}}\right) \Big/ s_n^2 = 0$, where $s_n^2 = \sum_{i=1}^{n}\sigma_i^2$ and $\mathbf{1}_{\{\ldots\}}$ is the indicator function.

1.4.2. *Weak law of large numbers*

The WLLN states that under certain conditions, the sample mean of a sequence converges in probability to the average population mean:

$$\frac{1}{n}\sum_{i=1}^{n}\left(\eta_i - E(\eta_i)\right)\xrightarrow{p} 0,$$

where $\{\eta_1, \eta_2, \ldots, \eta_n\}$ is a sequence of independent random variables and $\xrightarrow{p}$ denotes convergence in probability.

A set of sufficient conditions for WLLN for i.n.i.d samples is presented in Chung (2005, Theorem 5.2.3):

**B.1**. $\eta_1, \eta_2, \ldots, \eta_n$ is a sequence of independent samples;

**B.2**. $\sum_{i=1}^{n} E\left(\eta_i 1_{\{|\eta_i|>n\}}\right) \to 0$ where $1_{\{\mathrm{A}\}}$ is an indicator function which equals 1 if the condition denoted by A holds and 0 otherwise;

**B.3**. $n^{-2}\sum_{i=1}^{n} E\left(\eta_i^2 1_{\{|\eta_i|\le n\}}\right) \to 0$.

In this dissertation, we apply conditions presented in Sections 1.4.1 and 1.4.2 for the CLT and the WLLN under i.n.i.d samples. We discuss more on the concrete forms of i.n.i.d sequences in Chapter 2.

This dissertation is organized as follows. In Chapter 2, we present the theoretical development that leads to the main result. In Chapter 3, we present numerical studies on three distinct problems to support the main theoretical result. In Chapter 4, we summarize the achievement in this dissertation and discuss potential future work to extend the results of this dissertation. This dissertation also includes two appendices on topics of relevance to stochastic systems. In Appendix A, we discuss optimal perturbation distributions for the simultaneous perturbation stochastic approximation (SPSA) algorithm. In Appendix B, we consider Monte Carlo methods for computing FIMs when closed forms are not attainable.

# Chapter 2

# Theoretical Analysis

In this chapter, we present the theoretical development in this dissertation on comparing the expected and observed FIM in estimating the covariance matrix of MLEs. In Section 2.1, we begin with notation definitions, followed by a discussion on a list of sufficient conditions used to achieve the theoretical conclusion in Section 2.2. In Section 2.3, we present preliminary results as a preparation for the main result. In Section 2.4, we present the main result.

## 2.1 Notation

As defined in Chapter 1, $\boldsymbol{X} = [\boldsymbol{X}_1, \boldsymbol{X}_2, \ldots, \boldsymbol{X}_n]$ is a collection of i.n.i.d. random vectors (variables) where $\boldsymbol{X}_i \in \mathbb{R}^q$, $i = 1, 2, \ldots, n$, and $q \geq 1$. Each $\boldsymbol{X}_i$ may contain discrete or continuous components. If we let $\boldsymbol{X}_i^d$ and $\boldsymbol{X}_i^c$ denote the sub-vectors of discrete and continuous components of $\boldsymbol{X}_i$, respectively, then $\dim(\boldsymbol{X}_i^d) + \dim(\boldsymbol{X}_i^c) = q$, where $\dim(\cdot)$ denotes the dimension of a vector. Either $\boldsymbol{X}_i^d$ or $\boldsymbol{X}_i^c$ may be a null sub-vector for a given $\boldsymbol{X}_i$, i.e., dim( $\boldsymbol{X}_i^d$ ) = 0 or dim( $\boldsymbol{X}_i^c$ ) = 0. And dim( $\boldsymbol{X}_i^d$ ) = 0 implies that all elements in $\boldsymbol{X}_i$ are continuous and vice versa.

Recalling the definitions in Chapter 1, the probability density/mass function and the negative log-likelihood function of $\boldsymbol{X}_i$ are $p_i(\boldsymbol{x}_i, \boldsymbol{\theta})$ and $l_i(\boldsymbol{\theta}, \boldsymbol{x}_i) \equiv -\log p_i(\boldsymbol{x}_i, \boldsymbol{\theta})$, respectively, where $\boldsymbol{\theta} = [t_1, t_2, \ldots, t_p]^T \in \Theta$ is a $p$-dimensional vector valued parameter. The joint density/mass function and the negative log-likelihood function of $\boldsymbol{X}$ are $p(\boldsymbol{x}, \boldsymbol{\theta}) \equiv \prod_{i=1}^n p_i(\boldsymbol{x}_i, \boldsymbol{\theta})$ and $l(\boldsymbol{\theta}, \boldsymbol{x}) \equiv \sum_{i=1}^n l_i(\boldsymbol{\theta}, \boldsymbol{x}_i) = -\sum_{i=1}^n \log p_i(\boldsymbol{x}_i, \boldsymbol{\theta})$, respectively. The MLE for $\boldsymbol{\theta}$ based on $\boldsymbol{X}$ is denoted as $\hat{\boldsymbol{\theta}}_n = [\hat{t}_{n1}, \hat{t}_{n2}, \ldots, \hat{t}_{np}]^T$ and the true value of $\boldsymbol{\theta}$ is $\boldsymbol{\theta}^* = [t_1^*, t_2^*, \ldots, t_p^*]^T$. Let $U_r^i$, $U_{rs}^i$, and $U_{rst}^i$ be the derivatives of $l_i(\boldsymbol{\theta}, \boldsymbol{x}_i)$ according to $U_r^i \equiv \partial l_i(\boldsymbol{\theta}, \boldsymbol{x}_i)/\partial t_r$, $U_{rs}^i \equiv \partial^2 l_i(\boldsymbol{\theta}, \boldsymbol{x}_i)/\partial t_r \partial t_s$, and $U_{rst}^i \equiv \partial^3 l_i(\boldsymbol{\theta}, \boldsymbol{x}_i)/\partial t_r \partial t_s \partial t_t$. Correspondingly, $U_r$, $U_{rs}$, and $U_{rst}$ are the derivatives of $l(\boldsymbol{\theta}, \boldsymbol{x})$ according to $U_r \equiv \partial l(\boldsymbol{\theta}, \boldsymbol{x})/\partial t_r$, $U_{rs} \equiv \partial^2 l(\boldsymbol{\theta}, \boldsymbol{x})/\partial t_r \partial t_s$, and $U_{rst} \equiv \partial^3 l(\boldsymbol{\theta}, \boldsymbol{x})/\partial t_r \partial t_s \partial t_t$. Note that $U_{rs}$ is the $(r, s)$ entry of $\boldsymbol{H}_n(\boldsymbol{\theta})$.

Let us define null-cumulants for each observation $\boldsymbol{X}_i$ (e.g. $\kappa_r^i, \kappa_{r,s}^i$, etc.) and average null-cumulants per observation (e.g. $\bar{\kappa}_r, \bar{\kappa}_{r,s}$, etc.) as follows: (All expectations are well defined and the word “null” refers to the fact that the twin processes of differentiation and averaging both take place at the same value: $\boldsymbol{\theta}^*$, see McCullagh (1987, page 201)):

$$\kappa_r^i \equiv E(U_r^i), \tag{2.1a}$$

$$\kappa_{rs}^{i} \equiv E(U_{rs}^{i}), \tag{2.1b}$$

$$\kappa_{rst}^{i} \equiv E(U_{rst}^{i}), \tag{2.1c}$$

$$\kappa_{r,s}^{i} \equiv E(U_r^i U_s^i) - E(U_r^i)E(U_s^i), \tag{2.1d}$$

$$\kappa_{rs,t}^{i} \equiv E(U_{rs}^i U_t^i) - E(U_{rs}^i)E(U_t^i), \tag{2.1e}$$

$$\bar{\kappa}_r \equiv \sum_{i=1}^{n} \kappa_r^i \Big/ n, \tag{2.2a}$$

$$\bar{\kappa}_{rs} \equiv \sum_{i=1}^{n} \kappa_{rs}^i \Big/ n, \tag{2.2b}$$

$$\bar{\kappa}_{rst} \equiv \sum_{i=1}^{n} \kappa_{rst}^i \Big/ n, \tag{2.2c}$$

$$\bar{\kappa}_{r,s} \equiv \sum_{i=1}^{n} \kappa_{r,s}^i \Big/ n, \tag{2.2d}$$

$$\bar{\kappa}_{rs,t} \equiv \sum_{i=1}^{n} \kappa_{rs,t}^i \Big/ n. \tag{2.2e}$$

The standardized likelihood scores, denoted by indexed $Z$'s, are the derivatives of the negative log-likelihood centered by its expectation and scaled by $n^{-1/2}$. That is,

$$Z_r \equiv n^{-1/2}\left(U_r(\boldsymbol{\theta}^*) - n\bar{\kappa}_r\right), \tag{2.3}$$

and

$$Z_{st} \equiv n^{-1/2}\left(U_{st}(\boldsymbol{\theta}^*) - n\bar{\kappa}_{st}\right). \tag{2.4}$$

We assume that the likelihood function is regular in the sense that necessary interchanges of differentiation and integration are valid (more details are provided in Section 2.2 below). Furthermore, given the notation of $\boldsymbol{X}_i^d$ and $\boldsymbol{X}_i^c$, $p_i(\boldsymbol{x}_i, \boldsymbol{\theta})$ can be decomposed as a product of two terms: $p_i(\boldsymbol{x}_i, \boldsymbol{\theta}) = p_i^c(\boldsymbol{x}_i^c, \boldsymbol{\theta} | \boldsymbol{x}_i^d) \times p_i^d(\boldsymbol{x}_i^d, \boldsymbol{\theta})$, where $p_i^c(\boldsymbol{x}_i^c, \boldsymbol{\theta} | \boldsymbol{x}_i^d)$ is the conditional density function of $\boldsymbol{X}_i^c$ given $\boldsymbol{X}_i^d = \boldsymbol{x}_i^d$, and $p_i^d(\boldsymbol{x}_i^d, \boldsymbol{\theta})$ is the marginal mass function of $\boldsymbol{X}_i^d$. Let $\boldsymbol{S}_i^d$ denote the support of $\boldsymbol{X}_i^d$ and $\boldsymbol{S}_i^c | \boldsymbol{x}_i^d$ denote the support of $\boldsymbol{X}_i^c$ given $\boldsymbol{X}_i^d = \boldsymbol{x}_i^d$. Now we are ready to show that with valid interchange of differentiation and integration, $E(U_r^i) = 0$, for $i$ = 1, 2, …, $n$ and $r$ = 1, 2, …, $p$. In fact,

$$E(U_r^i) = E\left(\frac{\partial l_i(\boldsymbol{\theta}, \boldsymbol{x}_i)}{\partial t_r}\right)$$

$$= \sum_{\boldsymbol{x}_i^d \in \boldsymbol{S}_i^d} \int_{\boldsymbol{x}_i^c \in \boldsymbol{S}_i^c | \boldsymbol{x}_i^d} \frac{\partial l_i(\boldsymbol{\theta}, \boldsymbol{x}_i)}{\partial t_r} p_i^c(\boldsymbol{x}_i^c, \boldsymbol{\theta} | \boldsymbol{x}_i^d) d\boldsymbol{x}_i^c \times p_i^d(\boldsymbol{x}_i^d, \boldsymbol{\theta})$$

$$= -\sum_{\boldsymbol{x}_i^d \in \boldsymbol{S}_i^d} \int_{\boldsymbol{x}_i^c \in \boldsymbol{S}_i^c | \boldsymbol{x}_i^d} \frac{1}{p_i(\boldsymbol{x}_i, \boldsymbol{\theta})} \frac{\partial p_i(\boldsymbol{x}_i, \boldsymbol{\theta})}{\partial t_r} p_i^c(\boldsymbol{x}_i^c, \boldsymbol{\theta} | \boldsymbol{x}_i^d) p_i^d(\boldsymbol{x}_i^d, \boldsymbol{\theta}) d\boldsymbol{x}_i^c$$

$$= -\sum_{\boldsymbol{x}_i^d \in \boldsymbol{S}_i^d} \int_{\boldsymbol{x}_i^c \in \boldsymbol{S}_i^c | \boldsymbol{x}_i^d} \frac{\partial p_i(\boldsymbol{x}_i, \boldsymbol{\theta})}{\partial t_r} d\boldsymbol{x}_i^c \quad \text{(cancellation of } p_i(\boldsymbol{x}_i, \boldsymbol{\theta})\text{)}$$

$$= -\frac{\partial}{\partial t_r} \sum_{\boldsymbol{x}_i^d \in \boldsymbol{S}_i^d} \int_{\boldsymbol{x}_i^c \in \boldsymbol{S}_i^c | \boldsymbol{x}_i^d} p_i(\boldsymbol{x}_i, \boldsymbol{\theta}) d\boldsymbol{x}_i^c$$

(interchange of differentiation and integration)

$$= -\frac{\partial}{\partial t_r} \sum_{\boldsymbol{x}_i^d \in \boldsymbol{S}_i^d} \int_{\boldsymbol{x}_i^c \in \boldsymbol{S}_i^c | \boldsymbol{x}_i^d} p_i^c(\boldsymbol{x}_i^c, \boldsymbol{\theta} | \boldsymbol{x}_i^d) d\boldsymbol{x}_i^c \times p_i^d(\boldsymbol{x}_i^d, \boldsymbol{\theta})$$

$$= \frac{\partial(-1)}{\partial t_r} \quad \text{(mass/density function integrates to 1)}$$

$$= 0.$$

Thus, $\kappa_r^i = 0$ and $\bar{\kappa}_r = 0$ for all $i$ and $r$; $Z_r = n^{-1/2} U_r$ for all $r$ according to the definition in (2.3).

Let $\bar{\kappa}^{v,u}$ be the ($v$, $u$) element of the inverse matrix of $\bar{\boldsymbol{\kappa}}$, where $\bar{\boldsymbol{\kappa}}$ is a $p \times p$ matrix whose ($s$, $t$) element is $\bar{\kappa}_{s,t}$, $s$, $t$ = 1, …, $p$. Throughout this paper, the double bar notation ($\overline{\overline{\cdot}}$) indicates a special summation operation. Specifically, for the argument under the double bar, summation is implied over any index repeated once as a superscript and once as a subscript. For example,

$$\overline{\overline{\bar{\kappa}_{st,v} \bar{\kappa}^{v,u} Z_u}} = \sum_{v=1}^{p} \sum_{u=1}^{p} \bar{\kappa}_{st,v} \bar{\kappa}^{v,u} Z_u \ ;$$

$$\overline{\overline{\bar{\kappa}^{r,t}\bar{\kappa}^{s,u}(U_{tu}^{i}-\bar{\kappa}_{tu}^{i}-\bar{\kappa}_{tu,v}\bar{\kappa}^{v,w}U_{w}^{i})}} = \sum_{t=1}^{p}\sum_{u=1}^{p}\bar{\kappa}^{r,t}\bar{\kappa}^{s,u}\left(U_{tu}^{i}-\bar{\kappa}_{tu}^{i}-\sum_{v=1}^{p}\sum_{w=1}^{p}\bar{\kappa}_{tu,v}\bar{\kappa}^{v,w}U_{w}^{i}\right).$$

This short-hand notation of summation is the same as the index notation used in McCullagh (1987) and Lindsay and Li (1997) except that we add the double bar notation to distinguish the summation from each individual summand.

To orthogonalize $Z_r$ and $Z_{st}$, we define

$$Y_{st} \equiv Z_{st} - \overline{\overline{\bar{\kappa}_{st,v}\bar{\kappa}^{v,u}Z_{u}}}.$$

Given the definition above, we have $\text{cov}(Z_r,\ Y_{st}) = 0$, $r$, $s$, $t = 1, \ldots, p$, which is an important property used in Sections 2.3 and 2.4. The uncorrelatedness is seen by noting:

$$\text{cov}(Z_r,\ Y_{st}) = \text{cov}(Z_r,\ Z_{st} - \overline{\overline{\bar{\kappa}_{st,v}\bar{\kappa}^{v,u}Z_{u}}}\ )$$

$$= \text{cov}\left(n^{-1/2}\sum_{i=1}^{n}U_{r}^{i}, n^{-1/2}\sum_{i=1}^{n}\left(U_{st}^{i}-\bar{\kappa}_{st}\right)-n^{-1/2}\overline{\overline{\bar{\kappa}_{st,v}\bar{\kappa}^{v,u}\sum_{i=1}^{n}U_{u}^{i}}}\right)$$

(definitions (2.3) and (2.4))

$$= \text{cov}\left(n^{-1/2}\sum_{i=1}^{n}U_{r}^{i}, n^{-1/2}\sum_{i=1}^{n}\left(U_{st}^{i}-\bar{\kappa}_{st}\right)\right)$$

$$-\text{cov}\left(n^{-1/2}\sum_{i=1}^{n}U_{r}^{i}, n^{-1/2}\overline{\overline{\bar{\kappa}_{st,v}\bar{\kappa}^{v,u}\sum_{i=1}^{n}U_{u}^{i}}}\right)$$

$$= n^{-1}\operatorname{cov}\left(\sum_{i=1}^{n} U_r^i, \sum_{i=1}^{n} U_{st}^i\right) - n^{-1}\overline{\overline{\bar{\kappa}_{st,v}\bar{\kappa}^{v,u}\operatorname{cov}\left(\sum_{i=1}^{n} U_u^i, \sum_{i=1}^{n} U_r^i\right)}}$$

$$= n^{-1}\sum_{i=1}^{n}\operatorname{cov}(U_r^i, U_{st}^i) - n^{-1}\overline{\overline{\bar{\kappa}_{st,v}\bar{\kappa}^{v,u}\sum_{i=1}^{n}\operatorname{cov}(U_u^i, U_r^i)}}$$

(independence between observations)

$$= n^{-1}\sum_{i=1}^{n}\kappa_{st,r}^{i} - n^{-1}\overline{\overline{\bar{\kappa}_{st,v}\bar{\kappa}^{v,u}\sum_{i=1}^{n}\kappa_{u,r}^{i}}} \quad \text{(definitions (3.1d) and (3.1e) )}$$

$$= \bar{\kappa}_{st,r} - \overline{\overline{\bar{\kappa}_{st,v}\bar{\kappa}^{v,u}\bar{\kappa}_{u,r}}} \quad \text{(definitions (3.2d) and (3.2e))}$$

$$= \bar{\kappa}_{st,r} - \bar{\kappa}_{st,r}$$

$$= 0.$$

In the discussion below, we frequently use the stochastic big-$O$ and little-$o$ terms: $O_d(n^{-1})$, $O_d^2(n^{-1})$, $o_p(n^{-1})$, and $o_p(1)$. Specifically, $O_d(n^{-1})$ denotes a stochastic term that converges in distribution to a random variable when multiplied by $n$; $O_d^2(n^{-1})$ denotes a product of two $O_d(n^{-1})$ terms; $o_p(n^{-1})$ is a stochastic term that converges in probability to zero when multiplied by $n$; and $o_p(1)$ is a stochastic term that converges in probability to zero, i.e., $o_p(1) = n \times o_p(n^{-1})$. In addition, for simplicity, we introduce $\tilde{O}_d(n^{-1})$ to denote a summation of a finite number of $O_d(n^{-1})$ terms and $\tilde{O}_d^2(n^{-1})$ to denote a summation of a finite number of $O_d^2(n^{-1})$ terms.

## 2.2 Conditions

In this section, we introduce sufficient conditions for the analytical development below. We provide some interpretation of the conditions immediately following the presentation of the conditions below.

**A1**. Necessary interchanges of differentiation and integration are valid for the following functions denoted generally as $g(\boldsymbol{x}_i, \boldsymbol{\theta})$:

I. $p_i(\boldsymbol{x}_i, \boldsymbol{\theta})$, $i = 1, 2, \ldots, n$;

II. $U_{rs}^{i} \exp\{-\, l_i(\boldsymbol{x}_i, \boldsymbol{\theta})\}$, $i = 1, 2, \ldots, n$ and $r, s = 1, 2, \ldots, p$.

Specifically, the following conditions hold for $g(\boldsymbol{x}_i, \boldsymbol{\theta})$:

A1(a). $g(\boldsymbol{x}_i, \boldsymbol{\theta})$ and $\partial g(\boldsymbol{x}_i, \boldsymbol{\theta})/\partial t_j$ are continuous on $\Theta \times \mathbb{R}^q$ for $j = 1, 2, \ldots, p$;

A1(b). There exist nonnegative functions $q_0(\boldsymbol{x}_i)$ and $q_1(\boldsymbol{x}_i)$ such that

$$| g(\boldsymbol{x}_i, \boldsymbol{\theta})| \le q_0(\boldsymbol{x}_i),\ |\partial g(\boldsymbol{x}_i, \boldsymbol{\theta})/\partial t_j| \le q_1(\boldsymbol{x}_i) \text{ for all } \boldsymbol{x}_i \in \mathbb{R}^q \text{ and } \boldsymbol{\theta} \in \Theta,$$

where $\sum_{\boldsymbol{x}_i^d \in \boldsymbol{S}_i^d} \int_{\boldsymbol{x}_i^c \in \boldsymbol{S}_i^c | \boldsymbol{x}_i^d} q_0(\boldsymbol{x}_i) d\boldsymbol{x}_i^c < \infty$ and $\sum_{\boldsymbol{x}_i^d \in \boldsymbol{S}_i^d} \int_{\boldsymbol{x}_i^c \in \boldsymbol{S}_i^c | \boldsymbol{x}_i^d} q_1(\boldsymbol{x}_i) d\boldsymbol{x}_i^c < \infty$.

**A2**. The negative log-likelihood function $l(\boldsymbol{x}, \boldsymbol{\theta})$ has continuous partial derivatives with respect to $\boldsymbol{\theta}$ up to the fourth order and all expectations in (2.1a–e) are well defined.

**A3**. $\boldsymbol{F}_n(\boldsymbol{\theta}^*)$ is positive definite, $\boldsymbol{F}(\boldsymbol{\theta}^*) \equiv \lim_{n\to\infty} \bar{\boldsymbol{F}}_n(\boldsymbol{\theta}^*)$ exists and is invertible.

**A4**. The following limits exist and are finite in magnitude:

A4(a): $\lim_{n\to\infty} \bar{\kappa}_{rst} = \lim_{n\to\infty} n^{-1}\sum_{i=1}^{n} \kappa_{rst}^{i}$ for $r, s, t = 1, 2, \ldots, p$;

A4(b): $\lim_{n\to\infty}\ n^{-1}\partial E(U_{rs})\big/\partial t_i\,|_{\boldsymbol{\theta}\,=\,\boldsymbol{\theta}^*}$ for $r, s, i = 1, 2, \ldots, p$;

A4(c): $\lim_{n\to\infty}\ n^{-1}\partial^2 E(U_{rs})\big/\partial t_i\,\partial t_j\,|_{\boldsymbol{\theta}\,=\,\boldsymbol{\theta}^*}$ for $r, s, i, j = 1, 2, \ldots, p$;

A4(d): $\lim_{n\to\infty} n^{-1}\, E(U_{rstv})|_{\boldsymbol{\theta}\,=\,\boldsymbol{\theta}^*}$ for $r, s, t, v = 1, 2, \ldots, p$.

**A5**. The Lindberg-Feller condition holds for the following independent sequences denoted generally as $\xi_1,\ldots,\xi_n$:

I. $\{U_r^i(\boldsymbol{\theta}^*)\}_{i=1}^n$ for $r = 1, 2, \ldots, p$.

II. $\{U_{rs}^i(\boldsymbol{\theta}^*)\}_{i=1}^n$ for $r, s = 1, 2, \ldots, p$.

III. $\{U_{rst}^i(\boldsymbol{\theta}^*)\}_{i=1}^n$ for $r, s, t = 1, 2, \ldots, p$.

Specifically, $\lim_{n\to\infty}\sum_{i=1}^n E\left(\left(\xi_i - E\left(\xi_i\right)\right)^2\cdot\mathbf{1}_{\left\{\left|\xi_i - E\left(\xi_i\right)\right|>\varepsilon s_n\right\}}\right)\Big/ s_n^2 = 0$ for every $\varepsilon > 0$,

where $s_n^2 = \sum_{i=1}^n \sigma_i^2$, $\sigma_i^2$ is the variance of $\xi_i$, and $\mathbf{1}_{\{\ldots\}}$ is the indicator function.

**A6**. Conditions for the WLLN hold for the following independent sequences denoted generally as $\xi_1,\ldots,\xi_n$:

I. $\{U_r^i(\boldsymbol{\theta})\}_{i=1}^n$ for $r = 1, 2, \ldots, p$ and $\boldsymbol{\theta}$ in a neighborhood of $\boldsymbol{\theta}^*$;

II. $\{U_{rs}^i(\boldsymbol{\theta})\}_{i=1}^n$ for $r, s = 1, 2, \ldots, p$ and $\boldsymbol{\theta}$ in a neighborhood of $\boldsymbol{\theta}^*$;

III. $\{U_{rst}^i(\boldsymbol{\theta})\}_{i=1}^n$ for $r, s, t = 1, 2, \ldots, p$ and $\boldsymbol{\theta}$ in a neighborhood of $\boldsymbol{\theta}^*$;

IV. $\{U_{rstv}^i(\boldsymbol{\theta})\}_{i=1}^n$ for $r, s, t, v = 1, 2, \ldots, p$ and $\boldsymbol{\theta}$ in a neighborhood of $\boldsymbol{\theta}^*$;

V. $\overline{\overline{\bar{\kappa}^{r,t}\bar{\kappa}^{s,u}(U_{tu}^i - \bar{\kappa}_{tu,v}\bar{\kappa}^{v,w}U_w^i)}}$ for $i = 1, 2, \ldots, n$ and $r, s, t, v = 1, 2, \ldots, p$;

VI. $\overline{\overline{\overline{\kappa}^{r,t}\overline{\kappa}^{s,u}(\overline{\kappa}_{tuv}-\overline{\kappa}_{tu,v})\overline{\kappa}^{v,w}U_w^i}}$ for $i = 1, 2, \ldots, n$ and $r, s, t, v = 1, 2, \ldots, p$.

Specifically, the following holds for the i.n.i.d sequence $\xi_1,...,\xi_n$:

A6(a). $\sum_{i=1}^{n} E\left(\xi_i 1_{\{|\xi_i|>n\}}\right) \to 0$;

A6(b). $n^{-2}\sum_{i=1}^{n} E\left(\xi_i^2 1_{\{|\xi_i|\le n\}}\right) \to 0$.

**A7**. The dominated convergence theorem (DCT) applies to all stochastic high order terms $o_p(1)$. As a result, for any stochastic term that converges in probability to zero, the corresponding expectation converges to zero as well. Specifically, all $o_p(1)$ terms throughout this paper are formed as a linear combination of a finite number of the following terms and each coefficient converges in probability to a constant:

I. $n^{-3/2}Z_rZ_vZ_{st}$ for $r, v, s, t = 1, 2, \ldots, p$;

II. $n^{-2}Z_rZ_vZ_{uw}Z_{st}$ for $r, v, s, t, u, w = 1, 2, \ldots, p$;

III. $n^{-1/2}Z_r(t_{ni}-t_i^*)\times(t_{nj}-t_j^*)$ for $r, i, j = 1, 2, \ldots, p$;

IV. $n^{-1}Z_rZ_{st}(t_{ni}-t_i^*)\times(t_{nj}-t_j^*)$ for $r, s, t, i, j = 1, 2, \ldots, p$;

V. $(t_{ni}-t_i^*)\times(t_{nj}-t_j^*)\times(t_{nk}-t_k^*)\times(t_{ng}-t_g^*)$ for $i, j, k, g = 1, 2, \ldots, p$;

VI. $\overline{\overline{\overline{\kappa}^{r,t}\overline{\kappa}^{s,u}(U_{tu}^i-\overline{\kappa}_{tu,v}\overline{\kappa}^{v,w}U_w^i)}}$ for $i = 1, 2, \ldots, n$ and $r, s, t, v = 1, 2, \ldots, p$;

VII. $\overline{\overline{\overline{\kappa}^{r,t}\overline{\kappa}^{s,u}(\overline{\kappa}_{tuv}-\overline{\kappa}_{tu,v})\overline{\kappa}^{v,w}U_w^i}}$ for $i = 1, 2, \ldots, n$ and $r, s, t, v = 1, 2, \ldots, p$.

**A8.** The null-cumulants defined in (2.2a–e) are bounded in magnitude for all $n$, i.e. $\limsup_{n\to\infty}|\overline{\kappa}| < \infty$, where $\overline{\kappa}$ represents $\overline{\kappa}_r$, $\overline{\kappa}_{rs}$, $\overline{\kappa}_{r,s}$, $\overline{\kappa}_{rs,t}$ for $r, s, t = 1, 2, \ldots, p$.

**A9.** There exist entries $(r,s)$ such that there is a subsequence $\{n_1, n_2, n_3, \ldots\}$ of $\{1, 2, 3, \ldots\}$ so that $\bar{\boldsymbol{F}}_n(\hat{\boldsymbol{\theta}}_n)^{-1}(r,s)$ and $\bar{\boldsymbol{H}}_n(\hat{\boldsymbol{\theta}}_n)^{-1}(r,s)$ differ for all $n = n_1, n_2, n_3, \ldots$. And for all such entries $(r,s)$, $\liminf_{n\to\infty} \left\{ n^{-1}\sum_{i=1}^{n} \mathrm{var}\left[ \overline{\bar{\kappa}^{r,t}\bar{\kappa}^{s,u}\left(U_{tu}^{i} - \kappa_{tu}^{i} - \bar{\kappa}_{tu,v}\bar{\kappa}^{v,w}U_{w}^{i}\right)} \right] \right\} > 0$.

Condition A1 ensures valid interchange of differentiation and integral on relevant functions, which is crucial in proving $E(U_r^i) = 0$ for $i = 1, 2, \ldots, n$, $r = 1, 2, \ldots, p$ and an intermediate result in Lemma 2 below. Sufficient conditions for interchange of differentiation and integral on likelihood functions are also discussed in Wilks (1962, pp. 408–411 and 418–419) and Bickel and Doksum (2007, p.179). Condition A2 is to guarantee that all null-cumulants are well defined in (2.1a–e). Condition A3 guarantees the limit of the Fisher information exists and is invertible. Limits in Condition A4 are to ensure necessary convergence in the proof of the lemmas below. Specifically, we assume finite limits for average null-cumulants and its derivatives with respect to components of the parameter. Condition A5 describes Lindberg-Feller condition of the CLT for i.n.i.d. samples. In our context, we assume that the CLT holds for sequences of 1st, 2nd, and 3rd derivatives of the log-likelihood function with respect to elements of the parameter. Note that we keep the analysis at individual element level, so we require CLT conditions only for scalar sequences, even though we consider multiple-dimension parameters in our context. Condition A6 presents sufficient conditions for the WLLN for i.n.i.d samples (Chung 2005, Theorem 5.2.3). Condition A6(a) implies that the relevant i.n.i.d samples should not have heavy tails; condition A6(b) indicates that the variance of the sequence cannot grow too fast. Condition A7 assumes that the DCT applies to relevant sequences, which guarantees that the rate of convergence in stochastic sense is preserved after expectation. This condition is implicitly used in Lindsay and Li (1997) and McCullagh (1987, Chapter 7). In condition A8, the imposed boundedness on null-cumulants is to guarantee that any stochastic term multiplied by these cumulants preserve the convergence rate. Condition A9 states that for any entry where $\bar{\boldsymbol{F}}_n(\hat{\boldsymbol{\theta}}_n)^{-1}$ and $\bar{\boldsymbol{H}}_n(\hat{\boldsymbol{\theta}}_n)^{-1}$ differ for a subsequence, the limit inferior of the variance term in the condition is positive. This condition is used to show the superiority of $\bar{\boldsymbol{F}}_n(\hat{\boldsymbol{\theta}}_n)^{-1}$ over $\bar{\boldsymbol{H}}_n(\hat{\boldsymbol{\theta}}_n)^{-1}$ in the main theorem. In fact, the term inside the variance function in condition A9 is random only through the 1st and 2nd derivatives of the log-likelihood function. The condition requires a certain level of variability for the 1st and 2nd derivatives of the log-likelihood function. This is not surprising because if the variability is too low, $\bar{\boldsymbol{F}}_n(\hat{\boldsymbol{\theta}}_n)^{-1}$ and $\bar{\boldsymbol{H}}_n(\hat{\boldsymbol{\theta}}_n)^{-1}$ are very close to each other or even identical. Besides, the concept of "subsequence" in condition A9 allows for the flexibility where $\bar{\boldsymbol{F}}_n(\hat{\boldsymbol{\theta}}_n)^{-1}$ and $\bar{\boldsymbol{H}}_n(\hat{\boldsymbol{\theta}}_n)^{-1}$ do not have to be different for every single term of the sequence. In fact, we only require that $\bar{\boldsymbol{F}}_n(\hat{\boldsymbol{\theta}}_n)^{-1}$ and $\bar{\boldsymbol{H}}_n(\hat{\boldsymbol{\theta}}_n)^{-1}$ be different for infinite terms. Obviously, condition A9 is applicable to situations where $\bar{\boldsymbol{F}}_n(\hat{\boldsymbol{\theta}}_n)^{-1}$ and $\bar{\boldsymbol{H}}_n(\hat{\boldsymbol{\theta}}_n)^{-1}$ differ for all $n$. All the above assumptions are assumed in this dissertation as sufficient conditions for the main result. As discussed above, these are reasonable assumptions that hold for a wide class of problems like other standard conditions.

### 2.3 Preliminary results

Before we present the main result, let us summarize some preliminary results that are essential to our analysis.

**Lemma 1**

For i.n.i.d sample data with conditions A1–A8 in Section 2.2, the estimation error of $\hat{\boldsymbol{\theta}}_n$ has the following form:

$$\hat{t}_{ni} - t_i^* = -n^{-1/2}\overline{\overline{\kappa^{i,u}} Z_u} + n\tilde{O}_d^2(n^{-1}), \tag{2.5}$$

for $r = 1, 2, \ldots, p$.

*Proof*: For $r = 1, 2, \ldots, p$, the MLE $\hat{\boldsymbol{\theta}}_n$ satisfies the equation $n^{-1}U_r(\hat{\boldsymbol{\theta}}_n) = 0$, which can be expanded in a Taylor's series around $\boldsymbol{\theta}^*$ as follows:

$$0 = n^{-1}U_r(\boldsymbol{\theta}^*) + n^{-1}\sum_{i=1}^{p} U_{ri}(\boldsymbol{\theta}^*) \times (\hat{t}_{ni} - t_i^*) + (2n)^{-1}\sum_{i=1}^{p}\sum_{j=1}^{p} U_{rij}(\tilde{\boldsymbol{\theta}}_n) \times (\hat{t}_{ni} - t_i^*) \times (\hat{t}_{nj} - t_j^*), \tag{2.6}$$

where $\tilde{\boldsymbol{\theta}}_n$ is an intermediate point between $\hat{\boldsymbol{\theta}}_n$ and $\boldsymbol{\theta}^*$. Let us write the error term in the following form

$$\hat{t}_{ni} - t_i^* = -n^{-1/2}\overline{\overline{\kappa^{i,u}} Z_u} + R_{ni}, \tag{2.7}$$

where $R_{ni}$ needs to be determined. In order to show that (2.5) is true, we now show that $R_{ni} = n\tilde{O}_d^2(n^{-1})$, for $i$ = 1, 2, …, $p$. Given (2.7), we show that (2.6) can be rewritten as follows:

$$0 = n^{-1}U_r(\boldsymbol{\theta}^*) + n^{-1}\sum_{i=1}^{p} U_{ri}(\boldsymbol{\theta}^*) \times (-n^{-1/2}\overline{\overline{\bar{\kappa}^{i,u}}}Z_u + R_i)$$

$$+(2n)^{-1}\sum_{i=1}^{p}\sum_{j=1}^{p} U_{rij}(\tilde{\boldsymbol{\theta}}_n) \times (\hat{t}_{ni} - t_i^*) \times (\hat{t}_{nj} - t_j^*)$$

$$= n^{-1/2}Z_r + n^{-1}\sum_{i=1}^{p}(n^{1/2}Z_{ri} + n\bar{\kappa}_{ri}) \times (-n^{-1/2}\overline{\overline{\bar{\kappa}^{i,u}Z_u}}) + n^{-1}\sum_{i=1}^{p}(n^{1/2}Z_{ri} + n\bar{\kappa}_{ri})R_{ni}$$

$$+(2n)^{-1}\sum_{i=1}^{p}\sum_{j=1}^{p} U_{rij}(\tilde{\boldsymbol{\theta}}_n) \times (\hat{t}_{ni} - t_i^*) \times (\hat{t}_{nj} - t_j^*)$$

$$= n^{-1/2}Z_r - n^{-1}\sum_{i=1}^{p} Z_{ri}\overline{\overline{\bar{\kappa}^{i,u}Z_u}} - n^{-1/2}\sum_{i=1}^{p}\bar{\kappa}_{r,i}\overline{\overline{\bar{\kappa}^{i,u}Z_u}} + n^{-1}\sum_{i=1}^{p}(n^{1/2}Z_{ri} + n\bar{\kappa}_{ri})R_{ni}$$

$$+(2n)^{-1}\sum_{i=1}^{p}\sum_{j=1}^{p} U_{rij}(\tilde{\boldsymbol{\theta}}_n) \times (\hat{t}_{ni} - t_i^*) \times (\hat{t}_{nj} - t_j^*)$$

$$= -n^{-1}\overline{\overline{Z_{ri}\bar{\kappa}^{i,u}Z_u}} + n^{-1}\sum_{i=1}^{p}(n^{1/2}Z_{ri} + n\bar{\kappa}_{ri})R_{ni},$$

$$+(2n)^{-1}\sum_{i=1}^{p}\sum_{j=1}^{p} U_{rij}(\tilde{\boldsymbol{\theta}}_n) \times (\hat{t}_{ni} - t_i^*) \times (\hat{t}_{nj} - t_j^*), \tag{2.8}$$

where the third equality follows from the fact that $\kappa^i_{rs} = E(U^i_{rs}) = E(U^i_r U^i_s) = \kappa^i_{r,s}$, which implies $\bar{\kappa}_{rs} = \bar{\kappa}_{r,s}$ for $r$, $s$ = 1, 2, …, $p$; and the last equality follows from a cancellation due to the fact that $-n^{-1/2}\sum_{i=1}^{p}\bar{\kappa}_{r,i}\overline{\overline{\bar{\kappa}^{i,u}Z_u}} = -n^{-1/2}Z_r$. By condition A5, we know that by the CLT for i.n.i.d samples, both $Z_{ri}$ and $Z_u$ converge in distribution to a normal random variable. Thus, the first term in the last equality of (2.8) is $n\tilde{O}^2_d(n^{-1})$. By condition A6, we know that by the WLLN for i.n.i.d samples, $n^{-1}U_{rij}(\tilde{\boldsymbol{\theta}}_n)$ converges in probability to $\lim_{n\to\infty} \bar{\kappa}_{rij}$, which is a constant by condition A4(a), for $i$, $j$ = 1, 2, …, $p$. Thus, by Slutsky's theorem, $n^{-1}U_{rij}(\tilde{\boldsymbol{\theta}}_n) \times (\hat{t}_{ni} - t^*_i) \times (\hat{t}_{nj} - t^*_j) = nO^2_p(n^{-1})$ for $i, j$ = 1, 2, …, $p$, and, consequently, the third term in the last equality of (2.8) is $n\tilde{O}^2_d(n^{-1})$. Now, (2.8) can be rewritten as

$$n^{-1}\sum_{i=1}^{p}(n^{1/2}Z_{ri} + n\bar{\kappa}_{ri})R_{ni} = n\tilde{O}^2_d(n^{-1}) \tag{2.9}$$

Equation (2.9) holds for $r$ = 1, 2, …, $p$, which can be presented in the following matrix form:

$$\left(n^{-1/2}\boldsymbol{Z} + \bar{\boldsymbol{F}}_n(\boldsymbol{\theta}^*)\right)\boldsymbol{R}_n = n\tilde{\boldsymbol{O}}^2_d(n^{-1}), \tag{2.10}$$

where $\boldsymbol{Z} \equiv [Z_{ri}]_{r,\, i = 1, 2, \ldots, p}$, $\bar{\boldsymbol{F}}_n(\boldsymbol{\theta}^*) \equiv [\,\bar{\kappa}_{ri}\,]_{r,\, i = 1, 2, \ldots, p}$, $\boldsymbol{R}_n \equiv [R_{n1}, R_{n2}, \ldots, R_{np}]^T$, $\boldsymbol{O}_d(n^{-1})$ $\equiv [O_d(n^{-1}), O_d(n^{-1}), \ldots, O_d(n^{-1})]^T$, and $\tilde{\boldsymbol{O}}_d^2(n^{-1}) = [\tilde{O}_d^2(n^{-1}), \tilde{O}_d^2(n^{-1}), \ldots, \tilde{O}_d^2(n^{-1})]^T$. With condition A6, we know that by the WLLN for i.n.i.d samples, $n^{-1/2}\boldsymbol{Z}$ converges to zero in probability. Thus by the continuous mapping theorem (Mann and Wald, 1943) and condition A3, $\left(n^{-1/2}\boldsymbol{Z} + \bar{\boldsymbol{F}}_n(\boldsymbol{\theta}^*)\right)^{-1}$ converges in probability to $\boldsymbol{F}(\boldsymbol{\theta}^*)^{-1}$, i.e. $\left(n^{-1/2}\boldsymbol{Z} + \bar{\boldsymbol{F}}_n(\boldsymbol{\theta}^*)\right)^{-1} = \boldsymbol{F}(\boldsymbol{\theta}^*)^{-1} + \boldsymbol{o}_p(1)$, where $\boldsymbol{o}_p(1)$ is a $p \times p$ matrix with each entry being a $o_p(1)$ term. As a result,

$$
\begin{aligned}
\boldsymbol{R}_n &= \left(n^{-1/2}\boldsymbol{Z} + \bar{\boldsymbol{F}}_n(\boldsymbol{\theta}^*)\right)^{-1} \times n\tilde{\boldsymbol{O}}_d^2(n^{-1}) \\
&= \boldsymbol{F}^{-1}(\boldsymbol{\theta}^*) \times n\tilde{\boldsymbol{O}}_d^2(n^{-1}) + \boldsymbol{o}_p(1) \times n\tilde{\boldsymbol{O}}_d^2(n^{-1}) \\
&= n\tilde{\boldsymbol{O}}_d^2(n^{-1})
\end{aligned}
$$

Thus, $R_{ni} = n\tilde{O}_d^2(n^{-1})$, which, combined with (2.7), produces (2.5).

**Lemma 2**

For i.n.i.d sample data with conditions A1–A8 in Section 2.2, the inverse of the Fisher information matrix $\bar{\boldsymbol{F}}_n(\hat{\boldsymbol{\theta}}_n)^{-1}$ has the following expansion:

$$
\bar{\boldsymbol{F}}_n(\hat{\boldsymbol{\theta}}_n)^{-1}(r,s) = \bar{\kappa}^{r,s} + n^{-1/2}\overline{\overline{\bar{\kappa}^{r,j}\bar{\kappa}^{s,k}\bar{\kappa}^{i,u}(\bar{\kappa}_{jki} - \bar{\kappa}_{jk,i})Z_u}} + n\tilde{O}_d^2(n^{-1}), \tag{2.11}
$$

where $\bar{\boldsymbol{F}}_n(\hat{\boldsymbol{\theta}}_n)^{-1}(r,s)$ is the (*r*, *s*) element of $\bar{\boldsymbol{F}}_n(\hat{\boldsymbol{\theta}}_n)^{-1}$.

*Proof*: By a Taylor expansion around $\boldsymbol{\theta}^*$, the (*r*, *s*) element of $\bar{\boldsymbol{F}}_n(\hat{\boldsymbol{\theta}}_n)$ can be expressed as:

$$\bar{\boldsymbol{F}}_n(\hat{\boldsymbol{\theta}}_n)(r,s) = n^{-1}E(U_{rs})\Big|_{\boldsymbol{\theta}=\boldsymbol{\theta}^*} + n^{-1}\sum_{i=1}^{p}\frac{\partial E(U_{rs})}{\partial t_i}\Big|_{\boldsymbol{\theta}=\boldsymbol{\theta}^*} \times (\hat{t}_{ni} - t_i^*)$$

$$+(2n)^{-1}\sum_{i=1}^{p}\sum_{j=1}^{p}\frac{\partial^2 E(U_{rs})}{\partial t_i \partial t_j}\Big|_{\boldsymbol{\theta}=\tilde{\tilde{\boldsymbol{\theta}}}_n} \times (\hat{t}_{ni} - t_i^*) \times (\hat{t}_{nj} - t_j^*)$$

$$= n^{-1}E(U_{rs})\Big|_{\boldsymbol{\theta}=\boldsymbol{\theta}^*} - n^{-3/2}\sum_{i=1}^{p}\frac{\partial E(U_{rs})}{\partial t_i}\Big|_{\boldsymbol{\theta}=\boldsymbol{\theta}^*} \times \overline{\overline{\kappa^{i,u} Z_u}}$$

$$+n^{-1}\sum_{i=1}^{p}\frac{\partial E(U_{rs})}{\partial t_i}\Big|_{\boldsymbol{\theta}=\boldsymbol{\theta}^*} \times n\tilde{O}_d^2(n^{-1})$$

$$+(2n)^{-1}\sum_{i=1}^{p}\sum_{j=1}^{p}\frac{\partial^2 E(U_{rs})}{\partial t_i \partial t_j}\Big|_{\boldsymbol{\theta}=\tilde{\tilde{\boldsymbol{\theta}}}_n} \times (\hat{t}_{ni} - t_i^*) \times (\hat{t}_{nj} - t_j^*), \tag{2.12}$$

where $\tilde{\tilde{\boldsymbol{\theta}}}_n$ is an intermediate point between $\hat{\boldsymbol{\theta}}_n$ and $\boldsymbol{\theta}^*$ and the second equality follows from the result of Lemma 1. Notice that $n^{-1}\partial E(U_{rs})/\partial t_i\,|_{\boldsymbol{\theta}=\boldsymbol{\theta}^*}$ converges deterministically by condition A4(b), the third term in (2.12) after the second equality is $n\tilde{O}_d^2(n^{-1})$. With condition A6, we know that by the WLLN for i.n.i.d. samples, $n^{-1}\,\partial^2 E(U_{rs})\big/\partial t_i \partial t_j\Big|_{\boldsymbol{\theta}=\tilde{\tilde{\boldsymbol{\theta}}}_n}$ converges in probability to $\lim_{n\to\infty} n^{-1}\partial E(U_{rs})\big/\partial t_i\,|_{\boldsymbol{\theta}\,=\,\boldsymbol{\theta}^*}$, which is a constant by condition A4(b), for *i*, *j* = 1, 2, …, *p*. Thus, by Slutsky's theorem,

$n^{-1}\,\partial^2 E(U_{rs})\big/\partial t_i \partial t_j \Big|_{\boldsymbol{\theta}=\tilde{\tilde{\boldsymbol{\theta}}}_n} \times(\hat{t}_{ni}-t_i^*)\times(\hat{t}_{nj}-t_j^*) \;= nO_d^2(n^{-1})$ for $i$, $j$ = 1, 2, …, $p$ and consequently, the fourth term in the last equality of (2.8) is $n\tilde{O}_d^2(n^{-1})$. As a result, expression (2.12) is equivalent to the following:

$$\bar{\boldsymbol{F}}_n(\hat{\boldsymbol{\theta}}_n)(r,s) = n^{-1}E(U_{rs})\Big|_{\boldsymbol{\theta}=\boldsymbol{\theta}^*} - n^{-3/2}\sum_{i=1}^{p}\frac{\partial E(U_{rs})}{\partial t_i}\Big|_{\boldsymbol{\theta}=\boldsymbol{\theta}^*} \times \overline{\overline{\bar{\kappa}^{i,u} Z_u}} + n\tilde{O}_d^2(n^{-1}). \qquad (2.13)$$

We now claim the following two facts:

(i) $n^{-1}E(U_{rs})\Big|_{\boldsymbol{\theta}=\boldsymbol{\theta}^*} = \bar{\kappa}_{r,s}$;

(ii) $n^{-3/2}\sum_{i=1}^{p}\frac{\partial E(U_{rs})}{\partial t_i}\Big|_{\boldsymbol{\theta}=\boldsymbol{\theta}^*} \times \overline{\overline{\bar{\kappa}^{i,u} Z_u}} = n^{-1/2}\overline{\overline{\bar{\kappa}^{i,u}(\bar{\kappa}_{rsi} - \bar{\kappa}_{rs,i})Z_u}}$.

First, (i) follows from the definition of $U_{rs}$ and the equivalent form of FIM in (1.2):

$$n^{-1}E(U_{rs})\Big|_{\boldsymbol{\theta}=\boldsymbol{\theta}^*} = n^{-1}\sum_{j=1}^{n}E(U_{rs}^j)\Big|_{\boldsymbol{\theta}=\boldsymbol{\theta}^*} = n^{-1}\sum_{j=1}^{n}E(U_r^j U_s^j)\Big|_{\boldsymbol{\theta}=\boldsymbol{\theta}^*} = n^{-1}\sum_{j=1}^{n}\kappa_{r,s}^j = \bar{\kappa}_{r,s}.$$

To show (ii), we first rewrite $\partial E(U_{rs})\big/\partial t_i$ by definition:

$$\frac{\partial E(U_{rs})}{\partial t_i} = \frac{\partial}{\partial t_i}\sum_{j=1}^{n}E(U_{rs}^j)$$

$$= \sum_{j=1}^{n} \frac{\partial E(U_{rs}^{j})}{\partial t_i}$$

$$= \sum_{j=1}^{n} \frac{\partial}{\partial t_i} \left( \sum_{\boldsymbol{x}_j^d \in \boldsymbol{S}_j^d} \left[ \int_{\boldsymbol{x}_j^c \in \boldsymbol{S}_j^c | \boldsymbol{x}_j^d} U_{rs}^{j} p_j^c(\boldsymbol{x}_j^c, \boldsymbol{\theta} | \boldsymbol{x}_j^d) d\boldsymbol{x}_j^c \right] p_j^d(\boldsymbol{x}_j^d, \boldsymbol{\theta}) \right)$$

Furthermore, by condition A1,

$$\frac{\partial}{\partial t_i} \left( \sum_{\boldsymbol{x}_j^d \in \boldsymbol{S}_j^d} \left[ \int_{\boldsymbol{x}_j^c \in \boldsymbol{S}_j^c | \boldsymbol{x}_j^d} U_{rs}^{j} p_j^c(\boldsymbol{x}_j^c, \boldsymbol{\theta} | \boldsymbol{x}_j^d) d\boldsymbol{x}_j^c \right] p_j^d(\boldsymbol{x}_j^d, \boldsymbol{\theta}) \right)$$

$$= \frac{\partial}{\partial t_i} \left( \sum_{\boldsymbol{x}_j^d \in \boldsymbol{S}_j^d} \int_{\boldsymbol{x}_j^c \in \boldsymbol{S}_j^c | \boldsymbol{x}_j^d} U_{rs}^{j} p_j^c(\boldsymbol{x}_j^c, \boldsymbol{\theta} | \boldsymbol{x}_j^d) p_j^d(\boldsymbol{x}_j^d, \boldsymbol{\theta}) d\boldsymbol{x}_j^c \right)$$

$$= \frac{\partial}{\partial t_i} \left( \sum_{\boldsymbol{x}_j^d \in \boldsymbol{S}_j^d} \int_{\boldsymbol{x}_j^c \in \boldsymbol{S}_j^c | \boldsymbol{x}_j^d} U_{rs}^{j} \exp\left(-l(\boldsymbol{x}_j, \boldsymbol{\theta})\right) d\boldsymbol{x}_j^c \right)$$

$$= \sum_{\boldsymbol{x}_j^d \in \boldsymbol{S}_j^d} \int_{\boldsymbol{x}_j^c \in \boldsymbol{S}_j^c | \boldsymbol{x}_j^d} \left( U_{rsi}^{j} \exp\left(-l(\boldsymbol{x}_j, \boldsymbol{\theta})\right) - U_{rs}^{j} U_i^{j} p_j(\boldsymbol{x}_j, \boldsymbol{\theta}) \right) d\boldsymbol{x}_i^c$$

(interchange of differentiation and integration)

$$= \sum_{\boldsymbol{x}_j^d \in \boldsymbol{S}_j^d} \int_{\boldsymbol{x}_j^c \in \boldsymbol{S}_j^c | \boldsymbol{x}_j^d} \left( U_{rsi}^{j} - U_{rs}^{j} U_i^{j} \right) p_j(\boldsymbol{x}_j, \boldsymbol{\theta}) d\boldsymbol{x}_i^c$$

$$= \sum_{\boldsymbol{x}_j^d \in \boldsymbol{S}_j^d} \int_{\boldsymbol{x}_j^c \in \boldsymbol{S}_j^c | \boldsymbol{x}_j^d} \left( U_{rsi}^{j} - U_{rs}^{j} U_i^{j} \right) p_j^c(\boldsymbol{x}_j^c, \boldsymbol{\theta} | \boldsymbol{x}_j^d) d\boldsymbol{x}_j^c \times p_j^d(\boldsymbol{x}_j^d, \boldsymbol{\theta})$$

$$= E\left( U_{rsi}^{j} - U_{rs}^{j} U_i^{j} \right).$$

Thus,

$$n^{-3/2}\sum_{i=1}^{p}\frac{\partial E(U_{rs})}{\partial t_i}\Big|_{\boldsymbol{\theta}=\boldsymbol{\theta}^*} = n^{-3/2}\sum_{i=1}^{p}\sum_{j=1}^{n}E\left(U_{rsi}^{j} - U_{rs}^{j}U_{i}^{j}\right)\Big|_{\boldsymbol{\theta}=\boldsymbol{\theta}^*}$$

$$= n^{-3/2}\sum_{i=1}^{p}\sum_{j=1}^{n}\left(\kappa_{rsi}^{j} - \kappa_{rs,i}^{j}\right)$$

$$= n^{-1/2}\sum_{i=1}^{p}(\overline{\kappa}_{rsi} - \overline{\kappa}_{rs,i}),$$

$$n^{-3/2}\sum_{i=1}^{p}\frac{\partial E(U_{rs})}{\partial t_i}\Big|_{\boldsymbol{\theta}=\boldsymbol{\theta}^*}\times\overline{\overline{\overline{\kappa}^{i,u}Z_u}} = n^{-1/2}\sum_{i=1}^{p}(\overline{\kappa}_{rsi} - \overline{\kappa}_{rs,i})\times\overline{\overline{\overline{\kappa}^{i,u}Z_u}} = n^{-1/2}\overline{\overline{\overline{\kappa}^{i,u}(\overline{\kappa}_{rsi} - \overline{\kappa}_{rs,i})Z_u}}\,.$$

Given (i), (ii), we re-express $\overline{\boldsymbol{F}}_n(\hat{\boldsymbol{\theta}}_n)$ in (2.13) as follows:

$$\overline{\boldsymbol{F}}_n(\hat{\boldsymbol{\theta}}_n)(r,s) = \overline{\kappa}_{r,s} - n^{-1/2}\overline{\overline{\overline{\kappa}^{i,u}(\overline{\kappa}_{rsi} - \overline{\kappa}_{rs,i})Z_u}} + n\tilde{O}_d^2(n^{-1})\,. \tag{2.14}$$

By the definition of matrix inverse,

$$\sum_{s=1}^{p}\overline{\boldsymbol{F}}_n(\hat{\boldsymbol{\theta}}_n)(r,s)\times\overline{\boldsymbol{F}}_n(\hat{\boldsymbol{\theta}}_n)^{-1}(s,t) = \begin{cases}1 & \text{if } r = t\\ 0 & \text{if } r \neq t.\end{cases} \tag{2.15}$$

We now develop the form for $\bar{\boldsymbol{F}}_n(\hat{\boldsymbol{\theta}}_n)^{-1}(r,s)$ in order to satisfy (2.15). Given the expression in (2.11), let us suppose $\bar{\boldsymbol{F}}_n(\hat{\boldsymbol{\theta}}_n)^{-1}(r,s)$ has the following representation:

$$\bar{\boldsymbol{F}}_n(\hat{\boldsymbol{\theta}}_n)^{-1}(r,s) = \bar{\kappa}^{r,s} + n^{-1/2}\overline{\overline{\bar{\kappa}^{r,j}\bar{\kappa}^{s,k}\bar{\kappa}^{i,u}(\bar{\kappa}_{jki}-\bar{\kappa}_{jk,i})Z_u}} + W_n(r,s), \tag{2.16}$$

where $W_n\,(r, s)$ is to be determined. In fact, we want to show that $W_n\,(r, s) = n\tilde{O}_d^2(n^{-1})$, for $r, s = 1, 2, \ldots, p$.

By plugging (2.14) and (2.16) into (2.15), we have:

$$\sum_{s=1}^{p}\bar{\boldsymbol{F}}_n(\hat{\boldsymbol{\theta}}_n)(r,s)\times\bar{\boldsymbol{F}}_n(\hat{\boldsymbol{\theta}}_n)^{-1}(s,t)$$

$$=\sum_{s=1}^{p}\left[\left(\bar{\kappa}_{r,s} - n^{-1/2}\overline{\overline{\bar{\kappa}^{i,u}(\bar{\kappa}_{rsi}-\bar{\kappa}_{rs,i})Z_u}} + n\tilde{O}_d^2(n^{-1})\right)\right.$$

$$\left.\times\left(\bar{\kappa}^{s,t} + n^{-1/2}\overline{\overline{\bar{\kappa}^{s,j}\bar{\kappa}^{t,k}\bar{\kappa}^{i,u}(\bar{\kappa}_{jki}-\bar{\kappa}_{jk,i})Z_u}} + W_n(s,t)\right)\right]$$

$$=\overline{\overline{\bar{\kappa}_{r,s}\bar{\kappa}^{s,t}}} - n^{-1/2}\sum_{s=1}^{p}\bar{\kappa}^{s,t}\overline{\overline{\bar{\kappa}^{i,u}(\bar{\kappa}_{rsi}-\bar{\kappa}_{rs,i})Z_u}} + n^{-1/2}\sum_{s=1}^{p}\bar{\kappa}_{r,s}\overline{\overline{\bar{\kappa}^{s,j}\bar{\kappa}^{t,k}\bar{\kappa}^{i,u}(\bar{\kappa}_{jki}-\bar{\kappa}_{jk,i})Z_u}}$$

$$-\sum_{s=1}^{p}n^{-1/2}\overline{\overline{\bar{\kappa}^{i,u}(\bar{\kappa}_{rsi}-\bar{\kappa}_{rs,i})Z_u}}\times W_n(s,t) + n^{-1/2}\sum_{s=1}^{p}\overline{\overline{\bar{\kappa}^{s,j}\bar{\kappa}^{t,k}\bar{\kappa}^{i,u}(\bar{\kappa}_{jki}-\bar{\kappa}_{jk,i})Z_u}}\times n\tilde{O}_d^2(n^{-1})$$

$$-n^{-1}\sum_{s=1}^{p}\overline{\overline{\bar{\kappa}^{i,u}(\bar{\kappa}_{rsi}-\bar{\kappa}_{rs,i})Z_u}}\times\overline{\overline{\bar{\kappa}^{s,j}\bar{\kappa}^{t,k}\bar{\kappa}^{i,u}(\bar{\kappa}_{jki}-\bar{\kappa}_{jk,i})Z_u}} + \sum_{s=1}^{p}\bar{\kappa}_{r,s}\times W_n(s,t)$$

$$+\sum_{s=1}^{p}\bar{\kappa}^{s,t}\times n\tilde{O}_d^2(n^{-1})+\sum_{s=1}^{p}W_n(s,t)\times n\tilde{O}_d^2(n^{-1})$$

(expand the product and pass the summation sign to each individual term)

$$=\overline{\overline{\bar{\kappa}_{r,s}\bar{\kappa}^{s,t}}}-n^{-1/2}\overline{\overline{\bar{\kappa}^{s,t}\bar{\kappa}^{i,u}(\bar{\kappa}_{rsi}-\bar{\kappa}_{rs,i})Z_u}}+n^{-1/2}\overline{\overline{\bar{\kappa}^{t,k}\bar{\kappa}^{i,u}(\bar{\kappa}_{rki}-\bar{\kappa}_{rk,i})Z_u}}$$

$$-n^{-1/2}\sum_{s=1}^{p}\overline{\overline{\bar{\kappa}^{i,u}(\bar{\kappa}_{rsi}-\bar{\kappa}_{rs,i})Z_u\times W_n(s,t)}}$$

$$+n^{-1/2}\sum_{s=1}^{p}\overline{\overline{\bar{\kappa}^{s,j}\bar{\kappa}^{t,k}\bar{\kappa}^{i,u}(\bar{\kappa}_{jki}-\bar{\kappa}_{jk,i})Z_u\times n\tilde{O}_d^2(n^{-1})}}$$

$$-n^{-1}\sum_{s=1}^{p}\overline{\overline{\bar{\kappa}^{i,u}(\bar{\kappa}_{rsi}-\bar{\kappa}_{rs,i})Z_u}}\times\overline{\overline{\bar{\kappa}^{s,j}\bar{\kappa}^{t,k}\bar{\kappa}^{i,u}(\bar{\kappa}_{jki}-\bar{\kappa}_{jk,i})Z_u}}+\sum_{s=1}^{p}\bar{\kappa}_{r,s}\times W_n(s,t)$$

$$+\sum_{s=1}^{p}\bar{\kappa}^{s,t}\times n\tilde{O}_d^2(n^{-1})+\sum_{s=1}^{p}W_n(s,t)\times n\tilde{O}_d^2(n^{-1})$$

$$=\overline{\overline{\bar{\kappa}_{r,s}\bar{\kappa}^{s,t}}}-n^{-1/2}\sum_{s=1}^{p}\overline{\overline{\bar{\kappa}^{i,u}(\bar{\kappa}_{rsi}-\bar{\kappa}_{rs,i})Z_u}}\times W_n(s,t)$$

$$+n^{-1/2}\sum_{s=1}^{p}\overline{\overline{\bar{\kappa}^{s,j}\bar{\kappa}^{t,k}\bar{\kappa}^{i,u}(\bar{\kappa}_{jki}-\bar{\kappa}_{jk,i})Z_u}}\times n\tilde{O}_d^2(n^{-1})$$

$$-n^{-1}\sum_{s=1}^{p}\overline{\overline{(\bar{\kappa}_{rsi}-\bar{\kappa}_{rs,i})\bar{\kappa}^{i,u}Z_u}}\times\overline{\overline{\bar{\kappa}^{s,j}\bar{\kappa}^{t,k}\bar{\kappa}^{i,u}(\bar{\kappa}_{jki}-\bar{\kappa}_{jk,i})Z_u}}+\sum_{s=1}^{p}\bar{\kappa}_{r,s}\times W_n(s,t)$$

$$+\sum_{s=1}^{p}\bar{\kappa}^{s,t}\times n\tilde{O}_d^2(n^{-1})+\sum_{s=1}^{p}W_n(s,t)\times n\tilde{O}_d^2(n^{-1}). \qquad (2.17)$$

By the definition of matrix inverse, the first term after the last equality above is:

$$\overline{\overline{\bar{\kappa}_{r,s}\bar{\kappa}^{s,t}}} = \begin{cases} 1 & \text{if } r = t \\ 0 & \text{if } r \neq t. \end{cases}$$

As a result, in order for expression (2.17) to equal the right hand side of (2.15), we must have the rest of the terms in the last equation of (2.17) sum up to zero, i.e.,

$$-n^{-1/2}\sum_{s=1}^{p}\overline{\overline{\bar{\kappa}^{i,u}(\bar{\kappa}_{rsi}-\bar{\kappa}_{rs,i})Z_u}}\times W_n(s,t)+n^{-1/2}\sum_{s=1}^{p}\overline{\overline{\bar{\kappa}^{s,j}\bar{\kappa}^{t,k}\bar{\kappa}^{i,u}(\bar{\kappa}_{jki}-\bar{\kappa}_{jk,i})Z_u}}\times n\tilde{O}_d^2(n^{-1})$$

$$-n^{-1}\sum_{s=1}^{p}\overline{\overline{(\bar{\kappa}_{rsi}-\bar{\kappa}_{rs,i})\bar{\kappa}^{i,u}Z_u}}\times\overline{\overline{\bar{\kappa}^{s,j}\bar{\kappa}^{t,k}\bar{\kappa}^{i,u}(\bar{\kappa}_{jki}-\bar{\kappa}_{jk,i})Z_u}}+\sum_{s=1}^{p}\bar{\kappa}_{r,s}\times W_n(s,t)$$

$$+\sum_{s=1}^{p}\bar{\kappa}^{s,t}\times n\tilde{O}_d^2(n^{-1})+\sum_{s=1}^{p}W_n(s,t)\times n\tilde{O}_d^2(n^{-1})=0. \tag{2.18}$$

Group the left hand side of (2.18) by $W_n(s, t)$, we have:

$$\sum_{s=1}^{p}W_n(s,t)\times\left(-n^{-1/2}\overline{\overline{\bar{\kappa}^{i,u}(\bar{\kappa}_{rsi}-\bar{\kappa}_{rs,i})Z_u}}+\bar{\kappa}_{r,s}+n\tilde{O}_d^2(n^{-1})\right)$$

$$=-n^{-1/2}\sum_{s=1}^{p}\overline{\overline{\bar{\kappa}^{s,j}\bar{\kappa}^{t,k}\bar{\kappa}^{i,u}(\bar{\kappa}_{jki}-\bar{\kappa}_{jk,i})Z_u}}\times n\tilde{O}_d^2(n^{-1})-\sum_{s=1}^{p}\bar{\kappa}^{s,t}\times n\tilde{O}_d^2(n^{-1})$$

$$+n^{-1}\sum_{s=1}^{p}\overline{\overline{(\bar{\kappa}_{rsi}-\bar{\kappa}_{rs,i})\bar{\kappa}^{i,u}Z_u}}\times\overline{\overline{\bar{\kappa}^{s,j}\bar{\kappa}^{t,k}\bar{\kappa}^{i,u}(\bar{\kappa}_{jki}-\bar{\kappa}_{jk,i})Z_u}}. \tag{2.19}$$

For the left hand side of (2.19), $-n^{-1/2}\overline{\overline{\bar{\kappa}^{i,u}(\bar{\kappa}_{rsi}-\bar{\kappa}_{rs,i})Z_u}} + n\tilde{O}_d^2(n^{-1}) = o_p(1)$, which follows from Slutsky's theorem. For the right hand side of (2.19), by Slutsky's therorem, the first term $-n^{-1/2}\sum_{s=1}^{p}\overline{\overline{\bar{\kappa}^{s,j}\bar{\kappa}^{t,k}\bar{\kappa}^{i,u}(\bar{\kappa}_{jki}-\bar{\kappa}_{jk,i})Z_u}} \times n\tilde{O}_d^2(n^{-1}) = o_p(n^{-1})$. The second term $-\sum_{s=1}^{p}\bar{\kappa}^{s,t} \times n\tilde{O}_d^2(n^{-1}) = n\tilde{O}_d^2(n^{-1})$ by the fact that $\bar{\kappa}^{s,t} = O(1)$ by condition A3. The last term $n^{-1}\sum_{s=1}^{p}\overline{\overline{(\bar{\kappa}_{rsi}-\bar{\kappa}_{rs,i})\bar{\kappa}^{i,u}Z_u}} \times \overline{\overline{\bar{\kappa}^{s,j}\bar{\kappa}^{t,k}\bar{\kappa}^{i,u}(\bar{\kappa}_{jki}-\bar{\kappa}_{jk,i})Z_u}} = n\tilde{O}_d^2(n^{-1})$. As a result, (2.19) can be rewritten as

$$\sum_{s=1}^{p} W_n(s,t)\times\left(\bar{\kappa}_{r,s}+o_p(1)\right) = n\tilde{O}_d^2(n^{-1}), \tag{2.20}$$

where the left hand side of (2.20) uses the comment below (2.19) and the right hand side follows from the analysis directly above. Equation (2.20) holds for all $r$, $t$ = 1, 2, …, $p$. A matrix form representation is:

$$\boldsymbol{W}_n \times\left(\bar{\boldsymbol{F}}_n(\boldsymbol{\theta}^*)+\boldsymbol{o}_p(1)\right) = n\tilde{\boldsymbol{O}}_d^2(n^{-1}), \tag{2.21}$$

where $\boldsymbol{W}_n \equiv \{W_n(s, t)\}_{s,\, t\, =\, 1,\, 2,\, \ldots,\, p}$, $\bar{\boldsymbol{F}}_n(\boldsymbol{\theta}^*) \equiv \{\kappa_{r,\, s}\}_{r,\, s\, =\, 1,\, 2,\, \ldots,\, p}$, $\boldsymbol{o}_p(1)$ is a $p \times p$ matrix with each element being $o_p(1)$, and $\tilde{\boldsymbol{O}}_d^2(n^{-1})$ is a $p \times p$ matrix with each element being $\tilde{O}_d^2(n^{-1})$. By condition A3, $\bar{\boldsymbol{F}}_n(\boldsymbol{\theta}^*)$ converges to an invertible constant matrix $\boldsymbol{F}(\boldsymbol{\theta}^*)$, which implies that $\left(\bar{\boldsymbol{F}}_n(\boldsymbol{\theta}^*)+\boldsymbol{o}_p(1)\right)^{-1}$ converges in probability to $\boldsymbol{F}(\boldsymbol{\theta}^*)^{-1}$. By the continuous mapping theorem (Mann and Wald, 1943), we have the following:

$$\boldsymbol{W}_n = n\tilde{\boldsymbol{O}}_d^2(n^{-1})\times\left(\boldsymbol{F}(\boldsymbol{\theta}^*)^{-1}+\boldsymbol{o}_p(1)\right)$$

$$= n\tilde{\boldsymbol{O}}_d^2(n^{-1})\times\boldsymbol{F}(\boldsymbol{\theta}^*)^{-1}+n\tilde{\boldsymbol{O}}_d^2(n^{-1})\times\boldsymbol{o}_p(1).$$

As a result, each element in $\boldsymbol{W}_n$ is a linear combination of $n\tilde{O}_d^2(n^{-1})$ terms plus some higher order terms which are dominated by $n\tilde{O}_d^2(n^{-1})$. Thus, $W_n(r, s) = n\tilde{O}_d^2(n^{-1})$, for all $r, s = 1, 2, \ldots, p$. Consequently, by (2.16),

$$\bar{\boldsymbol{F}}_n(\hat{\boldsymbol{\theta}}_n)^{-1}(r,s)=\bar{\kappa}^{r,s}+n^{-1/2}\overline{\overline{\bar{\kappa}^{r,j}\bar{\kappa}^{s,k}\bar{\kappa}^{i,u}(\bar{\kappa}_{jki}-\bar{\kappa}_{jk,i})Z_u}}+n\tilde{O}_d^2(n^{-1}).\qquad\square$$

**Lemma 3**

For i.n.i.d sample data with conditions A1–A8 in Section 2.2, the inverse of the observed Fisher information matrix $\bar{\boldsymbol{H}}_n(\hat{\boldsymbol{\theta}}_n)^{-1}$ has the following expansion:

$$\bar{\boldsymbol{H}}_n(\hat{\boldsymbol{\theta}}_n)^{-1}(r,s)=\bar{\kappa}^{r,s}+n^{-1/2}\overline{\overline{\bar{\kappa}^{r,j}\bar{\kappa}^{s,k}\left(\bar{\kappa}^{i,u}(\bar{\kappa}_{jki}-\bar{\kappa}_{jk,i})Z_u-Y_{jk}\right)}}+n\tilde{O}_d^2(n^{-1}),\qquad(2.22)$$

where $\bar{\boldsymbol{H}}_n(\hat{\boldsymbol{\theta}}_n)^{-1}(r,s)$ is the $(r, s)$ element of $\bar{\boldsymbol{H}}_n(\hat{\boldsymbol{\theta}}_n)^{-1}$.

*Proof*: By Lemma 1 and the Taylor's expansion, the $(r,s)$ element of $\bar{\boldsymbol{H}}_n(\hat{\boldsymbol{\theta}}_n)$ can be expanded as follows:

$$\bar{\boldsymbol{H}}_n(\hat{\boldsymbol{\theta}}_n)(r,s) = n^{-1}U_{rs}(\boldsymbol{\theta}^*) + n^{-1}\sum_{i=1}^{p} U_{rsi}(\boldsymbol{\theta}^*)\times(\hat{t}_{ni}-t_i^*)$$

$$+(2n)^{-1}\sum_{i=1}^{p} U_{rsij}(\tilde{\tilde{\boldsymbol{\theta}}}_n)\times(\hat{t}_{ni}-t_i^*)\times(\hat{t}_{nj}-t_j^*)$$

$$= n^{-1}\left(n^{1/2}Z_{rs} + n\bar{\kappa}_{rs}\right) - n^{-3/2}\sum_{i=1}^{p} U_{rsi}(\boldsymbol{\theta}^*)\overline{\overline{\kappa^{i,u}}Z_u} + n\tilde{O}_d^2(n^{-1})\times n^{-1}\sum_{i=1}^{p} U_{rsi}(\boldsymbol{\theta}^*)$$

$$+(2n)^{-1}\sum_{i=1}^{p}\sum_{j=1}^{p} U_{rsij}(\tilde{\tilde{\boldsymbol{\theta}}}_n)\times(\hat{t}_{ni}-t_i^*)\times(\hat{t}_{nj}-t_j^*) \qquad (2.23)$$

where $\tilde{\tilde{\boldsymbol{\theta}}}_n$ is an intermediate point between $\hat{\boldsymbol{\theta}}_n$ and $\boldsymbol{\theta}^*$. With condition A6, we know that by the WLLN for i.n.i.d samples, $n^{-1}\sum_{i=1}^{p} U_{rsi}(\boldsymbol{\theta}^*)$ converges in probability to $\sum_{i=1}^{p} \lim_{n\to\infty} \bar{\kappa}_{rsi}$, which exists and is constant by condition A4(a). Thus, the third term after the last equality in (2.23) is $n\tilde{O}_d^2(n^{-1})$ by Slutsky's theorem. In addition, again by condition A6 and the WLLN for i.n.i.d samples, $n^{-1}U_{rsij}(\tilde{\tilde{\boldsymbol{\theta}}}_n)$ converges in probability to $\lim_{n\to\infty} n^{-1} E(U_{rsij})|_{\boldsymbol{\theta}=\boldsymbol{\theta}^*}$, which is a constant by condition A4(d), for $i, j = 1, 2, \ldots, p$. Thus, by Slutsky's theorem, $n^{-1}U_{rsij}(\tilde{\tilde{\boldsymbol{\theta}}}_n)\times(\hat{t}_{ni}-t_i^*)\times(\hat{t}_{nj}-t_j^*) = nO_d^2(n^{-1})$ for $i, j = 1,$ 2, …, $p$ and consequently, the fourth term in the last equality of (2.23) is $n\tilde{O}_d^2(n^{-1})$. As a result, expression (2.23) is equivalent to the following:

$$\bar{\boldsymbol{H}}_n(\hat{\boldsymbol{\theta}}_n)(r,s) = n^{-1}\left(n^{1/2}Z_{rs} + n\bar{\kappa}_{rs}\right) - n^{-3/2}\sum_{i=1}^{p}\left(U_{rsi}(\boldsymbol{\theta}^*) - n\bar{\kappa}_{rsi} + n\bar{\kappa}_{rsi}\right)\overline{\overline{\bar{\kappa}^{i,u}Z_u}} + n\tilde{O}_d^2(n^{-1})$$

$$= n^{-1/2}Z_{rs} + \bar{\kappa}_{rs} - n^{-1/2}\sum_{i=1}^{p}\overline{\overline{\bar{\kappa}_{rsi}\bar{\kappa}^{i,u}Z_u}}$$

$$-n^{-3/2}\sum_{i=1}^{p}\left(U_{rsi}(\boldsymbol{\theta}^*) - n\bar{\kappa}_{rsi}\right)\overline{\overline{\bar{\kappa}^{i,u}Z_u}} + n\tilde{O}_d^2(n^{-1})$$

$$= \bar{\kappa}_{rs} + n^{-1/2}\left(Z_{rs} - \overline{\overline{\bar{\kappa}_{rsi}\bar{\kappa}^{i,u}Z_u}}\right)$$

$$-n^{-3/2}\sum_{i=1}^{p}\left(U_{rsi}(\boldsymbol{\theta}^*) - n\bar{\kappa}_{rsi}\right)\overline{\overline{\bar{\kappa}^{i,u}Z_u}} + n\tilde{O}_d^2(n^{-1})$$

$$= \bar{\kappa}_{r,s} + n^{-1/2}\left(Z_{rs} - \overline{\overline{\bar{\kappa}_{rsi}\bar{\kappa}^{i,u}Z_u}}\right) + n\tilde{O}_d^2(n^{-1}), \qquad (2.24)$$

where the last equation follows from the facts that $\bar{\kappa}_{rs} = \sum_{i=1}^{n}\kappa_{rs}^{i}/n = \sum_{i=1}^{n}E(U_{rs}^{i})/n$ $= \sum_{i=1}^{n}E\left(U_r^iU_s^i\right)/n = \sum_{i=1}^{n}\kappa_{r,s}^{i}/n = \bar{\kappa}_{r,s}$ and that $n^{-1/2}\left(U_{rsi}(\boldsymbol{\theta}^*) - n\bar{\kappa}_{rsi}\right)$ converges in distribution to a normal random variable, according to condition A5 and the CLT for i.n.i.d data. Thus $n^{-3/2}\sum_{i=1}^{p}\left(U_{rsi}(\boldsymbol{\theta}^*) - n\bar{\kappa}_{rsi}\right) \times\overline{\overline{\bar{\kappa}^{i,u}Z_u}} = n\tilde{O}_d^2(n^{-1})$.

To show (2.22), let us first assume that $\bar{\boldsymbol{H}}_n(\hat{\boldsymbol{\theta}}_n)^{-1}(r,s)$ has the following form:

$$\bar{\boldsymbol{H}}_n(\hat{\boldsymbol{\theta}}_n)^{-1}(r,s) = \bar{\kappa}^{r,s} + n^{-1/2}\overline{\overline{\bar{\kappa}^{r,j}\bar{\kappa}^{s,k}\left((\bar{\kappa}_{jki} - \bar{\kappa}_{jk,i})\bar{\kappa}^{i,u}Z_u - Y_{jk}\right)}} + V_n(r,s)\,, \qquad (2.25)$$

where $V_n(r,s)$ is to be determined. We now want to show that $V_n(r,s) = n\tilde{O}_d^2(n^{-1})$. By definition of matrix inverse, expression (2.25) must satisfy the following:

$$\sum_{s=1}^{p}\bar{\boldsymbol{H}}_n(\hat{\boldsymbol{\theta}}_n)(r,s)\times\bar{\boldsymbol{H}}_n(\hat{\boldsymbol{\theta}}_n)^{-1}(s,t) = \begin{cases} 1 & \text{if } r = t, \\ 0 & \text{if } r \neq t. \end{cases} \qquad (2.26)$$

Plugging (2.24) and (2.25) into (2.26), we have the following:

$$\begin{aligned}
&\sum_{s=1}^{p}\bar{\boldsymbol{H}}_n(\hat{\boldsymbol{\theta}}_n)(r,s)\times\bar{\boldsymbol{H}}_n(\hat{\boldsymbol{\theta}}_n)^{-1}(r,s) \\
&= \sum_{s=1}^{p}\left[\left(\bar{\kappa}_{r,s} + n^{-1/2}\left(Z_{rs} - \overline{\overline{\bar{\kappa}_{rsi}\bar{\kappa}^{i,u}Z_u}}\right) + n\tilde{O}_d^2(n^{-1})\right)\right. \\
&\qquad \left.\times\left(\bar{\kappa}^{s,t} + n^{-1/2}\overline{\overline{\bar{\kappa}^{s,j}\bar{\kappa}^{t,k}\left((\bar{\kappa}_{jki} - \bar{\kappa}_{jk,i})\bar{\kappa}^{i,u}Z_u - Y_{jk}\right)}} + V_n(r,s)\right)\right] \\
&= \sum_{s=1}^{p}\bar{\kappa}_{r,s}\bar{\kappa}^{s,t} + n^{-1/2}\sum_{s=1}^{p}\bar{\kappa}_{r,s}\overline{\overline{\bar{\kappa}^{s,j}\bar{\kappa}^{t,k}\left((\bar{\kappa}_{jki} - \bar{\kappa}_{jk,i})\bar{\kappa}^{i,u}Z_u - Y_{jk}\right)}} \\
&\quad + n^{-1/2}\sum_{s=1}^{p}\left(Z_{rs} - \overline{\overline{\bar{\kappa}_{rsi}\bar{\kappa}^{i,u}Z_u}}\right)\times\bar{\kappa}^{s,t} \\
&\quad + n^{-1}\sum_{s=1}^{p}\left(Z_{rs} - \overline{\overline{\bar{\kappa}_{rsi}\bar{\kappa}^{i,u}Z_u}}\right)\times\overline{\overline{\bar{\kappa}^{s,j}\bar{\kappa}^{t,k}\left((\bar{\kappa}_{jki} - \bar{\kappa}_{jk,i})\bar{\kappa}^{i,u}Z_u - Y_{jk}\right)}}
\end{aligned}$$

$$+n^{-1/2}\sum_{s=1}^{p}\left(Z_{rs}-\overline{\overline{\bar{\kappa}_{rsi}\bar{\kappa}^{i,u}Z_u}}\right)V_n(r,s)+\sum_{s=1}^{p}\bar{\kappa}_{r,s}\times V_n(r,s)+\sum_{s=1}^{p}\bar{\kappa}^{s,t}\times n\tilde{O}_d^2(n^{-1})$$

$$+n^{-1/2}\sum_{s=1}^{p}\overline{\overline{\bar{\kappa}^{s,j}\bar{\kappa}^{t,k}\left((\bar{\kappa}_{jki}-\bar{\kappa}_{jk,i})\bar{\kappa}^{i,u}Z_u-Y_{jk}\right)}}\times n\tilde{O}_d^2(n^{-1})+\sum_{s=1}^{p}V_n(r,s)\times n\tilde{O}_d^2(n^{-1})$$

(Expand the product and pass the summation sign to each individual term)

$$=\sum_{s=1}^{p}\bar{\kappa}_{r,s}\bar{\kappa}^{s,t}+n^{-1/2}\sum_{s=1}^{p}\bar{\kappa}_{r,s}\overline{\overline{\bar{\kappa}^{s,j}\bar{\kappa}^{t,k}\left((\bar{\kappa}_{jki}-\bar{\kappa}_{jk,i})\bar{\kappa}^{i,u}Z_u-Y_{jk}\right)}}$$

$$+n^{-1/2}\sum_{s=1}^{p}\left(Y_{rs}+\overline{\overline{\bar{\kappa}_{rs,i}\bar{\kappa}^{i,u}Z_u}}-\overline{\overline{\bar{\kappa}_{rsi}\bar{\kappa}^{i,u}Z_u}}\right)\times\bar{\kappa}^{s,t}+n^{-1/2}\sum_{s=1}^{p}\left(Z_{rs}-\overline{\overline{\bar{\kappa}_{rsi}\bar{\kappa}^{i,u}Z_u}}\right)\times V_n(r,s)$$

$$+n^{-1}\sum_{s=1}^{p}\left(Z_{rs}-\overline{\overline{\bar{\kappa}_{rsi}\bar{\kappa}^{i,u}Z_u}}\right)\times\overline{\overline{\bar{\kappa}^{s,j}\bar{\kappa}^{t,k}\left((\bar{\kappa}_{jki}-\bar{\kappa}_{jk,i})\bar{\kappa}^{i,u}Z_u-Y_{jk}\right)}}+\sum_{s=1}^{p}\bar{\kappa}_{r,s}\times V_n(r,s)$$

$$+\sum_{s=1}^{p}\bar{\kappa}^{s,t}\times n\tilde{O}_d^2(n^{-1})+n^{-1/2}\sum_{s=1}^{p}\overline{\overline{\bar{\kappa}^{s,j}\bar{\kappa}^{t,k}\left((\bar{\kappa}_{jki}-\bar{\kappa}_{jk,i})\bar{\kappa}^{i,u}Z_u-Y_{jk}\right)}}\times n\tilde{O}_d^2(n^{-1})$$

$$+\sum_{s=1}^{p}V_n(r,s)\times n\tilde{O}_d^2(n^{-1})$$

(Expand $Z_{rs}$ by definition in the third term)

$$=\overline{\overline{\bar{\kappa}_{r,s}\bar{\kappa}^{s,t}}}+n^{-1/2}\overline{\overline{\bar{\kappa}^{t,k}\left((\bar{\kappa}_{rki}-\bar{\kappa}_{rk,i})\bar{\kappa}^{i,u}Z_u-Y_{rk}\right)}}-n^{-1/2}\overline{\overline{\bar{\kappa}^{s,t}\left((\bar{\kappa}_{rsi}-\bar{\kappa}_{rs,i})\bar{\kappa}^{i,u}Z_u-Y_{rs}\right)}}$$

$$+n^{-1}\sum_{s=1}^{p}\left(Z_{rs}-\overline{\overline{\bar{\kappa}_{rsi}\bar{\kappa}^{i,u}Z_u}}\right)\times\overline{\overline{\bar{\kappa}^{s,j}\bar{\kappa}^{t,k}\left((\bar{\kappa}_{jki}-\bar{\kappa}_{jk,i})\bar{\kappa}^{i,u}Z_u-Y_{jk}\right)}}$$

$$+n^{-1/2}\sum_{s=1}^{p}\left(Z_{rs}-\overline{\overline{\bar{\kappa}_{rsi}\bar{\kappa}^{i,u}Z_u}}\right)\times V_n(r,s)$$

$$+\sum_{s=1}^{p}\bar{\kappa}_{r,s}\times V_n(r,s)+\sum_{s=1}^{p}\bar{\kappa}^{s,t}\times n\tilde{O}_d^2(n^{-1})+\sum_{s=1}^{p}V_n(r,s)\times n\tilde{O}_d^2(n^{-1})$$

$$+n^{-1/2}\sum_{s=1}^{p}\overline{\overline{\bar{\kappa}^{s,j}\bar{\kappa}^{t,k}\left((\bar{\kappa}_{jki}-\bar{\kappa}_{jk,i})\bar{\kappa}^{i,u}Z_u-Y_{jk}\right)}}\times n\tilde{O}_d^2(n^{-1})$$

(Rearrange the third term to find out that it can cancel out with the second term)

$$=\overline{\overline{\bar{\kappa}_{r,s}\bar{\kappa}^{s,t}}}+n^{-1}\sum_{s=1}^{p}\left(Z_{rs}-\overline{\overline{\bar{\kappa}_{rsi}\bar{\kappa}^{i,u}Z_u}}\right)\times\overline{\overline{\bar{\kappa}^{s,j}\bar{\kappa}^{t,k}\left((\bar{\kappa}_{jki}-\bar{\kappa}_{jk,i})\bar{\kappa}^{i,u}Z_u-Y_{jk}\right)}}$$

$$+n^{-1/2}\sum_{s=1}^{p}\left(Z_{rs}-\overline{\overline{\bar{\kappa}_{rsi}\bar{\kappa}^{i,u}Z_u}}\right)\times V_n(r,s)+\sum_{s=1}^{p}\bar{\kappa}_{r,s}\times V_n(r,s)+\sum_{s=1}^{p}\bar{\kappa}^{s,t}\times n\tilde{O}_d^2(n^{-1})$$

$$+n^{-1/2}\sum_{s=1}^{p}\overline{\overline{\bar{\kappa}^{s,j}\bar{\kappa}^{t,k}\left((\bar{\kappa}_{jki}-\bar{\kappa}_{jk,i})\bar{\kappa}^{i,u}Z_u-Y_{jk}\right)}}\times n\tilde{O}_d^2(n^{-1})+\sum_{s=1}^{p}V_n(r,s)\times n\tilde{O}_d^2(n^{-1}).$$

In order for the above expression to equal the right hand side of (2.26), we must have the collection of terms other than $\overline{\overline{\bar{\kappa}_{r,s}\bar{\kappa}^{s,t}}}$ satisfy :

$$n^{-1}\sum_{s=1}^{p}\left(Z_{rs}-\overline{\overline{\bar{\kappa}_{rsi}\bar{\kappa}^{i,u}Z_u}}\right)\times\overline{\overline{\bar{\kappa}^{s,j}\bar{\kappa}^{t,k}\left((\bar{\kappa}_{jki}-\bar{\kappa}_{jk,i})\bar{\kappa}^{i,u}Z_u-Y_{jk}\right)}}$$

$$+n^{-1/2}\sum_{s=1}^{p}\left(Z_{rs}-\overline{\overline{\bar{\kappa}_{rsi}\bar{\kappa}^{i,u}Z_u}}\right)\times V_n(r,s)+\sum_{s=1}^{p}\bar{\kappa}_{r,s}\times V_n(r,s)$$

$$+n^{-1/2}\sum_{s=1}^{p}\overline{\overline{\bar{\kappa}^{s,j}\bar{\kappa}^{t,k}\left((\bar{\kappa}_{jki}-\bar{\kappa}_{jk,i})\bar{\kappa}^{i,u}Z_u-Y_{jk}\right)}}\times n\tilde{O}_d^2(n^{-1})$$

$$+\sum_{s=1}^{p}\bar{\kappa}^{s,t}\times n\tilde{O}_d^2(n^{-1})+\sum_{s=1}^{p}V_n(r,s)\times n\tilde{O}_d^2(n^{-1})=0\,. \tag{2.27}$$

Grouping the left hand side of (2.27) by $V_n(s, t)$, we have:

$$\sum_{s=1}^{p} V_n(s,t)\times\left(n^{-1/2}\left(Z_{rs}-\overline{\overline{\kappa}_{rsi}\overline{\kappa}^{i,u}Z_u}\right)+\overline{\kappa}_{r,s}+n\tilde{O}_d^2(n^{-1})\right)$$

$$=-n^{-1/2}\sum_{s=1}^{p}\overline{\overline{\kappa}^{s,j}\overline{\kappa}^{t,k}\left((\overline{\kappa}_{jki}-\overline{\kappa}_{jk,i})\overline{\kappa}^{i,u}Z_u-Y_{jk}\right)}\times n\tilde{O}_d^2(n^{-1})-\sum_{s=1}^{p}\overline{\kappa}^{s,t}\times n\tilde{O}_d^2(n^{-1})$$

$$-n^{-1}\sum_{s=1}^{p}\left(Z_{rs}-\overline{\overline{\kappa}_{rsi}\overline{\kappa}^{i,u}Z_u}\right)\times\overline{\overline{\kappa}^{s,j}\overline{\kappa}^{t,k}\left((\overline{\kappa}_{jki}-\overline{\kappa}_{jk,i})\overline{\kappa}^{i,u}Z_u-Y_{jk}\right)}\ . \tag{2.28}$$

For the left hand side of (2.28), $n^{-1/2}\left(Z_{rs}-\overline{\overline{\kappa}_{rsi}\overline{\kappa}^{i,u}Z_u}\right)+n\tilde{O}_d^2(n^{-1})=o_p(1)$ follows from Slutsky's theorem. For the right hand side of (2.28), by Slutsky's theorem, we have $-n^{-1/2}\sum_{s=1}^{p}\overline{\overline{\kappa}^{s,j}\overline{\kappa}^{t,k}\left((\overline{\kappa}_{jki}-\overline{\kappa}_{jk,i})\overline{\kappa}^{i,u}Z_u-Y_{jk}\right)}\times n\tilde{O}_d^2(n^{-1})=o_p(n^{-1})$. The second term can be simplified as $-\sum_{s=1}^{p}\overline{\kappa}^{s,t}\times n\tilde{O}_d^2(n^{-1})=n\tilde{O}_d^2(n^{-1})$ by the fact that $\overline{\kappa}^{s,t}=O(1)$ followed from condition A3. The last term of the right hand side of (2.28) $-n^{-1}\sum_{s=1}^{p}\left(Z_{rs}-\overline{\kappa_{rsi}\kappa^{i,u}Z_u}\right)\times\overline{\kappa^{s,j}\kappa^{t,k}\left((\kappa_{jki}-\kappa_{jk,i})\kappa^{i,u}Z_u-Y_{jk}\right)}=n\tilde{O}_d^2(n^{-1})$. As a result, (2.28) can be rewritten as

$$\sum_{s=1}^{p}V_n(s,t)\times\left(\overline{\kappa}_{r,s}+o_p(1)\right)=n\tilde{O}_d^2(n^{-1})\,. \tag{2.29}$$

Equation (2.29) holds for $r, t = 1, 2, \ldots, p$. A matrix form representation is as follows:

$$\boldsymbol{V}_n \times \left( \bar{\boldsymbol{F}}_n(\boldsymbol{\theta}^*) + \boldsymbol{o}_p(1) \right) = n\tilde{\boldsymbol{O}}_d^2(n^{-1}), \tag{2.30}$$

where $\boldsymbol{V}_n \equiv \{V_n(r, s)\}_{r,\, s\, =\, 1,\, 2,\, \ldots,\, p}$, $\bar{\boldsymbol{F}}_n(\boldsymbol{\theta}^*) \equiv \{\kappa_{r,\, s}\}_{r,\, s\, =\, 1,\, 2,\, \ldots,\, p}$, $\boldsymbol{o}_p(1)$ is a $p \times p$ matrix with each element being $o_p(1)$, and $\tilde{\boldsymbol{O}}_d^2(n^{-1})$ is a $p \times p$ matrix with each element being $\tilde{O}_d^2(n^{-1})$. By condition A3, $\bar{\boldsymbol{F}}_n(\boldsymbol{\theta}^*)$ converges to an invertible constant matrix $\boldsymbol{F}(\boldsymbol{\theta}^*)$, which implies that $\left( \bar{\boldsymbol{F}}_n(\boldsymbol{\theta}^*) + \boldsymbol{o}_p(1) \right)^{-1}$ converges in probability to $\boldsymbol{F}(\boldsymbol{\theta}^*)^{-1}$. By the continuous mapping theorem (Mann and Wald, 1943), we have the following:

$$\begin{aligned} \boldsymbol{V}_n &= n\tilde{\boldsymbol{O}}_d^2(n^{-1}) \times \left( \boldsymbol{F}(\boldsymbol{\theta}^*)^{-1} + \boldsymbol{o}_p(1) \right) \\ &= n\tilde{\boldsymbol{O}}_d^2(n^{-1}) \times \boldsymbol{F}(\boldsymbol{\theta}^*)^{-1} + n\tilde{\boldsymbol{O}}_d^2(n^{-1}) \times \boldsymbol{o}_p(1). \end{aligned}$$

As a result, each element in $\boldsymbol{V}_n$ is a linear combination of $n\tilde{O}_d^2(n^{-1})$ terms plus some higher order terms which are dominated by $n\tilde{O}_d^2(n^{-1})$. Thus, $V_n(r, s) = n\tilde{O}_d^2(n^{-1})$, for $r$, $s = 1, 2, \ldots, p$. Consequently,

$$\bar{\boldsymbol{H}}_n(\hat{\boldsymbol{\theta}}_n)^{-1}(r,s) = \bar{\kappa}^{r,s} + n^{-1/2}\overline{\overline{\bar{\kappa}^{r,j}\bar{\kappa}^{s,k}\left( (\bar{\kappa}_{jki} - \bar{\kappa}_{jk,i})\bar{\kappa}^{i,u} Z_u - Y_{jk} \right)}} + n\tilde{O}_d^2(n^{-1}). \qquad \square$$

**Lemma 4**

For i.n.i.d sample data with conditions A1–A8 in Section 2.2, the covariance between $\hat{t}_{nr}$ and $\hat{t}_{ns}$, for any *r*, *s*, can be expressed as:

$$\operatorname{cov}(\hat{t}_{nr}, \hat{t}_{ns}) = n^{-1}\bar{\kappa}^{r,s} + o(n^{-1}).$$

*Proof*: By Lemma 1 and the definition of covariance,

$$
\begin{aligned}
\operatorname{cov}(\hat{t}_{nr}, \hat{t}_{ns}) &= \operatorname{cov}\Big((\hat{t}_{nr} - t_r^*), (\hat{t}_{ns} - t_s^*)\Big) \\
&= \operatorname{cov}\Big( (-n^{-1/2}\overline{\overline{\bar{\kappa}^{r,u} Z_u}} + n\tilde{O}_d^2(n^{-1})), (-n^{-1/2}\overline{\overline{\bar{\kappa}^{s,v} Z_v}} + n\tilde{O}_d^2(n^{-1}))\Big) \\
&= E\Big( (-n^{-1/2}\overline{\overline{\bar{\kappa}^{r,u} Z_u}} + n\tilde{O}_d^2(n^{-1})) \times (-n^{-1/2}\overline{\overline{\bar{\kappa}^{s,v} Z_v}} + n\tilde{O}_d^2(n^{-1}))\Big) \\
&\quad - E\Big( -n^{-1/2}\overline{\overline{\bar{\kappa}^{r,u} Z_u}} + n\tilde{O}_d^2(n^{-1})\Big) \times E\Big( -n^{-1/2}\overline{\overline{\bar{\kappa}^{s,v} Z_v}} + n\tilde{O}_d^2(n^{-1})\Big) \\
&= E\Big( n^{-1}\overline{\overline{\bar{\kappa}^{r,u} Z_u}} \times \overline{\overline{\bar{\kappa}^{s,v} Z_v}}\Big) - n^{-1/2} E\Big[\Big(\overline{\overline{\bar{\kappa}^{r,u} Z_u}} + \overline{\overline{\bar{\kappa}^{s,v} Z_v}}\Big) \times n\tilde{O}_d^2(n^{-1})\Big] \\
&\quad + E\Big[ n\tilde{O}_d^2(n^{-1}) \times n\tilde{O}_d^2(n^{-1})\Big] - E\Big(n\tilde{O}_d^2(n^{-1})\Big) \times E\Big(n\tilde{O}_d^2(n^{-1})\Big), \qquad (2.31)
\end{aligned}
$$

where the last equation follows from the fact that $E(\bar{\kappa}^{r,u} Z_u) = \bar{\kappa}^{r,u} E(Z_u) = 0$ for *r*, *u* = 1, 2, …, *p*. The first term after the last equality in (2.31) can be rewritten as $E\left(n^{-1}\overline{\overline{\bar{\kappa}^{r,u}Z_u}}\times\overline{\overline{\bar{\kappa}^{s,v}Z_v}}\right) = E\left(n^{-1}\overline{\overline{\bar{\kappa}^{r,u}\bar{\kappa}^{s,v}Z_uZ_v}}\right) = n^{-1}\overline{\overline{\bar{\kappa}^{r,u}\bar{\kappa}^{s,v}E(Z_uZ_v)}} = n^{-1}\overline{\overline{\bar{\kappa}^{r,u}\bar{\kappa}^{s,v}\bar{\kappa}_{u,v}}}$ $= n^{-1}\bar{\kappa}^{r,s}$. For the second term in the last equality of (2.31), with condition A6, we know that by the WLLN for i.n.i.d data, $n^{-1/2}Z_u \to 0$ in probability. Thus, by Slutsky's theorem, $n^{-1/2}\left(\overline{\overline{\bar{\kappa}^{r,u}Z_u}} + \overline{\overline{\bar{\kappa}^{s,u}Z_u}}\right) \times n\tilde{O}_d^2(n^{-1}) = o_p(n^{-1})$ and by condition A7 and the DCT, the second term in (2.31) is $o(n^{-1})$. Similarly, for the third term in the last equality of (2.31), $n \times n\tilde{O}_d^2(n^{-1}) \times n\tilde{O}_d^2(n^{-1}) \to 0$ in probability, which follows from Slutsky's theorem. Again by condition A7, the third term in (2.31) is $o(n^{-1})$. Now, we want to show the last term in the last equality of (2.31) is also $o(n^{-1})$. It suffices to show that $nE\left(n\tilde{O}_d^2(n^{-1})\right)$ $\times E\left(n\tilde{O}_d^2(n^{-1})\right) = E\left(\sqrt{n}\times n\tilde{O}_d^2(n^{-1})\right) \times E\left(\sqrt{n}\times n\tilde{O}_d^2(n^{-1})\right) \to 0$. As a matter of fact, $\sqrt{n}\times n\tilde{O}_d^2(n^{-1}) \to 0$ in probability by Slutsky's theorem. Condition A7 implies $E\left(\sqrt{n}\times n\tilde{O}_d^2(n^{-1})\right) \to 0$, and thus $nE\left(n\tilde{O}_d^2(n^{-1})\right) \times E\left(n\tilde{O}_d^2(n^{-1})\right) \to 0$, i.e. the last term in (2.31) is $o(n^{-1})$. Consequently, the covariance between $\hat{t}_{nr}$ and $\hat{t}_{ns}$, for any $r$, $s$, can be expressed as:

$$\operatorname{cov}(\hat{t}_{nr}, \hat{t}_{ns}) = n^{-1}\bar{\kappa}^{r,s} + o(n^{-1}). \qquad \square$$

**Lemma 5**

If we define

$$A_{nrs} \equiv n^{-1/2} \overline{\overline{\bar{\kappa}^{r,t} \bar{\kappa}^{s,u} Y_{tu}}}, \tag{2.32}$$

$$B_{nrs} \equiv n^{-1/2} \overline{\overline{\bar{\kappa}^{r,t} \bar{\kappa}^{s,u} (\bar{\kappa}_{tuv} - \bar{\kappa}_{tu,v}) \bar{\kappa}^{v,w} Z_w}}, \tag{2.33}$$

then the following are true:

**(a).** $E(A_{nrs}B_{nrs}) = 0$.

**(b).** $E\left((A_{nrs} + B_{nrs}) \times n\tilde{O}_d^2(n^{-1})\right) = o(n^{-1})$.

**(c).** $E\left(o(1) \times \left(-A_{nrs} + n\tilde{O}_d^2(n^{-1})\right)\right) = o(n^{-1})$.

*Proof*:

**(a).** By definition, $E(Y_{tu}) = 0$ and $E(Z_w) = 0$, which implies that $E(A_{nrs}) = E(B_{nrs}) = 0$; since $\mathrm{cov}(Z_r, Y_{st}) = 0$ for $r, s, t = 1, \ldots, p$, $E(A_{nrs}B_{nrs}) = \mathrm{cov}(A_{nrs}, B_{nrs}) = 0$.

**(b).** Rewrite $A_{nrs}$ as:

$$\begin{aligned}
A_{nrs} &= n^{-1/2} \overline{\overline{\bar{\kappa}^{r,t} \bar{\kappa}^{s,u} Y_{tu}}} \\
&= n^{-1/2} \overline{\overline{\bar{\kappa}^{r,t} \bar{\kappa}^{s,u} (Z_{tu} - \bar{\kappa}_{tu,v} \bar{\kappa}^{v,w} Z_w)}} \\
&= n^{-1} \overline{\overline{\bar{\kappa}^{r,t} \bar{\kappa}^{s,u} (U_{tu} - n\bar{\kappa}_{tu} - \bar{\kappa}_{tu,v} \bar{\kappa}^{v,w} U_w)}}
\end{aligned}$$

$$= n^{-1}\sum_{i=1}^{n} \overline{\overline{\bar{\kappa}^{r,t}\bar{\kappa}^{s,u}(U_{tu}^{i} - \kappa_{tu}^{i} - \bar{\kappa}_{tu,v}\bar{\kappa}^{v,w}U_{w}^{i})}}\,. \tag{2.34}$$

Equation in (2.34) reveals that $A_{nrs}$ is a sample mean of a sequence of independent random variables, each of which has mean zero. By condition A6, we know that the WLLN for i.n.i.d sample implies that $A_{nrs} \rightarrow 0$ in probability. Similarly,

$$B_{nrs} = n^{-1/2}\overline{\overline{\bar{\kappa}^{r,t}\bar{\kappa}^{s,u}(\bar{\kappa}_{tuv} - \bar{\kappa}_{tu,v})\bar{\kappa}^{v,w}Z_{w}}}$$

$$= n^{-1}\overline{\overline{\bar{\kappa}^{r,t}\bar{\kappa}^{s,u}(\bar{\kappa}_{tuv} - \bar{\kappa}_{tu,v})\bar{\kappa}^{v,w}U_{w}}}$$

$$= n^{-1}\sum\nolimits_{i=1}^{n} \overline{\overline{\bar{\kappa}^{r,t}\bar{\kappa}^{s,u}(\bar{\kappa}_{tuv} - \bar{\kappa}_{tu,v})\bar{\kappa}^{v,w}U_{w}^{i}}}\,. \tag{2.35}$$

Equation (2.35) indicates that $B_{nrs}$ is a sample mean of a sequence of independent random variables, each of which has mean zero. By condition A6 and the WLLN for i.n.i.d sample, $B_{nrs} \rightarrow 0$ in probability. Thus, $(A_{nrs} + B_{nrs}) \times n\tilde{O}_{d}^{2}(n^{-1}) = o_{p}(n^{-1})$ by Slutsky's theorem. As a result, by condition A7 and the DCT, $E[(A_{nrs}+B_{nrs}) \times n\tilde{O}_{d}^{2}(n^{-1})] = o(n^{-1})$ .

**(c).** From (2.34), we know that $A_{nrs}$ is a sample mean of a sequence of independent random variables, each of which has mean zero. Thus, $E(A_{nrs}) = 0$, implying $E(o(1) \times A_{nrs}) = o(1) \times E(A_{nrs}) = 0$. Furthermore, $o(1) \times n\tilde{O}_{d}^{2}(n^{-1}) = o_{p}(n^{-1})$ by Slutsky's theorem. Thus $E\left(o(1) \times n\tilde{O}_{d}^{2}(n^{-1})\right) = o(n^{-1})$ by condition A7 and the DCT. □

## 2.4 Main results

In this section, we present results that show the advantage of $\bar{\boldsymbol{F}}_n(\hat{\boldsymbol{\theta}}_n)^{-1}$ over $\bar{\boldsymbol{H}}_n(\hat{\boldsymbol{\theta}}_n)^{-1}$ in estimating $n\,\text{cov}(\hat{\boldsymbol{\theta}}_n)$, the scaled covariance matrix of $\hat{\boldsymbol{\theta}}_n$. In our scheme, we compare the two matrices for an arbitrary corresponding entry. Specifically, we show that asymptotically, $\bar{\boldsymbol{F}}_n(\hat{\boldsymbol{\theta}}_n)^{-1}(r,s)$ estimates $n\,\text{cov}(\hat{t}_{nr},\hat{t}_{ns})$ at least as well as $\bar{\boldsymbol{H}}_n(\hat{\boldsymbol{\theta}}_n)^{-1}(r,s)$ under the mean squared error criterion for all *r*, *s* = 1, 2, …, *p*. Hence, in a limit sense, $\bar{\boldsymbol{F}}_n(\hat{\boldsymbol{\theta}}_n)^{-1}$ is preferred to $\bar{\boldsymbol{H}}_n(\hat{\boldsymbol{\theta}}_n)^{-1}$ in estimating $n\,\text{cov}(\hat{\boldsymbol{\theta}}_n)$.

There are degenerate cases where $\bar{\boldsymbol{F}}_n(\hat{\boldsymbol{\theta}}_n) = \bar{\boldsymbol{H}}_n(\hat{\boldsymbol{\theta}}_n)$ for all *n*, and thus the equal performance of the two estimates: $\bar{\boldsymbol{F}}_n(\hat{\boldsymbol{\theta}}_n)^{-1}$ and $\bar{\boldsymbol{H}}_n(\hat{\boldsymbol{\theta}}_n)^{-1}$. The following lemma demonstrates situations when $\bar{F}_n(\hat{\theta}_n) = \bar{H}_n(\hat{\theta}_n)$ for one-parameter i.i.d. exponential families. Please note that such situations do not satisfy condition A9 in Section 2.2.

**Lemma 6**

If $\boldsymbol{X} = [X_1, X_2, \ldots, X_n]$ is a sequence of i.i.d scalar random variables whose density belongs to the one-parameter exponential family, i.e., for $i = 1, 2, \ldots, n$,

$$p_i(x_i, \theta) = h(x_i)\exp\{\eta(\theta)T(x_i) - A(\theta)\},$$

where θ is a scalar parameter, $h(\cdot)$, $T(\cdot)$, $\eta(\cdot)$, and $A(\cdot)$ are known functions. Then $\bar{F}_n(\hat{\theta}_n) = \bar{H}_n(\hat{\theta}_n)$ if and only if

$$\eta''(\hat{\theta}_n)\left\{E(T(X_1))\Big|_{\theta=\hat{\theta}_n} - n^{-1}\sum_{i=1}^{n} T(x_i)\right\} = 0\,,$$

where $\eta''(\hat{\theta}_n)$ denotes the second derivative of $\eta(\theta)$ with respect to θ evaluated at $\theta = \hat{\theta}_n$.

*Remarks*: Conditions above hold for the following situations:

I. $T(x) = x$, $E(T(X_1))\big|_{\theta=\hat{\theta}_n} = \hat{\theta}_n$, and $\hat{\theta}_n = \bar{X}$, where $\bar{X}$ denotes the sample mean. Examples that satisfy these conditions include the Poisson distribution, binomial distribution, and normal distribution with unknown mean.

II. $\eta''(\hat{\theta}_n) = 0$. An example that satisfies this condition is $p_i(x_i, \theta) = \exp\{-\theta x_i + \log\theta\}$.

*Proof*: The negative log-likelihood function of $\boldsymbol{X}$ is

$$l(\theta, \boldsymbol{x}) = -\sum_{i=1}^{n}\left[\log(h(x_i)) + \eta(\theta)T(x_i) - A(\theta)\right].$$

The second derivative of $l(\theta, \boldsymbol{x})$ with respect to θ is

$$\frac{\partial^2 l(\theta, \boldsymbol{x})}{\partial\theta^2} = -\eta''(\theta)\sum_{i=1}^{n} T(x_i) + nA''(\theta)\,.$$

Thus,

$$\bar{F}_n(\hat{\theta}_n) = -n^{-1}\eta''(\hat{\theta}_n)\sum_{i=1}^{n} E(T(X_i))\Big|_{\theta=\hat{\theta}_n} + A''(\hat{\theta}_n)\,,$$

$$\bar{H}_n(\hat{\theta}_n) = -n^{-1}\eta''(\hat{\theta}_n)\sum_{i=1}^{n} T(x_i) + A''(\hat{\theta}_n)\,,$$

and

$$\bar{F}_n(\hat{\theta}_n) - \bar{H}_n(\hat{\theta}_n) = -n^{-1}\eta''(\hat{\theta}_n)\left\{\sum_{i=1}^{n} E(T(X_i))\Big|_{\theta=\hat{\theta}_n} - \sum_{i=1}^{n} T(X_i)\right\}$$

$$= -\eta''(\hat{\theta}_n)\left\{E(T(X_i))\Big|_{\theta=\hat{\theta}_n} - n^{-1}\sum_{i=1}^{n} T(X_i)\right\},$$

where the second equality follows from the fact that $X_1, X_2, \ldots, X_n$ are i.i.d. As a result, $\bar{F}_n(\hat{\theta}_n) = \bar{H}_n(\hat{\theta}_n)$ if and only if

$$\eta''(\hat{\theta}_n)\left\{E(T(X_1))\Big|_{\theta=\hat{\theta}_n} - n^{-1}\sum_{i=1}^{n} T(x_i)\right\} = 0\,. \qquad \square$$

Now let us present the main theorem which shows the superiority of $\bar{\boldsymbol{F}}_n(\hat{\boldsymbol{\theta}}_n)^{-1}$ over $\bar{\boldsymbol{H}}_n(\hat{\boldsymbol{\theta}}_n)^{-1}$ in estimating $n\,\mathrm{cov}(\hat{\boldsymbol{\theta}}_n)$.

**Theorem 1**

Under conditions A1–A8 in Section 2.2, for every pair $(r,s)$, $r$, $s = 1, 2, \ldots, p$

$$\liminf_{n\to\infty} \frac{E\left[\left(\bar{\boldsymbol{H}}_n(\hat{\boldsymbol{\theta}}_n)^{-1}(r,s) - n\,\mathrm{cov}(\hat{t}_{nr}, \hat{t}_{ns})\right)^2\right]}{E\left[\left(\bar{\boldsymbol{F}}_n(\hat{\boldsymbol{\theta}}_n)^{-1}(r,s) - n\,\mathrm{cov}(\hat{t}_{nr}, \hat{t}_{ns})\right)^2\right]} \geq 1. \tag{2.36}$$

Furthermore, if condition A9 in Section 2.2 is satisfied, then the strict inequality (> 1) in (2.36) holds.

*Remark*: The inequality in (2.36) indicates that, asymptotically, $\bar{\boldsymbol{F}}_n(\hat{\boldsymbol{\theta}}_n)^{-1}$ produces no greater mean squared error than $\bar{\boldsymbol{H}}_n(\hat{\boldsymbol{\theta}}_n)^{-1}$ in estimating $n\,\mathrm{cov}(\hat{\boldsymbol{\theta}}_n)$ at each element level. In addition, if the difference between $\bar{\boldsymbol{F}}_n(\hat{\boldsymbol{\theta}}_n)^{-1}$ and $\bar{\boldsymbol{H}}_n(\hat{\boldsymbol{\theta}}_n)^{-1}$ is significant enough (see condition A9 and the corresponding comments in Section 2.2), then the strict inequality in (2.36) holds, i.e., $\bar{\boldsymbol{F}}_n(\hat{\boldsymbol{\theta}}_n)^{-1}$ produces strictly smaller mean squared error than $\bar{\boldsymbol{H}}_n(\hat{\boldsymbol{\theta}}_n)^{-1}$ asymptotically. Please note that condition A9 requires that the referred function of the first and the second derivative of the log-likelihood has variance strictly positive asymptotically. This is common in situations where $\bar{\boldsymbol{F}}_n(\hat{\boldsymbol{\theta}}_n)^{-1}$ and $\bar{\boldsymbol{H}}_n(\hat{\boldsymbol{\theta}}_n)^{-1}$ are unequal for all $n$. Condition A9 is general enough to allow for other settings, as well, given its requirement that $\bar{\boldsymbol{F}}_n(\hat{\boldsymbol{\theta}}_n)^{-1}$ and $\bar{\boldsymbol{H}}_n(\hat{\boldsymbol{\theta}}_n)^{-1}$ be non-identical on only a subsequence.

*Proof*: By Lemmas 2–4, we derive the following decomposition:

$$\left[\bar{\boldsymbol{H}}_n(\hat{\boldsymbol{\theta}}_n)^{-1}(r,s) - n\operatorname{cov}(\hat{t}_{nr},\hat{t}_{ns})\right]^2 - \left[\bar{\boldsymbol{F}}_n(\hat{\boldsymbol{\theta}}_n)^{-1}(r,s) - n\operatorname{cov}(\hat{t}_{nr},\hat{t}_{ns})\right]^2$$

$$= \left[\bar{\boldsymbol{F}}_n(\hat{\boldsymbol{\theta}}_n)^{-1}(r,s) + \bar{\boldsymbol{H}}_n(\hat{\boldsymbol{\theta}}_n)^{-1}(r,s) - 2n\operatorname{cov}(\hat{t}_{nr},\hat{t}_{ns})\right] \times \left[\bar{\boldsymbol{H}}_n(\hat{\boldsymbol{\theta}}_n)^{-1}(r,s) - \bar{\boldsymbol{F}}_n(\hat{\boldsymbol{\theta}}_n)^{-1}(r,s)\right]$$

$$= \left[ n^{-1/2}\overline{\overline{\bar{\kappa}^{r,t}\bar{\kappa}^{s,u}\left(2(\bar{\kappa}_{tuv} - \bar{\kappa}_{tu,v})\bar{\kappa}^{v,w}Z_w - Y_{tu}\right)}} + n\tilde{O}_d^2(n^{-1}) + o(1) \right]$$

$$\times \left[ -n^{-1/2}\overline{\overline{\bar{\kappa}^{r,t}\bar{\kappa}^{s,u}Y_{tu}}} + n\tilde{O}_d^2(n^{-1}) \right]$$

$$= \left(2B_{nrs} - A_{nrs} + n\tilde{O}_d^2(n^{-1}) + o(1)\right) \times \left(-A_{nrs} + n\tilde{O}_d^2(n^{-1})\right)$$

$$= A_{nrs}^2 - 2A_{nrs}B_{nrs} - (A_{nrs} + B_{nrs}) \times n\tilde{O}_d^2(n^{-1}) + o_p(n^{-1})$$

$$+ o(1) \times \left(-A_{nrs} + n\tilde{O}_d^2(n^{-1})\right), \tag{2.37}$$

where $A_{nrs}$ and $B_{nrs}$ are as defined in (2.32) and (2.33). The $o_p(n^{-1})$ term in the last equality of (2.37) follows from a product of two $n\tilde{O}_d^2(n^{-1})$ terms, each of which converges in probability to zero when multiplied by $\sqrt{n}$, which is implied by Slutsky's theorem. Consequently, $n \times n\tilde{O}_d^2(n^{-1}) \times n\tilde{O}_d^2(n^{-1})$ converges in probability to zero by Slutsky's theorem, indicating $n\tilde{O}_d^2(n^{-1}) \times n\tilde{O}_d^2(n^{-1}) = o_p(n^{-1})$. Taking expectation of both sides of (2.37), we now have

$$d_n(r,s) \equiv E\left[\left(\bar{\boldsymbol{H}}_n(\hat{\boldsymbol{\theta}}_n)^{-1}(r,s) - n\operatorname{cov}(\hat{t}_{nr},\hat{t}_{ns})\right)^2\right] - E\left[\left(\bar{\boldsymbol{F}}_n(\hat{\boldsymbol{\theta}}_n)^{-1}(r,s) - n\operatorname{cov}(\hat{t}_{nr},\hat{t}_{ns})\right)^2\right]$$

$$= E\left(A_{nrs}^2 - 2A_{nrs}B_{nrs} - (A_{nrs}+B_{nrs}) \times n\tilde{O}_d^2(n^{-1}) + o_p(n^{-1})\right)$$

$$+ E\left(o(1) \times \left(-A_{n,r,s} + n\tilde{O}_d^2(n^{-1})\right)\right)$$

$$= E(A_{nrs}^2) + o(n^{-1}), \tag{2.38}$$

where the last equality follows from Lemma 5. Consequently,

$$\liminf_{n\to\infty} \frac{E\left[\left(\bar{\boldsymbol{H}}_n(\hat{\boldsymbol{\theta}}_n)^{-1}(r,s) - n\operatorname{cov}(\hat{t}_{nr},\hat{t}_{ns})\right)^2\right]}{E\left[\left(\bar{\boldsymbol{F}}_n(\hat{\boldsymbol{\theta}}_n)^{-1}(r,s) - n\operatorname{cov}(\hat{t}_{nr},\hat{t}_{ns})\right)^2\right]}$$

$$= \liminf_{n\to\infty} \frac{E\left[\left(\bar{\boldsymbol{H}}_n(\hat{\boldsymbol{\theta}}_n)^{-1}(r,s) - n\operatorname{cov}(\hat{t}_{nr},\hat{t}_{ns})\right)^2\right] - E\left[\left(\bar{\boldsymbol{F}}_n(\hat{\boldsymbol{\theta}}_n)^{-1}(r,s) - n\operatorname{cov}(\hat{t}_{nr},\hat{t}_{ns})\right)^2\right]}{E\left[\left(\bar{\boldsymbol{F}}_n(\hat{\boldsymbol{\theta}}_n)^{-1}(r,s) - n\operatorname{cov}(\hat{t}_{nr},\hat{t}_{ns})\right)^2\right]} + 1$$

$$= \liminf_{n\to\infty} \frac{nd_n(r,s)}{nE\left[\left(\bar{\boldsymbol{F}}_n(\hat{\boldsymbol{\theta}}_n)^{-1}(r,s) - n\operatorname{cov}(\hat{t}_{nr},\hat{t}_{ns})\right)^2\right]} + 1$$

$$\geq 1, \tag{2.39}$$

where the last inequality follows from (2.38) and the fact that $\liminf_{n\to\infty} nd_n(r,s)$ $= \liminf_{n\to\infty}(nE(A_{nrs}^2) + o(1)) \geq 0$.

Now we want to demonstrate that if, in addition, condition A9 is satisfied, the strict inequality in (2.39) holds. In fact, we have

$$\liminf_{n\to\infty}\left[nE(A_{nrs}^2)\right]$$

$$=\liminf_{n\to\infty}\left(n\,\mathrm{var}(A_{nrs})\right)$$

$$=\liminf_{n\to\infty}\; n\,\mathrm{var}\left\{n^{-1}\sum_{i=1}^{n}\left[\overline{\overline{\kappa}^{r,t}\overline{\kappa}^{s,u}\left(U_{tu}^{i}-\kappa_{tu}^{i}-\overline{\kappa}_{tu,v}\overline{\kappa}^{v,w}U_{w}^{i}\right)}\right]\right\}$$

$$=\liminf_{n\to\infty}\left\{n^{-1}\sum_{i=1}^{n}\mathrm{var}\left[\overline{\overline{\kappa}^{r,t}\overline{\kappa}^{s,u}\left(U_{tu}^{i}-\kappa_{tu}^{i}-\overline{\kappa}_{tu,v}\overline{\kappa}^{v,w}U_{w}^{i}\right)}\right]\right\}$$

$$>0, \tag{2.40}$$

where the first equality follows from the fact that $E(A_{nrs}) = 0$; the third equality is due to the fact that observations across $i$ are independent; and the inequality at the end follows from condition A9. Consequently,

$$\liminf_{n\to\infty} nd_n(r,s)$$

$$=\liminf_{n\to\infty}\left(nE(A_{nrs}^2)+n\times o(n^{-1})\right)$$

$$=\liminf_{n\to\infty} nE(A_{nrs}^2)+\liminf_{n\to\infty} n\times o(n^{-1})$$

$$>0,$$

where the last inequality follows from (2.40) and the fact that $\liminf_{n\to\infty} n \times o(n^{-1})$ $= \lim_{n\to\infty} n \times o(n^{-1})$ = 0. Furthermore, we have

$$\liminf_{n\to\infty} \frac{E\left[\left(\bar{\boldsymbol{H}}_n(\hat{\boldsymbol{\theta}}_n)^{-1}(r,s) - n\,\mathrm{cov}(\hat{t}_{nr},\hat{t}_{ns})\right)^2\right]}{E\left[\left(\bar{\boldsymbol{F}}_n(\hat{\boldsymbol{\theta}}_n)^{-1}(r,s) - n\,\mathrm{cov}(\hat{t}_{nr},\hat{t}_{ns})\right)^2\right]}$$

$$= \liminf_{n\to\infty} \frac{n d_n(r,s)}{nE\left[\left(\bar{\boldsymbol{F}}_n(\hat{\boldsymbol{\theta}}_n)^{-1}(r,s) - n\,\mathrm{cov}(\hat{t}_{nr},\hat{t}_{ns})\right)^2\right]} + 1$$

$$>1. \qquad \square$$

In summary, Theorem 1 indicates that, asymptotically, $\bar{\boldsymbol{F}}_n(\hat{\boldsymbol{\theta}}_n)^{-1}$ performs at least as well as $\bar{\boldsymbol{H}}_n(\hat{\boldsymbol{\theta}}_n)^{-1}$ in estimating the scaled covariance matrix of $\hat{\boldsymbol{\theta}}_n$ for each matrix entry. An immediate practical problem is that in many situations, the closed analytical form of the Fisher information is not attainable (e.g. Example 1 in Section 3.1). Given the relation between expected and observed Fisher information in (1.2), one way to get around with this issue is to use numerical approximations. A few Monte Carlo-based techniques are introduced in Appendix B, which include a basic resampling method, a feedback-based method, and an independent perturbation per measurement method.

# Chapter 3

# Numerical Studies

In this chapter, we show three numerical examples to demonstrate the superiority of the expected FIM in estimating the covariance matrix of MLEs. The first example considers a mixture Gaussian model, which is popularly used in practice to deal with statistical populations with two or more subpopulations. The second example covers a signal-plus-noise situation, which is commonly seen in practical problems where statistical inferences are made in the presence of noise. The last example discusses a linear discrete-time state-space model, which has wide applications in areas such as engineering, economics, and finance.

Before we present the examples, let us introduce the common notation that is used throughout all three cases. To compare the performance of expected and observed FIM in estimating the covariance matrix of MLEs, we define discrepancy matrices $\boldsymbol{M_H}$ and $\boldsymbol{M_F}$ such that the $(r, s)$ entry of $\boldsymbol{M_H}$ is $\boldsymbol{M_H}(r, s) \equiv E\left[\left(\bar{\boldsymbol{H}}_n(\hat{\boldsymbol{\theta}}_n)^{-1}(r,s) - n\operatorname{cov}(\hat{t}_{nr}, \hat{t}_{ns})\right)\right]^2$ and the $(r, s)$ entry of $\boldsymbol{M_F}$ is $\boldsymbol{M_F}(r, s) \equiv E\left[\left(\bar{\boldsymbol{F}}_n(\hat{\boldsymbol{\theta}}_n)^{-1}(r,s) - n\,\text{cov}(\hat{t}_{nr}, \hat{t}_{ns})\right)\right]^2$. Correspondingly, we use $\boldsymbol{R_H}$ and $\boldsymbol{R_F}$ to denote matrices composed of relative square root of mean squared error for each component, i.e., $\boldsymbol{R_H}(r, s) \equiv \left|\sqrt{\boldsymbol{M_H}(r,s)} \big/ n\,\text{cov}(\hat{t}_{nr}, \hat{t}_{ns})\right|$ and $\boldsymbol{R_F}(r, s) \equiv \left|\sqrt{\boldsymbol{M_F}(r,s)} \big/ n\,\text{cov}(\hat{t}_{nr}, \hat{t}_{ns})\right|$. Notice that the performance is assessed at a component level, which is consistent with the approach used in Theorem 1 of Chapter 2. However, in the examples that follow, we are not able to provide true $\boldsymbol{M_H}$ (or $\boldsymbol{M_F}$) or $\boldsymbol{R_H}$ (or $\boldsymbol{R_F}$) because closed forms of the expectations are not attainable. We present numerical estimates as replacements, which are derived from an average of a large number of sample values.

For each example, we also show a typical value of both $\bar{\boldsymbol{F}}_n(\hat{\boldsymbol{\theta}}_n)^{-1}$ and $\bar{\boldsymbol{H}}_n(\hat{\boldsymbol{\theta}}_n)^{-1}$. We first generate 1001 independent values of $\bar{\boldsymbol{F}}_n(\hat{\boldsymbol{\theta}}_n)^{-1}$ or $\bar{\boldsymbol{H}}_n(\hat{\boldsymbol{\theta}}_n)^{-1}$. We then rank the 1001 matrices by their (Euclidean) distance to the true (or approximated) $n\,\text{cov}(\hat{\boldsymbol{\theta}}_n)$. The outcome with the median distance from $n\,\text{cov}(\hat{\boldsymbol{\theta}}_n)$ is picked as the typical outcome.

### 3.1 Example 1—Mixture Gaussian distribution

The mixture Gaussian distribution is of great interest and is popularly used in practical applications (Wang, 2001; Stein et al., 2002). In this study, we consider a mixture of two univariate Gaussian distributions. Specifically, let $\boldsymbol{X} = [x_1, x_2, \ldots, x_n]^T$ be an i.i.d. sequence with probability density function:

$$f(x,\boldsymbol{\theta}) = \lambda \exp(-(x-\mu_1)^2/(2\sigma_1^2))\Big/\sqrt{2\pi\sigma_1^2} + (1-\lambda)\exp(-(x-\mu_2)^2/(2\sigma_2^2))\Big/\sqrt{2\pi\sigma_2^2} \ ,$$

where $\boldsymbol{\theta} = [\lambda, \mu_1, \mu_2]^T$ and $\sigma_1$, $\sigma_2$ are known. There is no closed form for MLE in this case. We use Newton's method to achieve numerical approximation of $\hat{\boldsymbol{\theta}}_n$. The covariance matrix of $\hat{\boldsymbol{\theta}}_n$ is approximated by the sample covariance of $10^6$ values of $\hat{\boldsymbol{\theta}}_n$ from $10^6$ independent realizations of data. This is a good approximation of the true covariance matrix because the first four post-decimal digits do not change as the number of independent realizations increases. The analytical form of the true FIM is not attainable. But the closed form of the Hessian matrix is computable (see Boldea and Magnus, 2009). We approximate the true FIM using the sample average of the Hessian matrix over $10^5$ independent replications. This is a good approximation since the first four post-decimal digits do not vary as the amount of averaging increases beyond $10^5$.

In this study, we consider two cases when $\boldsymbol{\theta}^* = [0.5, 0, 4]^T$ with $n = 50$, and $\boldsymbol{\theta}^* = [0.5, 0, 2]^T$ with $n = 100$, where for both cases $\sigma_1 = \sigma_2 = 1$. For the second case, we use a bigger sample size $n$ to allow for adequate information to achieve reliable MLE when two individual Gaussian distributions have closer mean. We estimate $\boldsymbol{M_H}$ and $\boldsymbol{M_F}$ by sample averages over $10^5$ independent replications. This is a good approximation of the mean squared error matrix because the first three post-decimal digits do not change as the number of independent realizations increases. Simulation results are summarized in Table 3.1.

**Table 3.1**: Simulation results for Example 1. The scaled covariance matrix $n\,\text{cov}(\hat{\boldsymbol{\theta}}_n)$ is approximated by the sample covariance matrix of $10^6$ values of $\hat{\boldsymbol{\theta}}_n$ from $10^6$ independent realizations; $\boldsymbol{M_H}$ and $\boldsymbol{M_F}$ are approximated by sample averages over $10^5$ independent replications.

| | $\boldsymbol{\theta} = [0.5, 0, 4]$, $n = 50$ | $\boldsymbol{\theta} = [0.5, 0, 2]$, $n = 100$ |
|---|---|---|
| $n\,\text{cov}(\hat{\boldsymbol{\theta}}_n)$ | $\begin{bmatrix} 0.2719 & 0.1151 & 0.1036 \\ 0.1151 & 2.4006 & 0.4333 \\ 0.1036 & 0.4333 & 2.5389 \end{bmatrix}$ | $\begin{bmatrix} 1.3881 & 2.4472 & 2.4351 \\ 2.4472 & 7.7186 & 4.7643 \\ 2.4351 & 4.7643 & 7.7076 \end{bmatrix}$ |
| Typical $\bar{\boldsymbol{H}}_n(\hat{\boldsymbol{\theta}}_n)^{-1}$ | $\begin{bmatrix} 0.2601 & 0.0518 & 0.0464 \\ 0.0518 & 2.4353 & 0.1796 \\ 0.0464 & 0.1796 & 2.0442 \end{bmatrix}$ | $\begin{bmatrix} 0.7831 & 1.6690 & 1.0829 \\ 1.6690 & 8.3483 & 3.0017 \\ 1.0829 & 3.0017 & 3.7717 \end{bmatrix}$ |
| Typical $\bar{\boldsymbol{F}}_n(\hat{\boldsymbol{\theta}}_n)^{-1}$ | $\begin{bmatrix} 0.2762 & 0.1065 & 0.1056 \\ 0.1065 & 2.4619 & 0.4019 \\ 0.1056 & 0.4019 & 2.4314 \end{bmatrix}$ | $\begin{bmatrix} 1.3006 & 1.9168 & 1.8297 \\ 1.9168 & 5.8289 & 2.8067 \\ 1.8297 & 2.8067 & 5.9017 \end{bmatrix}$ |
| $\boldsymbol{M_H}$ | $\begin{bmatrix} 0.0020 & 0.0033 & 0.0054 \\ 0.0033 & 0.2637 & 0.0657 \\ 0.0054 & 0.0657 & 0.3703 \end{bmatrix}$ | $\begin{bmatrix} 1.0903 & 2.6794 & 4.0365 \\ 2.6794 & 15.1482 & 5.6708 \\ 4.0365 & 5.6708 & 9.7639 \end{bmatrix}$ |
| $\boldsymbol{M_F}$ | $\begin{bmatrix} 0.0020 & 0.0002 & 0.0002 \\ 0.0002 & 0.0081 & 0.0023 \\ 0.0002 & 0.0023 & 0.0082 \end{bmatrix}$ | $\begin{bmatrix} 0.1015 & 0.9894 & 0.4956 \\ 0.9894 & 10.1784 & 3.4294 \\ 0.4956 & 3.4294 & 5.2064 \end{bmatrix}$ |
| $\boldsymbol{M_H} - \boldsymbol{M_F}$ | $\begin{bmatrix} -5.1\times 10^{-6} & 0.0031 & 0.0052 \\ 0.0031 & 0.2556 & 0.0634 \\ 0.0052 & 0.0634 & 0.3621 \end{bmatrix}$ | $\begin{bmatrix} 0.9888 & 1.6900 & 3.5409 \\ 1.6900 & 4.9698 & 2.2414 \\ 3.5409 & 2.2414 & 4.5575 \end{bmatrix}$ |
| $\boldsymbol{R_H}$ | $\begin{bmatrix} 0.1640 & 0.5021 & 0.7104 \\ 0.5021 & 0.2139 & 0.5915 \\ 0.7104 & 0.5915 & 0.2396 \end{bmatrix}$ | $\begin{bmatrix} 0.7521 & 0.6688 & 0.8250 \\ 0.6688 & 0.5042 & 0.4998 \\ 0.8250 & 0.4998 & 0.4045 \end{bmatrix}$ |
| $\boldsymbol{R_F}$ | $\begin{bmatrix} 0.1643 & 0.1086 & 0.1505 \\ 0.1086 & 0.0375 & 0.1114 \\ 0.1505 & 0.1114 & 0.0356 \end{bmatrix}$ | $\begin{bmatrix} 0.2294 & 0.4064 & 0.2891 \\ 0.4064 & 0.4133 & 0.3886 \\ 0.2891 & 0.3886 & 0.2960 \end{bmatrix}$ |

Results in Table 3.1 are consistent with theoretical conclusion in Chapter 2. For $\boldsymbol{\theta}^*$ = $[0.5, 0, 2]^T$, every entry of $\bar{\boldsymbol{F}}_n(\hat{\boldsymbol{\theta}}_n)^{-1}$ has a lower MSE in estimating the corresponding component in $n\,\mathrm{cov}(\hat{\boldsymbol{\theta}}_n)$ than $\bar{\boldsymbol{H}}_n(\hat{\boldsymbol{\theta}}_n)^{-1}$. The difference in MSEs is quite significant. For $\boldsymbol{\theta}^*$ = $[0.5, 0, 4]^T$, almost every entry of $\bar{\boldsymbol{F}}_n(\hat{\boldsymbol{\theta}}_n)^{-1}$ has a lower MSE than the corresponding entry in $\bar{\boldsymbol{H}}_n(\hat{\boldsymbol{\theta}}_n)^{-1}$ except for the (1, 1) entry. However, the difference in MSEs for the (1, 1) entry is very small and we believe this might correspond to the higher order term that we ignore in the theoretical discussion. Moreover, for both values of $\boldsymbol{\theta}^*$, typical values of $\bar{\boldsymbol{F}}_n(\hat{\boldsymbol{\theta}}_n)^{-1}$ are closer to $n\,\mathrm{cov}(\hat{\boldsymbol{\theta}}_n)$ than the typical values of $\bar{\boldsymbol{H}}_n(\hat{\boldsymbol{\theta}}_n)^{-1}$. Specifically, typical values of $\bar{\boldsymbol{F}}_n(\hat{\boldsymbol{\theta}}_n)^{-1}$ produce smaller sum of squared errors over all entries than typical values of $\bar{\boldsymbol{H}}_n(\hat{\boldsymbol{\theta}}_n)^{-1}$.

### 3.2 Example 2—Signal-plus-noise problem

The signal-plus-noise situation represents a class of common problems in practice. Examples of application for this statistical model include estimation of the initial mean vector and covariance matrix in a state-space (Kalman filter) model from a cross-section of realizations (Shumway et al., 1981), dose response analysis (Hui and Berger, 1983), estimation of parameters for random-coefficient linear models (Sun, 1982), small area estimating in survey sampling (Ghosha and Rao, 1994), sensitivity studies (Spall, 1985a; Spall and Chin, 1990), and nuisance parameter analysis (Spall, 1989).

This study is a generalization of the numerical study in Cao and Spall (2009), where a scalar case of $\boldsymbol{\theta}$ is discussed. Consider a sequence of i.n.i.d random vectors $\boldsymbol{X}_1$, $\boldsymbol{X}_2$, …, $\boldsymbol{X}_n$. For each $i \in \{1, 2, \ldots, n\}$, $\boldsymbol{X}_i$ is multivariate normal:

$$\boldsymbol{X}_i \sim N(\boldsymbol{\mu}, \boldsymbol{\Sigma} + \boldsymbol{Q}_i),$$

where $\boldsymbol{\mu}$ is the common mean vector across observations, $\boldsymbol{\Sigma}$ is the common part of the covariance matrices and $\boldsymbol{Q}_i$ is the covariance matrix of noise for observation $i$. In practice, the $\boldsymbol{Q}_i$ are known and $\boldsymbol{\theta}$ contains unique elements in $\boldsymbol{\mu}$ and $\boldsymbol{\Sigma}$.

There are no closed forms for $\hat{\boldsymbol{\theta}}_n$ or its covariance matrix. We use Newton's method to find a numerical approximation of $\hat{\boldsymbol{\theta}}_n$ and estimate the covariance matrix based on $10^6$ MLEs from $10^6$ independent realizations. This is a good approximation of the true covariance matrix because the first four post-decimal digits do not change as the number of independent realizations increases. Closed forms of FIM for this signal-plus-noise Gaussian model are provided in Shumway (1982) and closed forms of the corresponding Hessian matrix are provided in Spall (1985b). Spall (2003) shows the same for the special case of scalar data.

In this study, we consider 4-dimensional $\boldsymbol{X}_i$ and diagonal $\boldsymbol{\Sigma}$: $\boldsymbol{\Sigma} = \text{diag}\{\Sigma_{11}, \Sigma_{22}, \Sigma_{33}, \Sigma_{44}\}$. Thus, $\boldsymbol{\theta} = [\mu_1, \mu_2, \mu_3, \mu_4, \Sigma_{11}, \Sigma_{22}, \Sigma_{33}, \Sigma_{44}]^T$ is an $8 \times 1$ vector. The underlying true value of the parameters in this study is $\boldsymbol{\theta}^* = [0, 0, 0, 0, 1, 1, 1, 1]^T$. The known $\boldsymbol{Q}_i$ matrices are constructed in the following way: let $\boldsymbol{U}$ be a 4 × 4 deterministic matrix where each entry is drawn from uniform (0, 0.1)-distribution. $\boldsymbol{Q}_i$ is defined as $\boldsymbol{Q}_i = \sqrt{i} \times \boldsymbol{U}^T\boldsymbol{U}$.

In our study, we use the following $\boldsymbol{U}^T\boldsymbol{U}$ matrix:

$$\boldsymbol{U}^T\boldsymbol{U} = \begin{bmatrix} 0.0289 & 0.0219 & 0.0120 & 0.0216 \\ 0.0219 & 0.0200 & 0.0068 & 0.0189 \\ 0.0120 & 0.0068 & 0.0076 & 0.0053 \\ 0.0216 & 0.0189 & 0.0053 & 0.0210 \end{bmatrix}$$

The sample size in this study is $n = 80$. We estimate $\boldsymbol{M_H}$ and $\boldsymbol{M_F}$ by sample averages over $10^5$ independent replications. This is a good approximation of the mean squared error matrix because the first three post-decimal digits do not change as the number of independent realizations increases. Simulation results are summarized in Table 3.2.

**Table 3.2**: Simulation results for Example 2. The scaled covariance matrix $n\,\mathrm{cov}(\hat{\boldsymbol{\theta}}_n)$ is approximated by the sample covariance matrix of $10^6$ MLEs from $10^6$ independent realizations. Both $\boldsymbol{M_H}$ and $\boldsymbol{M_F}$ are approximated by sample averages over $10^5$ independent replications.

| | |
|---|---|
| $n\,\mathrm{cov}(\hat{\boldsymbol{\theta}}_n)$ | $\begin{bmatrix} 0.46 & 0.07 & 0.06 & 0.05 & 0.17 & -0.08 & -0.00 & -0.08 \\ 0.07 & 0.33 & 0.07 & 0.04 & -0.23 & 0.27 & -0.06 & -0.12 \\ 0.06 & 0.07 & 0.40 & 0.07 & -0.29 & -0.08 & -0.01 & -0.05 \\ 0.05 & 0.04 & 0.07 & 0.32 & -0.24 & -0.10 & -0.06 & 0.33 \\ 0.17 & -0.23 & -0.29 & -0.24 & 5.32 & 0.81 & 0.60 & 0.61 \\ -0.08 & 0.27 & -0.08 & -0.10 & 0.81 & 2.14 & 0.49 & 0.43 \\ -0.00 & -0.06 & -0.01 & -0.06 & 0.60 & 0.49 & 0.73 & 0.52 \\ -0.08 & -0.12 & -0.05 & 0.33 & 0.61 & 0.43 & 0.52 & 2.27 \end{bmatrix}$ |
| Typical $\bar{\boldsymbol{H}}_n(\hat{\boldsymbol{\theta}}_n)^{-1}$ | $\begin{bmatrix} 1.66 & -0.24 & -0.28 & -0.25 & -0.58 & -0.20 & -0.20 & -0.33 \\ -0.24 & 0.85 & -0.52 & -0.49 & -1.19 & -0.20 & -0.30 & -0.47 \\ -0.28 & -0.52 & 0.73 & -0.54 & -1.14 & -0.23 & -0.28 & -0.39 \\ -0.25 & -0.49 & -0.54 & 0.77 & -1.22 & -0.29 & -0.31 & -0.37 \\ -0.58 & -1.19 & -1.14 & -1.22 & 8.84 & -2.85 & -3.02 & -4.15 \\ -0.20 & -0.20 & -0.23 & -0.29 & -2.85 & 3.10 & -0.66 & -0.93 \\ -0.20 & -0.30 & -0.28 & -0.31 & -3.02 & -0.66 & 2.76 & -1.16 \\ -0.33 & -0.47 & -0.39 & -0.37 & -4.15 & -0.93 & -1.16 & 3.85 \end{bmatrix}$ |
| Typical $\bar{\boldsymbol{F}}_n(\hat{\boldsymbol{\theta}}_n)^{-1}$ | $\begin{bmatrix} 1.87 & 0.12 & 0.07 & 0.12 & 0 & 0 & 0 & 0 \\ 0.12 & 1.43 & 0.03 & 0.11 & 0 & 0 & 0 & 0 \\ 0.07 & 0.03 & 1.32 & 0.03 & 0 & 0 & 0 & 0 \\ 0.12 & 0.11 & 0.03 & 1.43 & 0 & 0 & 0 & 0 \\ 0 & 0 & 0 & 0 & 6.84 & -0.02 & -0.00 & -0.02 \\ 0 & 0 & 0 & 0 & -0.02 & 4.00 & -0.00 & -0.02 \\ 0 & 0 & 0 & 0 & -0.00 & -0.00 & 3.48 & -0.00 \\ 0 & 0 & 0 & 0 & -0.02 & -0.02 & -0.00 & 4.04 \end{bmatrix}$ |

(Table continues next page)

(Table 3.2, continued)

| | |
|---|---|
| $\boldsymbol{M_H}$ | $\begin{bmatrix} 1.66 & 0.40 & 0.31 & 0.48 & 28.72 & 65.39 & 2.50 & 31.12 \\ 0.40 & 1.47 & 0.57 & 1.29 & 98.65 & 228.58 & 8.61 & 108.54 \\ 0.31 & 0.57 & 0.71 & 0.69 & 36.05 & 81.20 & 3.09 & 38.53 \\ 0.48 & 1.29 & 0.69 & 2.16 & 136.86 & 317.60 & 11.96 & 151.32 \\ 28.72 & 98.65 & 36.05 & 136.86 & 11730.81 & 27047.99 & 1022.28 & 12810.97 \\ 65.39 & 228.58 & 81.20 & 317.60 & 27047.99 & 61770.41 & 2367.40 & 29953.45 \\ 2.50 & 8.61 & 3.09 & 11.96 & 1022.28 & 2367.40 & 95.91 & 1127.44 \\ 31.12 & 108.54 & 38.53 & 151.32 & 12810.97 & 29953.45 & 1127.44 & 14469.97 \end{bmatrix}$ |
| $\boldsymbol{M_F}$ | $\begin{bmatrix} 0.86 & 0.00 & 0.00 & 0.00 & 0.02 & 0.00 & 0.00 & 0.00 \\ 0.00 & 1.41 & 0.00 & 0.00 & 0.05 & 0.07 & 0.00 & 0.01 \\ 0.00 & 0.00 & 0.24 & 0.00 & 0.08 & 0.00 & 0.00 & 0.00 \\ 0.00 & 0.00 & 0.00 & 1.38 & 0.06 & 0.01 & 0.00 & 0.11 \\ 0.02 & 0.05 & 0.08 & 0.06 & 9.92 & 0.70 & 0.37 & 0.41 \\ 0.00 & 0.07 & 0.00 & 0.01 & 0.70 & 12.18 & 0.24 & 0.20 \\ 0.00 & 0.00 & 0.00 & 0.00 & 0.37 & 0.24 & 5.76 & 0.27 \\ 0.00 & 0.01 & 0.00 & 0.11 & 0.41 & 0.20 & 0.27 & 5.49 \end{bmatrix}$ |
| $\boldsymbol{M_H - M_F}$ | $\begin{bmatrix} 0.79 & 0.40 & 0.31 & 0.48 & 28.69 & 65.38 & 2.50 & 31.12 \\ 0.40 & 0.06 & 0.56 & 1.28 & 98.59 & 228.50 & 8.60 & 108.53 \\ 0.31 & 0.56 & 0.46 & 0.69 & 35.97 & 81.19 & 3.09 & 38.53 \\ 0.48 & 1.28 & 0.69 & 0.78 & 136.80 & 317.59 & 11.95 & 151.21 \\ 28.69 & 98.59 & 35.97 & 136.80 & 11720.88 & 27047.28 & 1021.91 & 12810.55 \\ 65.38 & 228.50 & 81.19 & 317.59 & 27047.28 & 61758.22 & 2367.15 & 29953.24 \\ 2.50 & 8.60 & 3.09 & 11.95 & 1021.91 & 2367.15 & 90.14 & 1127.16 \\ 31.12 & 108.53 & 38.53 & 151.21 & 12810.55 & 29953.24 & 1127.16 & 14464.47 \end{bmatrix}$ |
| $\boldsymbol{R_H}$ | $\begin{bmatrix} 2.75 & 8.04 & 8.23 & 11.64 & 31.01 & 100.07 & 1439.00 & 68.70 \\ 8.04 & 3.62 & 10.39 & 23.19 & 42.79 & 54.38 & 44.06 & 82.03 \\ 8.23 & 10.39 & 2.08 & 10.74 & 20.60 & 102.86 & 91.12 & 110.06 \\ 11.64 & 23.19 & 10.74 & 4.53 & 46.96 & 167.33 & 51.38 & 36.90 \\ 31.01 & 42.79 & 20.60 & 46.96 & 20.34 & 202.64 & 53.12 & 182.64 \\ 100.07 & 54.38 & 102.86 & 167.33 & 202.64 & 115.90 & 97.97 & 398.04 \\ 1439.00 & 44.06 & 91.12 & 51.38 & 53.12 & 97.97 & 13.33 & 64.55 \\ 68.70 & 82.03 & 110.06 & 36.90 & 182.64 & 398.04 & 64.55 & 52.93 \end{bmatrix}$ |

(Table continues next page)

(Table 3.2, continued)

| | |
|---|---|
| $\boldsymbol{R_F}$ | $\begin{bmatrix} 1.98 & 0.61 & 0.02 & 1.10 & 1 & 1 & 1 & 1 \\ 0.61 & 3.54 & 0.45 & 1.25 & 1 & 1 & 1 & 1 \\ 0.02 & 0.45 & 1.23 & 0.60 & 1 & 1 & 1 & 1 \\ 1.10 & 1.25 & 0.60 & 3.62 & 1 & 1 & 1 & 1 \\ 1 & 1 & 1 & 1 & 0.59 & 1.03 & 1.01 & 1.04 \\ 1 & 1 & 1 & 1 & 1.03 & 1.62 & 1.00 & 1.04 \\ 1 & 1 & 1 & 1 & 1.01 & 1.00 & 3.26 & 1.00 \\ 1 & 1 & 1 & 1 & 1.04 & 1.04 & 1.00 & 1.03 \end{bmatrix}$ |

Numerical results in Table 3.2 are consistent with the theoretical conclusion in Chapter 2. Every entry of $\bar{\boldsymbol{F}}_n(\hat{\boldsymbol{\theta}}_n)^{-1}$ produces a smaller error in estimating the corresponding component in $n\,\mathrm{cov}(\hat{\boldsymbol{\theta}}_n)$ than $\bar{\boldsymbol{H}}_n(\hat{\boldsymbol{\theta}}_n)^{-1}$. Notice that large numbers appear in the 4 × 4 lower right sub-matrix in $\boldsymbol{M_H}$ and $\boldsymbol{R_H}$. This is due to the fact that enormous values happen in same realizations of the Hessian matrix $\bar{\boldsymbol{H}}_n(\hat{\boldsymbol{\theta}}_n)^{-1}$. However, we do not see large numbers in $\boldsymbol{M_F}$ or $\boldsymbol{R_F}$ due to the expectation effect in $\bar{\boldsymbol{F}}_n(\hat{\boldsymbol{\theta}}_n)^{-1}$, which avoids enormous values. Notice that there are entries in $\boldsymbol{R_F}$ that exactly equal to one. This is due to the fact that the corresponding entries in $\bar{\boldsymbol{F}}_n(\hat{\boldsymbol{\theta}}_n)^{-1}$ are zero, making the relative root of MSE being 100%. Moreover, the typical outcome of $\bar{\boldsymbol{F}}_n(\hat{\boldsymbol{\theta}}_n)^{-1}$ is closer in value to $n\,\mathrm{cov}(\hat{\boldsymbol{\theta}}_n)$ than the typical outcome of $\bar{\boldsymbol{H}}_n(\hat{\boldsymbol{\theta}}_n)^{-1}$. Specifically, the typical value of $\bar{\boldsymbol{F}}_n(\hat{\boldsymbol{\theta}}_n)^{-1}$ produces smaller sum of squared error over all entries than the typical value of $\bar{\boldsymbol{H}}_n(\hat{\boldsymbol{\theta}}_n)^{-1}$.

### 3.3 Example 3—State-space model

A state-space model is a mathematical description of a physical system as a set of input, output and state variables. Before we introduce the state-space model, let us define some necessary notation. Assume that $\boldsymbol{x}_t$ is an unobserved $l$-dimensional state process, $\boldsymbol{A}$ is an $l \times l$ transition matrix, and $\boldsymbol{w}_t$ is a vector of $l$ zero-mean, independent disturbances with covariance matrix $\boldsymbol{Q}$. Let $\boldsymbol{y}_t$ be an observed $m$-dimensional process, $\boldsymbol{C}$ be an $m \times l$ design matrix, and $\boldsymbol{v}_t$ be a vector of $m$ zero-mean, independent disturbances with covariance matrix $\boldsymbol{R}$. The mean and covariance matrix of $\boldsymbol{x}_0$ (the initial $\boldsymbol{x}_t$) are denoted by $\boldsymbol{\mu}$ and $\boldsymbol{\Sigma}$, respectively. It is assumed that $\boldsymbol{\mu}$ and $\boldsymbol{\Sigma}$ are known and that $\boldsymbol{x}_0$, $\boldsymbol{w}_t$, and $\boldsymbol{v}_t$ are mutually independent and multivariate normal.

The state-space model considered is defined by the equations

$$\boldsymbol{x}_t = \boldsymbol{A}\boldsymbol{x}_{t-1} + \boldsymbol{w}_t, \tag{3.1}$$

$$\boldsymbol{y}_t = \boldsymbol{C}\boldsymbol{x}_t + \boldsymbol{v}_t, \tag{3.2}$$

for $t = 1, 2, \ldots, n$ time periods.

In our context, we consider situations where $l = 3$, $m = 1$ and $\boldsymbol{A}$, $\boldsymbol{C}$, $\boldsymbol{R}$ are known. The unknown parameters of interest are the unique elements in diagonal $\boldsymbol{Q}$, i.e., $\boldsymbol{\theta} = [Q_{11}, Q_{22}, Q_{33}]^T$, where

$$\boldsymbol{Q} = \begin{bmatrix} Q_{11} & 0 & 0 \\ 0 & Q_{22} & 0 \\ 0 & 0 & Q_{33} \end{bmatrix}.$$

Given the definition of $\boldsymbol{\theta}$ above, the log-likelihood function $L(\boldsymbol{\theta})$ for the system described in (3.1) and (3.2) is (neglecting constant terms) (Gibson and Ninness, 2005):

$$L(\boldsymbol{\theta}) = -\frac{1}{2}\sum_{t=1}^{n}\log(\boldsymbol{C}\boldsymbol{P}_{t|t-1}\boldsymbol{C}^T + R) - \frac{1}{2}\sum_{t=1}^{n}\left(\boldsymbol{C}\boldsymbol{P}_{t|t-1}\boldsymbol{C}^T + R\right)^{-1}\varepsilon_t^2, \tag{3.3}$$

whose computation requires Kalman Filter equations:

$$\varepsilon_t = y_t - \boldsymbol{C}\hat{\boldsymbol{x}}_{t|t-1}; \tag{3.4}$$

$$\hat{\boldsymbol{x}}_{t|t-1} = \boldsymbol{A}\hat{\boldsymbol{x}}_{t-1|t-1}; \tag{3.5}$$

$$\hat{\boldsymbol{x}}_{t|t} = \hat{\boldsymbol{x}}_{t|t-1} + \boldsymbol{K}_t\varepsilon_t; \tag{3.6}$$

$$\boldsymbol{K}_t = \boldsymbol{P}_{t|t-1}\boldsymbol{C}^T\left(\boldsymbol{C}\boldsymbol{P}_{t|t-1}\boldsymbol{C}^T + R\right)^{-1}; \tag{3.7}$$

$$\boldsymbol{P}_{t|t-1} = \boldsymbol{A}\boldsymbol{P}_{t-1|t-1}\boldsymbol{A}^T + \boldsymbol{Q}; \tag{3.8}$$

$$\boldsymbol{P}_{t|t} = \boldsymbol{P}_{t|t-1} - \boldsymbol{K}_t\boldsymbol{C}\boldsymbol{P}_{t|t-1}. \tag{3.9}$$

Given (3.3)–(3.9), it is not hard to derive the Hessian matrix $\boldsymbol{H}_n(\boldsymbol{\theta})$ in a recursive form. In addition, the FIM $\boldsymbol{F}_n(\boldsymbol{\theta})$ is also attainable for this state-space model (Cavanaugh and Shumway, 1996). The closed form for the MLE $\hat{\boldsymbol{\theta}}_n$ is not available. We use stochastic search method Algorithm B in Spall (2003, pp. 43–45) to approximate $\hat{\boldsymbol{\theta}}_n$.

Our simulation is based on the specific model:

$$\boldsymbol{A} = \begin{bmatrix} 0 & 1 & 0 \\ 0 & 0 & 1 \\ 0.8 & 0.8 & -0.8 \end{bmatrix},$$

$$\boldsymbol{C} = [1 \;\; 0 \;\; 0],$$

$$R = 1,$$

$$\boldsymbol{\mu} = [0 \;\; 0 \;\; 0]^T,$$

$$\boldsymbol{\Sigma} = \begin{bmatrix} 0 & 0 & 0 \\ 0 & 0 & 0 \\ 0 & 0 & 0 \end{bmatrix}.$$

The true input for $\boldsymbol{\theta}$ is $\boldsymbol{\theta}^* = [1, 1, 1]^T$, i.e.,

$$\boldsymbol{Q} = \begin{bmatrix} 1 & 0 & 0 \\ 0 & 1 & 0 \\ 0 & 0 & 1 \end{bmatrix}.$$

The forms of $\boldsymbol{A}$ and $\boldsymbol{C}$ above are chosen according to the process described in Ljung (1999, Chapter 4). The transition matrix $\boldsymbol{A}$ is designed in such a way that the system is identifiable. Given that there is no closed form for its covariance matrix cov($\hat{\boldsymbol{\theta}}_n$), we approximate cov($\hat{\boldsymbol{\theta}}_n$) using the sample covariance matrix of $10^4$ independent estimates of $\hat{\boldsymbol{\theta}}_n$, where each $\hat{\boldsymbol{\theta}}_n$ is computed from a sequence of observations $y_1, y_2, \ldots, y_{100}$. This is a good approximation since the first three post-decimal digits do not change as the amount of averaging increases beyond $10^4$.

In this study, we consider two sample sizes: $n$ = 100 and $n$ = 200. We estimate $\boldsymbol{M_H}$ and $\boldsymbol{M_F}$ by sample averages over $10^4$ independent replications. This is a good approximation since the first three post-decimal digits do not change as the amount of averaging increases beyond $10^4$.The simulation results are summarized in Table 3.3 ($n$ = 100) and Table 3.4 ($n$ = 200).

**Table 3.3**: Simulation results for Example 3 ($n$ = 100). For $n$ = 100, the scaled covariance matrix $n\operatorname{cov}(\hat{\boldsymbol{\theta}}_n)$ is approximated by the sample covariance matrix of $10^4$ values of $\hat{\boldsymbol{\theta}}_n$ from $10^4$ independent realizations. Both $\boldsymbol{M_H}$ and $\boldsymbol{M_F}$ are approximated by sample averages over $10^4$ independent replications.

| | |
|---|---|
| $n\operatorname{cov}(\hat{\boldsymbol{\theta}}_n)$ | $\begin{bmatrix} 51.8934 & -24.4471 & -33.7404 \\ -24.4471 & 59.4544 & -3.36401 \\ -33.7404 & -3.36401 & 63.0565 \end{bmatrix}$ |
| Typical $\bar{\boldsymbol{H}}_n(\hat{\boldsymbol{\theta}}_n)^{-1}$ | $\begin{bmatrix} 26.1023 & -8.78498 & -11.5611 \\ -8.78498 & 10.6455 & -1.77378 \\ -11.5611 & -1.77378 & 20.0403 \end{bmatrix}$ |
| Typical $\bar{\boldsymbol{F}}_n(\hat{\boldsymbol{\theta}}_n)^{-1}$ | $\begin{bmatrix} 54.7463 & -34.3979 & -47.2827 \\ -34.3979 & 68.2098 & -2.11431 \\ -47.2827 & -2.11431 & 68.8055 \end{bmatrix}$ |
| $\boldsymbol{M_H}$ | $\begin{bmatrix} 11810.3615 & 2706.7101 & 20839.1818 \\ 2706.7101 & 11028.3976 & 1326.5203 \\ 20839.1818 & 1326.5203 & 57649.6871 \end{bmatrix}$ |
| $\boldsymbol{M_F}$ | $\begin{bmatrix} 292.8718 & 139.8859 & 195.8054 \\ 139.8859 & 1011.1910 & 4.1963 \\ 195.8054 & 4.1963 & 1371.368 \end{bmatrix}$ |
| $\boldsymbol{M_H} - \boldsymbol{M_F}$ | $\begin{bmatrix} 11517.4897 & 2566.8241 & 20643.3763 \\ 2566.8241 & 10017.2065 & 1322.3240 \\ 20643.3763 & 1322.3240 & 56278.3183 \end{bmatrix}$ |
| $\boldsymbol{R_H}$ | $\begin{bmatrix} 2.0942 & 2.1279 & 4.2784 \\ 2.1279 & 1.7663 & 10.8267 \\ 4.2784 & 10.8267 & 3.8077 \end{bmatrix}$ |
| $\boldsymbol{R_F}$ | $\begin{bmatrix} 0.3297 & 0.4837 & 0.4147 \\ 0.4837 & 0.5348 & 0.6089 \\ 0.4147 & 0.6089 & 0.5872 \end{bmatrix}$ |

**Table 3.4**: Simulation results for Example 3 ($n$ = 200). For $n$ = 200, the scaled covariance matrix $n\,\mathrm{cov}(\hat{\boldsymbol{\theta}}_n)$ is approximated by the sample covariance matrix of $10^4$ values of $\hat{\boldsymbol{\theta}}_n$ from $10^4$ independent realizations. Both $\boldsymbol{M_H}$ and $\boldsymbol{M_F}$ are approximated by sample averages over $10^4$ independent replications.

| | |
|---|---|
| $n\,\mathrm{cov}(\hat{\boldsymbol{\theta}}_n)$ | $\begin{bmatrix} 73.3308 & -34.6841 & -49.5769 \\ -34.6841 & 75.7118 & -12.8851 \\ -49.5769 & -12.8851 & 83.5198 \end{bmatrix}$ |
| Typical $\bar{\boldsymbol{H}}_n(\hat{\boldsymbol{\theta}}_n)^{-1}$ | $\begin{bmatrix} 55.3965 & -20.0538 & -57.3601 \\ -20.0538 & 35.4422 & -9.9754 \\ -57.3601 & -9.9754 & 144.8040 \end{bmatrix}$ |
| Typical $\bar{\boldsymbol{F}}_n(\hat{\boldsymbol{\theta}}_n)^{-1}$ | $\begin{bmatrix} 62.9489 & -19.3356 & -52.5227 \\ -19.3356 & 34.8782 & -11.0883 \\ -52.5227 & -11.0883 & 99.8797 \end{bmatrix}$ |
| $\boldsymbol{M_H}$ | $\begin{bmatrix} 1241.1417 & 712.4586 & 1081.9223 \\ 712.4586 & 3042.7638 & 268.2577 \\ 1081.9223 & 268.2577 & 3752.2690 \end{bmatrix}$ |
| $\boldsymbol{M_F}$ | $\begin{bmatrix} 484.4139 & 166.9961 & 246.2559 \\ 166.9961 & 1101.7465 & 29.8785 \\ 246.2559 & 29.8785 & 1343.5888 \end{bmatrix}$ |
| $\boldsymbol{M_H} - \boldsymbol{M_F}$ | $\begin{bmatrix} 756.7278 & 545.4625 & 835.6663 \\ 545.4625 & 1941.0173 & 238.3791 \\ 835.6663 & 238.3791 & 2408.6802 \end{bmatrix}$ |
| $\boldsymbol{R_H}$ | $\begin{bmatrix} 0.4804 & 0.7695 & 0.6634 \\ 0.7695 & 0.7285 & 1.2711 \\ 0.6634 & 1.2711 & 0.7334 \end{bmatrix}$ |
| $\boldsymbol{R_F}$ | $\begin{bmatrix} 0.3001 & 0.3725 & 0.3165 \\ 0.3725 & 0.4384 & 0.4242 \\ 0.3165 & 0.4242 & 0.4388 \end{bmatrix}$ |

Both Table 3.3 and Table 3.4 show significant advantage of $\bar{\boldsymbol{F}}_n(\hat{\boldsymbol{\theta}}_n)^{-1}$ over $\bar{\boldsymbol{H}}_n(\hat{\boldsymbol{\theta}}_n)^{-1}$ in estimating $n\text{cov}(\hat{\boldsymbol{\theta}}_n)$. For both sample sizes $n$ = 100 and $n$ = 200, $\bar{\boldsymbol{F}}_n(\hat{\boldsymbol{\theta}}_n)^{-1}$ has smaller MSE in estimating the corresponding component in $n\,\text{cov}(\hat{\boldsymbol{\theta}}_n)$ than $\bar{\boldsymbol{H}}_n(\hat{\boldsymbol{\theta}}_n)^{-1}$. Because the observed Fisher information (Hessian) matrix is sample dependent, even one enormous outcome can result in a big MSE for $\bar{\boldsymbol{H}}_n(\hat{\boldsymbol{\theta}}_n)^{-1}$. But this is not the case for the expected FIM due to the averaging effect that is embedded. In both tables above, typical values of $\bar{\boldsymbol{F}}_n(\hat{\boldsymbol{\theta}}_n)^{-1}$ presents a better estimate of $n\,\text{cov}(\hat{\boldsymbol{\theta}}_n)$ than typical values of $\bar{\boldsymbol{H}}_n(\hat{\boldsymbol{\theta}}_n)^{-1}$. Specifically, the typical value of $\bar{\boldsymbol{F}}_n(\hat{\boldsymbol{\theta}}_n)^{-1}$ produces smaller sum of squared error over all entries than the typical value of $\bar{\boldsymbol{H}}_n(\hat{\boldsymbol{\theta}}_n)^{-1}$.

Comparing Table 3.3 and Table 3.4, we find that for the larger $n$, the difference between the MSEs of $\bar{\boldsymbol{F}}_n(\hat{\boldsymbol{\theta}}_n)^{-1}$ and $\bar{\boldsymbol{H}}_n(\hat{\boldsymbol{\theta}}_n)^{-1}$ is smaller. This is not surprising because as sample size grows, $\bar{\boldsymbol{H}}_n(\hat{\boldsymbol{\theta}}_n)^{-1}$ converges to $\bar{\boldsymbol{F}}_n(\hat{\boldsymbol{\theta}}_n)^{-1}$. Furthermore, for the larger sample size, the accuracy of $\bar{\boldsymbol{F}}_n(\hat{\boldsymbol{\theta}}_n)^{-1}$ increases. This makes sense because as we get more information (sample data), we have a better estimate of $n\,\text{cov}(\hat{\boldsymbol{\theta}}_n)$.

# Chapter 4

# Conclusions and Future Work

In Section 4.1, we summarize the research contribution of this dissertation. In Section 4.2, we discuss potential extensions of the work presented in this dissertation. We propose a few approaches which are preliminary ideas but are likely to be explored as future work.

### 4.1 Conclusions

In this dissertation, we compare the relative performance of the expected and observed Fisher information in estimating the covariance matrix of MLE. The discussion throughout this work applies broadly to many contexts with i.n.i.d samples and multi-dimensional parameters of interest. We demonstrate that under a set of reasonable conditions, the inverse expected Fisher information outperforms the inverse observed Fisher information in estimating the covariance matrix of MLE. Specifically, in estimating each entry of the covariance matrix of the MLE, the corresponding entry of the inverse Fisher information (evaluated at the MLE) has no greater mean squared error than the corresponding entry of the observed Fisher information (evaluated at the MLE) in an asymptotic sense, i.e.,

$$\liminf_{n\to\infty} \frac{E\left[\left(\bar{\boldsymbol{H}}_n(\hat{\boldsymbol{\theta}}_n)^{-1}(r,s) - n\operatorname{cov}(\hat{t}_{nr},\hat{t}_{ns})\right)^2\right]}{E\left[\left(\bar{\boldsymbol{F}}_n(\hat{\boldsymbol{\theta}}_n)^{-1}(r,s) - n\operatorname{cov}(\hat{t}_{nr},\hat{t}_{ns})\right)^2\right]} \geq 1,$$

for $r$, $s = 1, 2, \ldots, p$. Note that zero difference in the mean squared errors occurs when the corresponding entries of the inverse expected and the inverse observed Fisher information (evaluated at the MLE) are identical. This can happen even if the two matrices are not identical.

This dissertation provides the theoretical foundation as well as numerical demonstration to support the conclusion above. In Chapter 2, we present detailed theoretical analysis that we developed to reach the final conclusion. All analysis is done at element level, even though the expected and observed Fisher information under consideration are in matrix form. In Chapter 3, three numerical examples are illustrated to support the theoretical conclusion. We first consider an i.i.d mixture Gaussian distribution with three unknown parameters, which is a degenerate case of i.n.i.d samples. The second example demonstrates the theory in a signal-plus-noise situation, where each observation is independent but comes with a different level of noise. The last example considers system identification and parameter estimation in a state-space model, which is of great interest in engineering and other fields. All three examples show the advantage of the expected Fisher information over the observed Fisher information in estimating the covariance matrix of the MLE.

The conclusion of this dissertation provides a theoretical foundation for the choice between expected and observed Fisher information in estimating the covariance matrix of MLE. The development of such foundation has been missing in the literature, though there is great need for constructing accurate approximations to the covariance matrix. Due to the popularity of the MLE as a standard estimation method, people in practice are also interested in the variance/covariance of the MLE. However, there was no solid theory readily available in the literature to provide guideline in choosing a good estimate of the variance/covariance of the MLE. Consequently, people often chose whichever works easier for their problems, regardless of the accuracy of the estimate chosen in estimating the variance/covariance of the MLE. Motivated by the fact that the theoretical foundation for choosing a good estimate for the covariance of MLEs is of great interest in the literature, this dissertation successfully develops theoretical guideline for the choice of a good estimate. We demonstrate that the expected Fisher information performs better under reasonable conditions.

The conclusion of this dissertation may sound contradictory to some known results in the literature. For example, both Efron and Hinkley (1978) and Lindsay and Li (1997) favor the observed Fisher information over the expected Fisher, which the opposite of our conclusion. However, we need to be aware that the context and problem of interest are different in the three cases. In Efron and Hinkley (1978), the variance of the MLE for scalar parameters is discussed in the context of ancillary statistics. Specifically, the problem of interest considers the conditional variance of the MLE given an ancillary statistics. In our discussion, the covariance matrix calculation is in an unconditional setting where no conditional statistics are needed, which is of broader interest in practice. In fact, the reliance on ancillarity imposes a practical limitation on Efron and Hinkley's result. In many situations, ancillarity statistics are difficult to define and in some cases, the definition is not unique. In addition, discussions in Efron and Hinkley (1978) are limited to problems with scalar parameters. And the theoretical analysis in only provided for translation families. Both of these facts impose strong further limits on the practical application of Efron and Hinkley's conclusion. For Lindsay and Li (1997), there is no condition on ancillary statistics and no limitation to scalar parameter and translation families. However, the problem of interest is the realized mean squared error of MLE rather than the covariance. And by definition, the latter equals the expectation of the former. In other words, the estimation target in Lindsay and Li (1997) is an observation-dependent quantity. Thus, it is not surprising that the observed Fisher information is preferred to the expected Fisher in estimating the realized squared error. In contrast, this dissertation considers the unconditioned covariance matrix of MLE for any i.n.i.d observations, which, to our knowledge, has not been discussed in theoretical depth in the literature.

This dissertation includes two appendices. In Appendix A, we discuss the optimal perturbation distribution for small-sample SPSA. We show that if the number of observations is small, the segmented uniform distribution may outperform the asymptotically optimal Bernoulli $\pm 1$ distribution in generating the perturbation vectors for this stochastic algorithm. In Appendix B, Monte-Carlo based approximating techniques are discussed for computing the FIM for complex problems. To elaborate, in the main part of this dissertation, we have shown that under certain conditions, the expected Fisher information is preferred in estimating the covariance matrix of MLE. An immediate practical problem is that in many situations, the closed analytical form of the Fisher information is not attainable (e.g. Example 1 in Section 3.1). Given the relation between expected and observed Fisher information in (1.2), one way to get around with this issue is to use numerical approximations. A few approximating techniques are introduced in Appendix B, which include a basic resampling method, a feedback-based method, and an independent perturbation per measurement method.

### 4.2 Future work

This dissertation has been focusing on comparing the relative performance of two estimates, the inverse expected and inverse observed Fisher information matrix (both evaluated at the MLE), for approximating the covariance matrix of the MLE. It is also of interest to explore other estimates that can possibly obtain better estimation accuracy under different conditions. We introduce a few approaches that we can possibly take to extend the result of this dissertation. Note that these are preliminary thoughts and need to be explored more in the future work.

Although we conclude that under certain conditions/circumstances, the inverse expected Fisher information outperforms the inverse observed Fisher information in estimating the covariance matrix of MLE. This does not imply that the expected Fisher information is the best among *all* estimates. In fact, for some situations, the observed Fisher information or even a mixture of both estimates may be a better estimate. As such, one generalization of the estimation method is to consider linear combinations of the inverse expected Fisher information and the inverse observed Fisher information. This allows for flexibility in constructing a good estimate by assigning appropriate weights to each element (expected Fisher information or observed Fisher information) under various conditions. The problem of interest is now an optimization problem with two scalar variables which are the coefficients of each components of the linear combination. Similar discussion on mixture of expected and observed Fisher information has been seen in Jiang (2005), where the data is generated from mixed linear models.

A more ambitious extension of the problem is to consider all possible functions of the observations as an estimate of the covariance matrix of MLE. In other words, we are interested in extending the discussion to solving a functional optimization problem to find the best estimate of the covariance matrix of MLE. Specifically, we are looking for the solution $T(\boldsymbol{X})$ which solves the optimization problem (1.8), where $T(\boldsymbol{X})$ can be any feasible function of the observations. Here feasibility means that the matrix $T(\boldsymbol{X})$ should be positive semi-definite as an estimate of a covariance matrix. Solving a functional optimization problem is challenging because the dimension of the space of feasible functions is infinite. In other words, any form of function can be applied to the observations as long as the resulting matrix is positive semi-definite. As such, many tools developed in finite-dimensional optimization are not applicable.

Given the challenge of finding the analytical solution of a functional optimization problem, we can start with sub-optimal solutions through approximate functional optimization methods; see Daniel (1971) and Gelfand and Fomin (1963). For example, we can exploit linear approximation schemes based on a certain number of basis functions. Specifically, each entry of $T(\boldsymbol{X})$ can be expressed as a linear combination of a set of basis functions such as polynomial, sines, and cosines, as long as the resulting matrix $T(\boldsymbol{X})$ is positive semi-definite. In such an approach, the original functional optimization problem is reduced to a nonlinear programming problem, where the objective is optimized only through the coefficients of the linear combination. The rationale behind this sub-optimal approach is that when the number of basis functions becomes sufficiently large, the resulting sub-optimal solution should resemble the properties of the optimal solution of the original functional optimization problem. Other than the approaches mentioned above, there are other methods to solve the above functional optimization as well. We have not yet pursued these more general possibilities.

In summary, this dissertation has shown the advantage of the inverse expected FIM over the inverse observed FIM in estimating the covariance matrix of MLEs. In the future work, we may attempt to extend the results by considering other estimation methods besides the inverse expected FIM and the inverse observed FIM. We would first focus on the two possible approaches discussed above in finding the sub-optimal solution of the functional optimization problem. Furthermore, we will explore other possible approaches in solving the functional optimization problem.

# Appendix A

## Non-Bernoulli Perturbation Distributions for Small Samples in Simultaneous Perturbation Stochastic Approximation

Stochastic approximation methods are a family of iterative stochastic optimization algorithms that attempt to find zeroes or extrema of functions which cannot be computed directly, but only estimated via noisy observations. Among various approximation methods, simultaneous perturbation stochastic approximation (SPSA) is a commonly used method because it is easy to implement and it has very nice asymptotic properties. In this appendix, we discuss the optimal distribution for perturbation vectors in SPSA, which is a crucial component of this algorithm. Specifically, we talk about small-sample SPSA, where a limited number of function evaluations are allowed.

## A1. Introduction

Simultaneous perturbation stochastic approximation (SPSA) has proven to be an efficient stochastic approximation approach, see Spall (1992, 1998, and 2009). It has wide applications in areas such as signal processing, system identification and parameter estimation, see www.jhuapl.edu/SPSA/, Bhatnagar (2011), and Spall (2003). The merit of SPSA follows from the construction of the gradient approximation, where only two function evaluations are needed for each step of the gradient approximation regardless of the dimension of the unknown parameter. As a result, SPSA reduces computation demand as compared to the finite difference (FD) method, which requires $2p$ function evaluations to achieve each step of the gradient approximation, where $p$ is the dimension of the problem, see Spall (2003, Chapters 6 and 7). Obviously, the savings in computation with SPSA is more significant as $p$ gets large.

The implementation of SPSA involves perturbation vectors. Typically, the Bernoulli $\pm 1$ distribution is used for the components of the perturbation vectors. This distribution is easy to implement and has been proven asymptotically most efficient, see Sadegh and Spall (1998). As a result, for large-sample SPSA, the Bernoulli distribution is the best choice for the perturbation vectors. However, one might be curious if this optimality remains when only small-sample stochastic approximation (SA) is allowed. Small-sample SA appears commonly in practice where it is expensive, either physically or computationally, to evaluate system performances. For example, it might be very costly to run experiments on a complicated control system. Under such circumstances, a limited number of function evaluations are available for SA. Unlike with large-sample SPSA, one might not be confident that the Bernoulli distribution is still the best candidate for the perturbation vectors in small-sample SPSA.

In this appendix, we discuss the effective perturbation distributions for SPSA with limited samples. Specifically, we consider the segmented uniform (SU) distribution as a representative of non-Bernoulli distributions. The SU distribution has nice properties of easy manipulation both analytically and numerically. For instance, it has both a density function and a distribution function in closed form, making analytical computations possible. Moreover, it does not take much effort to generate SU random variables due to the nature of the SU density, resulting in time-efficient numerical analysis. In our discussion, we focus on one-iteration SPSA, which is a special case of small-sample SPSA. As a finite-sample analogue to asymptotic cases, the one-iteration case is a good starting point as it is easier to analyze and still captures insightful properties of general small-sample SPSA. Along with the analysis of the one-iteration scenario, we gain insights on the behavior of other small samples in the hope that the analysis can be generalized to more than one iteration cases. In fact, we demonstrate numerically that the one-iteration theoretical conclusions do apply to more than one iteration situations.

Discussion and research on non-Bernoulli perturbation distributions in SPSA have been found in the literature, see Bhatnagar et al. (2003) and Hutchison (2002). In Bhatnagar et al. (2003), numerical experiments along with rigorous convergence proofs indicate that deterministic perturbation sequences show promise for significantly faster convergence under certain circumstances; while in Hutchison (2002), conjecture is made based on empirical results that the Bernoulli distribution maintains optimality for small-sample analysis given an optimal choice of parameters. However, no theoretical foundation is provided to validate this conjecture. The application of non-Bernoulli perturbations in SPSA is discussed in Maeda and De Figueiredo (1997) and Spall (2003, Section 7.3).

## A2. Methodology

A2.1 *Problem formulation*

Let $\boldsymbol{\theta} \in \Theta \subseteq R^p$ denote a vector-valued parameter of interest, where $\Theta$ is the parameter space and $p$ is the dimension of $\boldsymbol{\theta}$. Let $L(\boldsymbol{\theta})$ be the loss function, which is observed in the presence of noise: $y(\boldsymbol{\theta}) = L(\boldsymbol{\theta}) + \varepsilon$, where $\varepsilon$ is i.i.d noise, with mean zero and variance $\sigma^2$ ; $y(\boldsymbol{\theta})$ is the observation of $L(\boldsymbol{\theta})$ with noise $\varepsilon$ . The problem is to

$$\min_{\boldsymbol{\theta}\in\Theta} L(\boldsymbol{\theta}) . \tag{A.1}$$

The stochastic optimization algorithm to solve (A.1) is given by the following iterative scheme:

$$\hat{\boldsymbol{\theta}}_{k+1} = \hat{\boldsymbol{\theta}}_k - a_k \hat{\boldsymbol{g}}_k(\hat{\boldsymbol{\theta}}_k) , \tag{A.2}$$

where $\hat{\boldsymbol{\theta}}_k$ is the estimate of $\boldsymbol{\theta}$ at iteration $k$ and $\hat{\boldsymbol{g}}_k(\bullet) \in R^p$ represents an estimate of the gradient of $L$ at iteration $k$. The scalar-valued step-size sequence $\{a_k\}$ is nonnegative, decreasing, and converging to zero. The generic iterative form of (A.2) is analogous to the steepest descent algorithm for deterministic problems.

A2.2 *Perturbation distribution for SPSA*

SPSA uses simultaneous perturbation to estimate the gradient of $L$. The efficiency of this method is that it requires only two function evaluations at each iteration, as compared to $2p$ for the FD method, see Spall (2003, Chapters 6 and 7). Let $\boldsymbol{\Delta}_k$ be a vector of $p$ scalar-valued independent random variables at iteration $k$:

$$\Delta_k = [\Delta_{k1}, \Delta_{k2}, \ldots, \Delta_{kp}]^T .$$

Let $c_k$ be a sequence of positive scalars. The standard simultaneous perturbation form for the gradient estimate is as follows:

$$\hat{\boldsymbol{g}}_k(\hat{\theta}_k) = \begin{bmatrix} \dfrac{y(\hat{\theta}_k + c_k\Delta_k) - y(\hat{\theta}_k - c_k\Delta_k)}{2c_k\Delta_{k1}} \\ \vdots \\ \dfrac{y(\hat{\theta}_k + c_k\Delta_k) - y(\hat{\theta}_k - c_k\Delta_k)}{2c_k\Delta_{kp}} \end{bmatrix}. \tag{A.3}$$

To guarantee the convergence of the algorithm, certain assumptions on $\boldsymbol{\Delta}_k$ should be satisfied:

**I.** $\{\Delta_{ki}\}$ are independent for all $k$, $i$, and identically distributed for all $i$ at each $k$.

**II.** $\{\Delta_{ki}\}$ are symmetrically distributed about zero and uniformly bounded in magnitude for all $k$, $i$.

**III.** $E\left[(y(\hat{\theta}_k \pm c_k \Delta_k)/\Delta_{ki})^2\right]$ is uniformly bounded over $k$ and $i$.

Condition I has an important relationship with the finite inverse moments of the elements of $\mathbf{\Delta}_k$, see Spall (2003, p. 184). An important part of SPSA is the bounded inverse moments condition for the $\Delta_{ki}$ .Valid distributions include the Bernoulli ±1, the segmented uniform, the U-shape distribution and many others, see Spall (2003, p. 185). Two common mean-zero distributions that do not satisfy the bounded inverse moments condition are the symmetric uniform and the mean-zero normal distributions. The failure of both these distributions is a consequence of the amount of probability mass near zero.

In the discussion that follows, we compare the segmented uniform (SU) distribution with the Bernoulli ±1 distribution. To guarantee that the two distributions have the same mean and variance, the domain of SU, following from basic statistics and simple algebra, is given as

$$\left(-(19+3\sqrt{13})/20,\, -(19-3\sqrt{13})/20\right) \cup \left((19-3\sqrt{13})/20,\, (19+3\sqrt{13})/20\right),$$

which is approximately $(-1.4908, -0.4092) \cup (0.4092, 1.4908)$, see Figure A.1. In our analysis, the sequences $\{a_k\}$ and $\{c_k\}$ take standard forms: $a_k = a/(k+2)^{0.602}$, $c_k = c/(k+1)^{0.101}$, where $a$ and $c$ are predetermined constants.

Moments of perturbations under the two distributions are summarized below in Table A.1. These moments will be used in Section A3. Subscripts $i$ and $j$ denote the elements of $\Delta_0$ and $i \neq j$. The derivation follows from basic statistics and simple algebra.

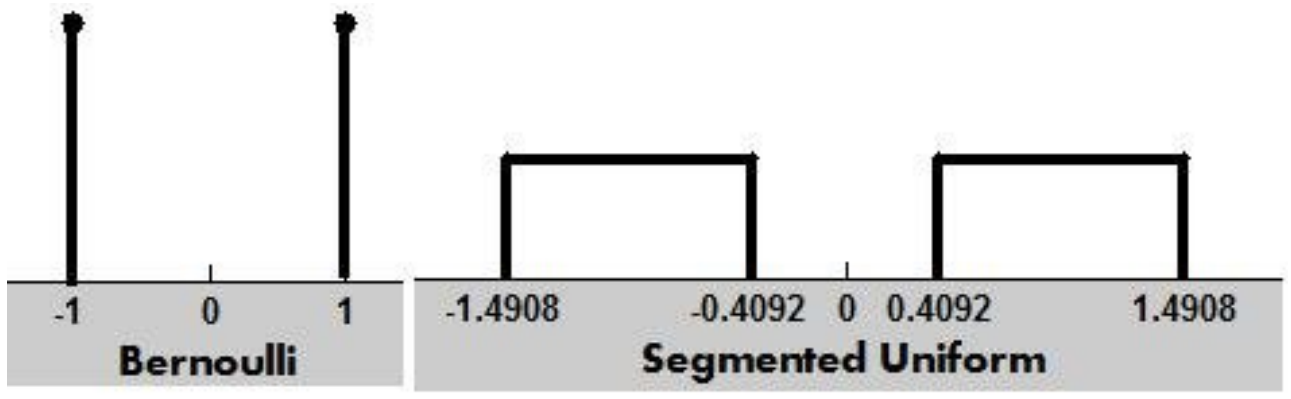

**Figure A.1**: Mass/ probability density functions of the Bernoulli ±1 and the segmenteduniform distributions. Both distributions have mean 0 and variance 1.

**Table A.1**: Moments of perturbations under two distributions

| Expectation | Bernoulli | SU |
|---|---|---|
| $E\left(\Delta_{0i}\right)$ | 0 | 0 |
| $E\left(\Delta_{0i} / \Delta_{0j}\right)$ | 0 | 0 |
| $E\left(\Delta_{0i}^2 / \Delta_{0j}^2\right)$ | 1 | 100/61 |
| $E\left(1 / \Delta_{0i}^2\right)$ | 1 | 100/61 |

**A3. Theoretical analysis**

In this section, we provide conditions under which SU outperforms the Bernoulli distribution. To specifically analyze the development of the algorithm, we consider the extreme example of small samples where only one iteration takes place in SPSA, that is, $k = 1$. We start with this simple case as a basis for possible generalization for larger values of $k$, where the analysis is more complicated. In our analysis, mean squared error (MSE) is used to compare the performance of two distributions.

Before we present the results, let us define necessary notations. Subscripts $S$, $B$ denote SU and the Bernoulli distribution, respectively, e.g. $a_{0S}$ denotes the value of $a_0$ under SU distribution; $L_i$ is the first derivatives of $L$ with respect to the $i$th component of $\boldsymbol{\theta}$, all first derivatives are evaluated at the starting point $\hat{\boldsymbol{\theta}}_0$ ; $\hat{\theta}_{0i}$ and $\theta_i^*$ are the $i$th component of $\hat{\boldsymbol{\theta}}_0$ and $\boldsymbol{\theta}^*$, respectively, where $\boldsymbol{\theta}^*$ is the true value of $\boldsymbol{\theta}$. Following the theorem statement below, we provide some interpretation of the main condition.

**Theorem A.1**

Consider loss function $L(\boldsymbol{\theta})$ with continuous third derivatives. For one iteration of SPSA, the SU distribution produces a smaller MSE between $\hat{\boldsymbol{\theta}}_1$ and $\boldsymbol{\theta}^*$ than the Bernoulli ±1 distribution if the starting point and the relevant coefficients ($a_0$, $c_0$, $\sigma^2$) are such that the following is true:

$$\left[\left(\frac{100}{61}p-\frac{39}{61}\right)a_{0S}^2-pa_{0B}^2\right]\sum_{i=1}^{p}L_i^2$$

$$+(a_{0S}-a_{0B})\left[\frac{p\sigma^2}{2c_{0B}^2}(a_{0S}+a_{0B})-2\sum_{i=1}^{p}(\hat{\theta}_{0i}-\theta_i^*)L_i\right]$$

$$-pa_{0S}^2\sigma^2\left(\frac{1}{2c_{0B}^2}-\frac{50}{61c_{0S}^2}\right)+O(c_0^2)\ <0, \tag{A.4}$$

where the $O(c_0^2)$ term is due to the higher order Taylor expansion.

*Remark A.1*: The choice of the coefficients is not arbitrary. For example, $a_0$ and $c_0$ should be picked according to the standard tuning process, see Spall (2003, Section 7.5); the starting point should stay in a reasonable range given any prior information for the problem. To best use the result of Theorem A.1, one should follow these standards rather than arbitrarily picking the coefficients to make (A.4) true.

*Remark A.2*: If the gains $c_{0S}$ and $c_{0B}$ are small enough such that $O(c_0^2)$ is negligible, the following conditions ((a) and (b)) are sufficient for (A.4) to hold:

(a) The ratios of the gain sequences have the following relations:

$$a_{0S}/a_{0B}<\sqrt{p/(100p/61-39/61)}<1;$$

$$c_{0B}/c_{0S}<\sqrt{61/100}\approx 0.781.$$

(b) In particular, the following inequality is true:

$$2\sum_{i=1}^{p}(\hat{\theta}_{0i}-\theta_i^*)L_i < \frac{p\sigma^2}{2c_{0B}^2}(a_{0S}+a_{0B})\,.$$

The above inequality means that the function is relatively flat and the starting point is not too far away from the true minimum.

*Proof*: By (A.2) and (A.3), the updated estimate of **θ** after one iteration is

$$\hat{\boldsymbol{\theta}}_1 = \hat{\boldsymbol{\theta}}_0 - a_0\frac{y(\hat{\boldsymbol{\theta}}_0 + c_0\boldsymbol{\Delta}_0) - y(\hat{\boldsymbol{\theta}}_0 - c_0\boldsymbol{\Delta}_0)}{2c_0}\times\left[\Delta_{01}^{-1},\ldots,\Delta_{0p}^{-1}\right]^T$$

$$= \hat{\boldsymbol{\theta}}_0 - a_0\frac{L(\hat{\boldsymbol{\theta}}_0 + c_0\boldsymbol{\Delta}_0) + \varepsilon^+ - (L(\hat{\boldsymbol{\theta}}_0 - c_0\boldsymbol{\Delta}_0) + \varepsilon^-)}{2c_0}\times\left[\Delta_{01}^{-1},\ldots,\Delta_{0p}^{-1}\right]^T, \tag{A.5}$$

where $\varepsilon^+$ and $\varepsilon^-$ are the corresponding noise. By a Taylor expansion of the third order,

$$L(\hat{\boldsymbol{\theta}}_0 + c_0\boldsymbol{\Delta}_0) - L(\hat{\boldsymbol{\theta}}_0 - c_0\boldsymbol{\Delta}_0) = 2c_0\sum_{i=1}^{p}L_i\Delta_{0i} + O(c_0^3), \tag{A.6}$$

where the $O(c_0^3)$ term is due to the higher order Taylor expansion. Specifically,

$$O(c_0^3) = \frac{1}{6}c_0^3\sum_{i=1}^{p}\sum_{j=1}^{p}\sum_{k=1}^{p}\left[\left(L_{ijk}(\tilde{\boldsymbol{\theta}}) + L_{ijk}(\tilde{\tilde{\boldsymbol{\theta}}})\right)\Delta_{0i}\Delta_{0j}\Delta_{0k}\right], \tag{A.7}$$

where $L_{ijk}$ denotes the third derivatives of $L$ with respect to the elements $i, j, k$ of $\boldsymbol{\theta}$, $\tilde{\boldsymbol{\theta}}$ and $\tilde{\tilde{\boldsymbol{\theta}}}$ are the intermediate points between $\hat{\boldsymbol{\theta}}_0$ and $\hat{\boldsymbol{\theta}}_0 + c_0\boldsymbol{\Delta}_0$, $\hat{\boldsymbol{\theta}}_0$ and $\hat{\boldsymbol{\theta}}_0 - c_0\boldsymbol{\Delta}_0$, respectively.

Given (A.5), (A.6) and (A.7), and following from algebraic calculation and necessary rearrangements, we compute the difference in MSE $E\left(\left\|\hat{\boldsymbol{\theta}}_1 - \boldsymbol{\theta}^*\right\|^2\right) = E\left(\left(\hat{\boldsymbol{\theta}}_1 - \boldsymbol{\theta}^*\right)^T \left(\hat{\boldsymbol{\theta}}_1 - \boldsymbol{\theta}^*\right)\right)$ under two distributions as follows:

$$E_S\left(\left\|\hat{\boldsymbol{\theta}}_1 - \boldsymbol{\theta}^*\right\|^2\right) - E_B\left(\left\|\hat{\boldsymbol{\theta}}_1 - \boldsymbol{\theta}^*\right\|^2\right)$$

$$= \left[\left(\frac{100}{61}p - \frac{39}{61}\right)a_{0S}^2 - pa_{0B}^2\right]\sum_{i=1}^{p} L_i^2 + (a_{0S} - a_{0B})\left[\frac{p\sigma^2}{2c_{0B}^2}(a_{0S} + a_{0B}) - 2\sum_{i=1}^{p}(\hat{\theta}_{0i} - \theta_i^*)L_i\right]$$

$$-pa_{0S}^2\sigma^2\left(\frac{1}{2c_{0B}^2} - \frac{50}{61c_{0S}^2}\right) + O(c_0^2)\,. \tag{A.8}$$

The derivation of (A.8) involves the computation of relevant moments, which are summarized in Table A.1.

Condition (A.8) in Theorem A.1 may be hard to check for general problems due to the unknown analytical form of the higher order term $O(c_0^2)$. However, if we know more information about the loss function $L$, condition (A.8) can be replaced by a sufficient condition, which is easier to manipulate in practice.

**Corollary A.1**

If we assume an upper bound for the magnitude of the third derivatives of $L$, say, $\left|L_{ijk}(\bullet)\right| \le M$ for all $i, j, k$, where $M$ is a constant, we can establish an upper bound $U$ for the term $O(c_0^2)$ in (A.8), i.e. $O(c_0^2) \le U$. As a result, a more conservative condition for the superiority of SU is

$$\left[\left(\frac{100}{61}p - \frac{39}{61}\right)a_{0S}^2 - pa_{0B}^2\right]\sum_{i=1}^{p} L_i^2 + (a_{0S} - a_{0B})\left[\frac{p\sigma^2}{2c_{0B}^2}(a_{0S} + a_{0B}) - 2\sum_{i=1}^{p}(\hat{\theta}_{0i} - \theta_i^*)L_i\right]$$

$$-pa_{0S}^2\sigma^2\left(\frac{1}{2c_{0B}^2} - \frac{50}{61c_{0S}^2}\right) + U < 0, \tag{A.9}$$

where $U$ is defined as:

$$U = (4a_{0S}c_{0S}^2 + a_{0B}c_{0B}^2)M\left(\sum_{i=1}^{p}\left|\hat{\theta}_{0i} - \theta_i^*\right|\right)(p-1)^2 + \frac{1}{20}a_{0S}^2c_{0S}^4M^2p^7a_{0S}$$

$$+\frac{1}{3}(a_{0S}^2c_{0S}^3 + a_{0B}^2c_{0B}^3)Mp^5 \max_i L_i(\hat{\boldsymbol{\theta}}_0). \tag{A.10}$$

*Proof*: Given (A.7) and the assumption that $\left|L_{ijk}(\bullet)\right| \le M$ for all $i, j, k$, we derive an upper bound $U$ for the term $O(c_0^2)$ as in (A.10). To derive (A.10), we should first find the explicit form of the term $O(c_0^2)$ in (A.8) as follows:

$$O(c^2) = E_S\left[a_{0S}^2 V_S{}^T V_S - 2a_{0S}(\hat{\boldsymbol{\theta}}_0 - \boldsymbol{\theta}^* - a_{0S}W_S)^T V_S\right]$$

$$-E_B\left[a_{0B}^2 V_B{}^T V_B - 2a_{0B}(\hat{\boldsymbol{\theta}}_0 - \boldsymbol{\theta}^* - a_{0B}W_B)^T V_B\right],$$

where for $h = S$ or $B$, as appropriate,

$$W_h = \frac{2c_{0h}\sum_{i=1}^{p} L_i \Delta_{0i} + \varepsilon^+ - \varepsilon^-}{2c_{0h}} \times \left[\Delta_{01}^{-1}, \ldots, \Delta_{0p}^{-1}\right]^T;$$

$$V_h = \frac{1}{12}c_{0h}^2 \sum_{i=1}^{p}\sum_{j=1}^{p}\sum_{k=1}^{p}\left[\left(L_{ijk}(\tilde{\theta}) + L_{ijk}(\tilde{\tilde{\theta}})\right)\Delta_{0i}\Delta_{0j}\Delta_{0k}\right] \times \left[\Delta_{01}^{-1}, \ldots, \Delta_{0p}^{-1}\right]^T.$$

Given the upper bound in (A.10), it follows immediately that (A.9) is a sufficient and more conservative condition for the superiority (smaller MSE) of SU. □

Notice that if $L$ is quadratic, the higher order terms in (A.6) and (A.8) vanish, resulting in the following simpler form of the condition in Theorem A.1.

**Corollary A.2**

For a quadratic loss function $L$, the SU distribution produces a smaller MSE between $\hat{\boldsymbol{\theta}}_1$ and $\boldsymbol{\theta}^*$ than the Bernoulli ±1 distribution for one-iteration SPSA if the following holds:

$$\left[\left(\frac{100}{61}p-\frac{39}{61}\right)a_{0S}^2-pa_{0B}^2\right]\sum_{i=1}^{p}L_i^2+(a_{0S}-a_{0B})\left[\frac{p\sigma^2}{2c_{0B}^2}(a_{0S}+a_{0B})-2\sum_{i=1}^{p}(\hat{\theta}_{0i}-\theta_i^*)L_i\right]$$

$$-pa_{0S}^2\sigma^2\left(\frac{1}{2c_{0B}^2}-\frac{50}{61c_{0S}^2}\right)<0.$$

If $p = 2$, the special form of Corollary A.2 becomes the following, which we use in the numerical example A4.1 below.

**Corollary A.3**

For a quadratic loss function with $p = 2$, the SU distribution produces a smaller MSE between $\hat{\boldsymbol{\theta}}_1$ and $\boldsymbol{\theta}^*$ than the Bernoulli ±1 distribution for one-iteration SPSA if the following holds:

$$\left[\frac{161}{61}a_{0S}^2-2a_{0B}^2\right]\left(L_1^2+L_2^2\right)$$

$$+(a_{0S}-a_{0B})\left[\frac{\sigma^2}{c_{0B}^2}(a_{0S}+a_{0B})-2\left(L_1(\hat{\theta}_{01}-\theta_1^*)+L_2(\hat{\theta}_{02}-\theta_2^*)\right)\right]$$

$$-a_{0S}^2\sigma^2\left(\frac{1}{c_{0B}^2}-\frac{100}{61c_{0S}^2}\right)<0. \tag{A.11}$$

## A4. Numerical examples

A4.1 *Quadratic loss function*

Consider the quadratic loss function $L(\boldsymbol{\theta}) = t_1^2 - t_1 t_2 + t_2^2$, where $\boldsymbol{\theta} = [t_1, t_2]^T$, $\sigma^2 = 1$, $\hat{\theta}_0 = [0.3, 0.3]^T$, $a_S$ = 0.00167, $a_B$ = 0.01897, $c_S$ = $c_B$ = 0.1, i.e. $a_{0S}$ = $a_S$ / $(0{+}2)^{0.602}$ = 0.0011, $a_{0B}$ = $a_B / (0{+}2)^{0.602}$ = 0.01252, $c_{0S}$ = $c_S$ / $(0{+}1)^{0.101}$ = 0.1, $c_{0B}$ = $c_B$ / $(0{+}1)^{0.101}$ = 0.1, i.e., the parameters are chosen according to the tuning process, see Spall (2003, Section 7.5). The left hand side of (A.11) is calculated as −0.0114, which satisfies the condition of Corollary A.3, meaning SU outperforms the Bernoulli for $k$ = 1. Now let us check this result with numerical simulation. We approximate the MSEs by averaging over $3\times10^7$ independent sample squared errors. Results are summarized in Table A.2.

In Table A.2, for each iteration count $k$, the MSEs $E\left(\left\|\hat{\boldsymbol{\theta}}_1 - \boldsymbol{\theta}^*\right\|^2\right)$ are approximated by averaging over $3\times10^7$ independent sample squared errors. $P$-values are derived from standard matched-pairs $t$-tests for comparing two population means, which in this case are the MSEs for the Bernoulli and SU. For $k$ = 1, the difference between MSEs under SU and the Bernoulli is −0.0115 (as compared to the theoretical value of −0.0114 computed from the expression in (A.11)), with the corresponding $P$-value being almost 0, which shows a strong indication that SU is preferred to the Bernoulli for $k$ = 1.

**Table A.2**: Results for quadratic loss functions

| Number of iterations | MSE for Bernoulli | MSE for SU | $P$-value |
|---|---|---|---|
| $k$=1 | 0.1913 | 0.1798 | $<10^{-10}$ |
| $k$=5 | 0.2094 | 0.1796 | $<10^{-10}$ |
| $k$=10 | 0.1890 | 0.1786 | $<10^{-10}$ |
| $k$=1000 | 0.0421 | 0.1403 | $>1-10^{-10}$ |

We also notice that the advantage of SU holds for $k$ = 5 and $k$ = 10 in this example. In fact, the better performance of SU for $k$ > 1 has been observed in other examples as well (e.g., Maeda and De Figueiredo, 1997; Spall, 2003, Exercise 7.7). Thus, even though this paper only provides the theoretical foundation for the $k$ = 1 case, it might be possible to generalize the theory to $k$ > 1 provided that $k$ is not too large a number.

A4.2 *Non-quadratic loss function*

Consider the loss function $L(\boldsymbol{\theta}) = t_1^4 + t_1^2 + t_1 t_2 + t_2^2$, where $\boldsymbol{\theta}$ = $[t_1, t_2]^T$, $\sigma^2 = 1$, $\hat{\boldsymbol{\theta}}_0 = [1,1]^T$, the tuning process (see Spall, 2003, Section 7.5) results in $a_S$ = 0.05, $a_B$ = 0.15, $c_S = c_B$ = 1. We estimate the MSEs by averaging over $10^6$ independent sample squared errors. Results are summarized in Table A.3.

**Table A.3**: Results for non-quadratic loss functions

| Number of iterations | MSE for Bernoulli | MSE for SU |
|---|---|---|
| $k$=1 | 1.7891 | 1.5255 |
| $k$=2 | 1.2811 | 1.2592 |
| $k$=5 | 0.6500 | 0.9122 |
| $k$=1000 | 0.0024 | 0.0049 |

In Table A.3, for each iteration count $k$, the MSEs $E(\| \hat{\boldsymbol{\theta}}_1 - \boldsymbol{\theta}^* \|^2)$ are approximated by averaging over $10^6$ independent sample squared errors. Results show that for $k$ = 1, there is a significant advantage of SU over the Bernoulli. But as the sample size increases, this advantage fades out, as we expect given the theory of the asymptotic optimality of the Bernoulli distribution.

**A5. Conclusion**

In this work, we investigate the performance of a non-Bernoulli distribution (specifically, the segmented uniform) for perturbation vectors in one step of SPSA. We show that for certain choices of parameters, non-Bernoulli will be preferred to the Bernoulli as the perturbation distribution for one-iteration SPSA. Furthermore, results in numerical examples indicate that we may generalize the above conclusion to other small sample sizes too, i.e., to two or more iterations of SPSA. In all, this paper gives a theoretical foundation for choosing an effective perturbation distribution when $k$ = 1, and numerical experience indicates favorable results for a limited range of values of $k > 1$. This will be useful for SPSA-based optimization process for which available sample sizes are necessarily small in number.

# Appendix B

## Demonstration of Enhanced Monte Carlo Computation of the Fisher Information for Complex Problems

In practice, it is often the case that closed forms of the Fisher information matrices are not attainable. To solve this problem, we use numerical approximations of the Fisher information matrices. In this appendix, we demonstrate some Monte Carlo methods in computing the Fisher information matrices for complex problems.

**B1. Introduction**

The Fisher information matrix plays an essential role in statistical modeling, system identification and parameter estimation, see Ljung (1999) and Bickel and Doksum (2007, Section 3.4). Consider a collection of $n$ random vectors $\boldsymbol{Z} = [z_1, z_2, \ldots, z_n]^T$, where each $z_i$ is a vector for $i = 1, 2, \ldots, n$. These vectors are not necessarily independent and identically distributed. Let us assume that the probability density/mass function for $\boldsymbol{Z}$ is $p_{\boldsymbol{Z}}(\boldsymbol{\zeta}|\boldsymbol{\theta})$, where $\boldsymbol{\zeta}$ is a dummy matrix representing a possible realization of $\boldsymbol{Z}$; $\boldsymbol{\theta}$ is the unknown $p \times 1$ parameter vector. The corresponding likelihood function is

$$l(\boldsymbol{\theta}|\boldsymbol{\zeta}) = p_{\boldsymbol{Z}}(\boldsymbol{\zeta}|\boldsymbol{\theta}).$$

Letting $L(\boldsymbol{\theta}) = -\log l(\boldsymbol{\theta}|\boldsymbol{Z})$ be the negative log-likelihood function, the $p \times p$ Fisher information matrix $\boldsymbol{F}(\boldsymbol{\theta})$ for a differentiable $L$ is given by

$$\boldsymbol{F}(\boldsymbol{\theta}) \equiv E\left(\frac{\partial L}{\partial \boldsymbol{\theta}} \times \frac{\partial L}{\partial \boldsymbol{\theta}^T}\right), \tag{B.1}$$

where the expectation is taken with respect to the data set $\boldsymbol{Z}$.

Except for relatively simple problems, however, the definition of $\boldsymbol{F}(\boldsymbol{\theta})$ in (B.1) is generally not useful in practical calculation of the information matrix. Computing the expectation of a product of multivariate nonlinear functions is usually a formidable task. A well-known equivalent form follows from the assumption that $L$ is twice continuously differentiable in $\boldsymbol{\theta}$. That is, the Hessian matrix

$$\boldsymbol{H}(\boldsymbol{\theta}) = \frac{\partial^2 L}{\partial \boldsymbol{\theta} \partial \boldsymbol{\theta}^T}$$

is assumed to exist. Furthermore, assume that $L$ is regular in the sense that standard conditions such as in Wilks (1962, pp. 408–411 and 418–419) or Bickel and Doksum (2007, p. 179) hold. Under such conditions, the information matrix is related to the Hessian matrix of $L$ through:

$$\boldsymbol{F}(\boldsymbol{\theta}) = E\left(\frac{\partial^2 L}{\partial \boldsymbol{\theta} \partial \boldsymbol{\theta}^T}\right), \tag{B.2}$$

where the expectation is taken with respect to the data set $\boldsymbol{Z}$. The form of $\boldsymbol{F}(\boldsymbol{\theta})$ in (B.2) is usually more amenable to calculate than the product-based form in (B.1).

In many practical problems, however, closed forms of $\boldsymbol{F}(\boldsymbol{\theta})$ do not exist. In such cases, we need to estimate the Fisher information numerically, see Al-Hussaini and Ahmad (1984), Lei (2010), and Mainassara et al. (2011). Given the equivalent form of $\boldsymbol{F}(\boldsymbol{\theta})$ in (B.2), we can estimate $\boldsymbol{F}(\boldsymbol{\theta})$ using measurements of $\boldsymbol{H}(\boldsymbol{\theta})$. The conventional approach uses resampling-based method to approximate $\boldsymbol{F}(\boldsymbol{\theta})$. In this paper, we demonstrate two other enhanced Monte Carlo methods: feedback-based approach and independent perturbation approach; see Spall (2008). The Monte Carlo computation of $\boldsymbol{F}(\boldsymbol{\theta})$ is discussed in other scenarios too, see Das et al. (2010) where prior information of $\boldsymbol{F}(\boldsymbol{\theta})$ is used in estimation. The remainder of the paper is organized as follows: in Section B2, we introduce methodology of three different approaches discussed in this paper; some relevant theory is summarized in Section B3; section B4 includes two numerical examples and discussions on relative performance of the three methods; a brief conclusion is made in section B5.

## B2. Methodology

B2.1. *Basic resampling-based approach*

We first give a brief review of a Monte Carlo resampling-based approach to compute $\boldsymbol{F}(\boldsymbol{\theta})$, as given in Spall (2005). Let $\boldsymbol{Z}_{\text{pseudo}}(i)$ be a collection of Monte Carlo generated random vectors from the assumed distribution based on the parameters $\boldsymbol{\theta}$. Note that $\boldsymbol{Z}_{\text{pseudo}}(i)$ is one realization of the collection of $n$ random vectors $\boldsymbol{Z}$. Let $\hat{\boldsymbol{H}}_{k|i}$ represent the $k$th estimate of $\boldsymbol{H}(\boldsymbol{\theta})$ at the data set $\boldsymbol{Z}_{\text{pseudo}}(i)$. We generate $\hat{\boldsymbol{H}}_{k|i}$ via efficient simultaneous perturbation (SPSA) principles:

$$\hat{\boldsymbol{H}}_{k|i} = \frac{1}{2}\left[\frac{\delta \boldsymbol{g}_{k|i}}{2c}(\boldsymbol{\Delta}_{k|i}^{-1})^T + \left(\frac{\delta \boldsymbol{g}_{k|i}}{2c}(\boldsymbol{\Delta}_{k|i}^{-1})^T\right)^T\right], \tag{B.3}$$

where $\delta \boldsymbol{g}_{k|i} = \boldsymbol{g}\left(\boldsymbol{\theta} + c\boldsymbol{\Delta}_{k|i} \middle| \boldsymbol{Z}_{\text{pseudo}}(i)\right) - \boldsymbol{g}\left(\boldsymbol{\theta} - c\boldsymbol{\Delta}_{k|i} \middle| \boldsymbol{Z}_{\text{pseudo}}(i)\right)$, $\boldsymbol{g}(\bullet)$ is the exact or estimated gradient function of $L$, depending on the information available; $\boldsymbol{\Delta}_{k|i} = \left[\Delta_{k1|i}, \Delta_{k2|i}, \ldots, \Delta_{kp|i}\right]^T$ is a mean-zero random vector such that the scalar elements are i.i.d. symmetrically distributed random variables that are uniformly bounded and satisfy $E\left(\left|1/\Delta_{kj|i}\right|\right) < \infty$, $\boldsymbol{\Delta}_{k|i}^{-1}$ denotes the vector of inverses of the $p$ individual elements of $\boldsymbol{\Delta}_{k|i}$, and $c > 0$ is a "small" constant. Each $k$ represents different draw of random perturbation vectors $\boldsymbol{\Delta}_{k|i}$. Notice that the second term in the summation in (B.3) is

simply the transpose of the first term. It is deliberately designed this way so that the resulting $\hat{\boldsymbol{H}}_{k|i}$ is symmetric.

The Monte Carlo approach of Spall (2005) is based on a double averaging scheme. The first "inner" average forms Hessian estimates at a given $\boldsymbol{Z}_{\text{pseudo}}(i)$ $(i = 1, 2, \ldots, N)$ from $k = 1, 2, \ldots, M$ values of $\hat{\boldsymbol{H}}_{k|i}$ and the second "outer" average combines these sample mean Hessian estimates across the $N$ values of pseudo data. Therefore, the "basic" Monte Carlo resampling-based estimate of $\boldsymbol{F}(\boldsymbol{\theta})$ in Spall (2005), denoted as $\bar{\boldsymbol{F}}_{M,N}(\boldsymbol{\theta})$, is:

$$\bar{\boldsymbol{F}}_{M,N}(\boldsymbol{\theta}) = \frac{1}{N}\sum_{i=1}^{N}\frac{1}{M}\sum_{k=1}^{M}\hat{\boldsymbol{H}}_{k|i}.$$

This resampling-based estimation method is easy to implement and works well in practice (Spall, 2005). However, this basic Monte Carlo approach could be improved by some extra effort. In the next two subsections, we introduce the use of feedback information and independent perturbation, respectively.

B2.2 *Enhancements through use of feedback*

The feedback ideas for FIM estimation in Spall (2008) are related to the feedback ideas presented with the most updates in Spall (2009), as applied to stochastic approximation. From Spall (2009), it is known that $\hat{\boldsymbol{H}}_{k|i}$ in (B.3) can be decomposed into three parts:

$$\hat{\boldsymbol{H}}_{k|i} = \boldsymbol{H}(\boldsymbol{\theta}) + \boldsymbol{\Psi}_{k|i} + O(c^2), \tag{B.4}$$

where $\boldsymbol{\Psi}_{k|i}$ is a $p \times p$ matrix of terms dependent on $\boldsymbol{H}(\boldsymbol{\theta})$ and $\boldsymbol{\Delta}_{k|i}$. Specifically,

$$\boldsymbol{\Psi}_{k|i}(\boldsymbol{H}) = \frac{1}{2}\boldsymbol{H}\boldsymbol{D}_{k|i} + \frac{1}{2}\boldsymbol{D}_{k|i}^T\boldsymbol{H},$$

where $\boldsymbol{D}_{k|i} = \boldsymbol{\Delta}_{k|i}(\boldsymbol{\Delta}_{k|i}^{-1})^T - \boldsymbol{I}_p$ and $\boldsymbol{I}_p$ is the $p \times p$ identity matrix.

Notice that for any value of $\boldsymbol{H}$, $E(\boldsymbol{\Psi}_{k|i}(\boldsymbol{H})) = \boldsymbol{0}$. If we subtract both sides of (B.4) by $\boldsymbol{\Psi}_{k|i}$ and use $\hat{\boldsymbol{H}}_{k|i} - \boldsymbol{\Psi}_{k|i}$ as an estimate of $\boldsymbol{H}(\boldsymbol{\theta})$, we end up with reduced variance of the Hessian estimate while the expectation of the estimate remains the same. Ultimately, the variance of the estimate of $\boldsymbol{F}(\boldsymbol{\theta})$ is also reduced. Based on this idea, Spall (2008) introduces a feedback-based method to improve the accuracy of the estimate of $\boldsymbol{F}(\boldsymbol{\theta})$. The recursive (in $i$) form of the feedback-based form of the estimate of $\boldsymbol{F}(\boldsymbol{\theta})$, say $\bar{\boldsymbol{F}}'_{M,N}(\boldsymbol{\theta})$, is

$$\bar{\boldsymbol{F}}'_{M,i}(\boldsymbol{\theta}) = \frac{i-1}{i}\bar{\boldsymbol{F}}'_{M,i-1}(\boldsymbol{\theta}) + \frac{1}{iM}\sum_{k=1}^{M}\left[\hat{\boldsymbol{H}}_{k|i} - \boldsymbol{\Psi}_{k|i}(\bar{\boldsymbol{F}}'_{M,i-1}(\boldsymbol{\theta}))\right], \tag{B.5}$$

where $\bar{\boldsymbol{F}}'_{M,0}(\boldsymbol{\theta}) = \boldsymbol{0}$. More recent work regarding the feedback-based approach includes Spall (2009), where the feedback ideas are applied to stochastic approximation.

B2.3 *Enhancements through use of independent perturbation per measurement*

If the $n$ vectors entering each $\boldsymbol{Z}_{\text{pseudo}}(i)$ are mutually independent, the estimation of $\boldsymbol{F}(\boldsymbol{\theta})$ can be improved by exploiting this independence. In particular, for the basic resampling-based approach, the variance of the elements of the individual Hessian estimates $\hat{\boldsymbol{H}}_{k|i}$ can be reduced by decomposing $\hat{\boldsymbol{H}}_{k|i}$ into a sum of $n$ independent estimates, each corresponding to one of the data vectors. A separate perturbation vector can then be applied to each of the independent estimates, which produces variance reduction in the resulting estimate $\bar{\boldsymbol{F}}_{M,N}(\boldsymbol{\theta})$. The independent perturbations above reduce the variance of the elements in the estimate of $\boldsymbol{F}(\boldsymbol{\theta})$ from $O(1/N)$ to $O(1/nN)$.

Similarly, this independent perturbation idea can be applied to the feedback-based approach as well. Besides applying separate perturbation vectors to each of the independent estimates of $\hat{\boldsymbol{H}}_{k|i}$, we also decompose the $\bar{\boldsymbol{F}}'_{M,i-1}(\boldsymbol{\theta})$ in (B.5) into a sum of $n$ independent estimates and then apply the $\boldsymbol{\Psi}_{k|i}$ function to individual estimates to gain feedback information to improve the corresponding independent estimates of $\hat{\boldsymbol{H}}_{k|i}$.

## B3. Theory

The following results are given in Spall (2008) as a theoretical validation for the advantage of the feedback-based approach.

**Lemma B.1**

For some open neighborhood of $\boldsymbol{\theta}^*$, suppose the forth derivative of the log-likelihood function $L''''(\boldsymbol{\theta})$ exists continuously and that $E\left(\left\|L''''(\boldsymbol{\theta})\right\|^2\right)$ is bounded in magnitude. Furthermore, let $E\left(\left\|\hat{\boldsymbol{H}}_{k|i}\right\|^2\right) < \infty$, then for any fixed $M \geq 1$ and all $c$ sufficiently small,

$$E\left(\left\|\bar{\boldsymbol{F}}'_{M,N} - \boldsymbol{F}^* - \boldsymbol{B}(\boldsymbol{\theta}^*)\right\|^2\right) \to 0 \text{ as } N \to \infty,$$

where $\boldsymbol{B}(\boldsymbol{\theta}^*)$ is a bias matrix satisfying $\boldsymbol{B}(\boldsymbol{\theta}^*) = O(c^2)$.

**Theorem B.1**

Suppose that the conditions of the Lemma hold, $p \geq 2$, $E\left(\left\|\boldsymbol{H}(\boldsymbol{\theta}^*)\right\|^2\right) < \infty$, $\boldsymbol{F}^* \geq 0$, and $\boldsymbol{F}^* \neq 0$. Further, suppose that for some $\delta > 0$ and $\delta' > 0$ such that $(1+\delta)^{-1} + (1+\delta')^{-1} = 1$, $E\left(\left\|L''''(\boldsymbol{\theta})\right\|^{2+2\delta'}\right)$ is uniformly bounded in magnitude for all $\boldsymbol{\theta}$ in an open neighborhood of $\boldsymbol{\theta}^*$, $E\left(\left|1/\Delta_{kj|i}^{2+2\delta}\right|\right) < \infty$. Then the accuracy of $\bar{\boldsymbol{F}}'_{M,N}(\boldsymbol{\theta})$ is greater than the accuracy of $\bar{\boldsymbol{F}}_{M,N}(\boldsymbol{\theta})$ in the sense that

$$\lim_{N\to\infty} \frac{E\left[\left\|\bar{\boldsymbol{F}}'_{M,N} - \boldsymbol{F}^*\right\|^2\right]}{E\left[\left\|\bar{\boldsymbol{F}}_{M,N} - \boldsymbol{F}^*\right\|^2\right]} \le 1 + O(c^2). \qquad \text{(B.6)}$$

**Corollary B.1**

Suppose that the conditions of the Theorem hold, $\text{rank}(\boldsymbol{F}^*) \ge 2$, and the elements of $\boldsymbol{\Delta}_{k|i}$ are generated according to the Bernoulli $\pm 1$ distribution. Then, the inequality in (B.6) is strict.

**B4. Numerical study**

In this section, we show the merit of the enhanced Monte Carlo methods over the basic Monte Carlo resampling method. The performance of the estimation is measured by the relative norm of the deviation matrix: $\left\|\boldsymbol{F}_{\text{est}}(\boldsymbol{\theta}) - \boldsymbol{F}_n(\boldsymbol{\theta})\right\| / \left\|\boldsymbol{F}_n(\boldsymbol{\theta})\right\|$, where the standard spectral norm (the largest singular value) is used, $\boldsymbol{F}_n(\boldsymbol{\theta})$ is the true information matrix, and $\boldsymbol{F}_{\text{est}}(\boldsymbol{\theta})$ stands for the estimated information matrix via either the basic or the enhanced Monte Carlo approach, as appropriate. For the purpose of comparison, we test under the cases where the true Fisher information is achievable or the exact Hessian matrix is computable, which are not the type of problems we would actually deal with in practice with these estimation methods.

B4.1 *Example 1—Multivariate normal distribution in a signal-plus-noise setting*

Suppose that the $z_i$ are independently distributed $N(\boldsymbol{\mu}, \boldsymbol{\Sigma}+\boldsymbol{P}_i)$ for all $i$, where $\boldsymbol{\mu}$ and $\boldsymbol{\Sigma}$ are to be estimated and the $\boldsymbol{P}_i$'s are known. This corresponds to a signal-plus-noise setting where the $N(\boldsymbol{\mu}, \boldsymbol{\Sigma})$-distributed signal is observed in the presence of independent $N(\mathbf{0}, \boldsymbol{P}_i)$-distributed noise. The varying covariance matrix for the noise may reflect different quality measurements of the signal. This setting arises, for example, in estimating the initial mean vector and covariance matrix in a state-space model from a cross-section of realizations (Shumway, Olsen, and Levy, 1981), in estimating parameters for random-coefficient linear models (Sun 1982), in small area estimating in survey sampling (Ghosha and Rao 1994), in sensitivity studies (Spall, 1985a; Spall and Chin, 1990), and in nuisance parameter analysis (Spall, 1989).

Let us consider the following scenario: $\dim(z_i) = 4$, $n = 30$, and $\boldsymbol{P}_i = \sqrt{i}\boldsymbol{U}^T\boldsymbol{U}$, where $\boldsymbol{U}$ is generated according to a 4×4 matrix of uniform (0, 1) random variables (so the $\boldsymbol{P}_i$'s are identical except for the scale factor $\sqrt{i}$). Note that once $\boldsymbol{U}$ is generated, it stays constant throughout the study. Let $\boldsymbol{\theta}$ represent the unique elements in $\boldsymbol{\mu}$ and $\boldsymbol{\Sigma}$; hence, $p = 4+4(4+1)/2 = 14$. So, there are $14(14+1)/2 = 105$ unique terms in $\boldsymbol{F}_n(\boldsymbol{\theta})$ that are to be estimated via the Monte Carlo methods (basic or enhanced approaches). The value of $\boldsymbol{\theta}$ used to generate the data is also used as the value of interest in evaluating $\boldsymbol{F}_n(\boldsymbol{\theta})$. This value corresponds to $\boldsymbol{\mu} = \mathbf{0}$ and $\boldsymbol{\Sigma}$ being a matrix with 1's on the diagonal and 0.5's on the off-diagonals. The gradient of the log-likelihood function and the analytical form of the FIM are available in this problem (see Shumway, Olsen, and Levy, 1981).

Throughout the study, elements in perturbation $\mathbf{\Delta}_{k|i}$ have symmetric Bernoulli ± 1 distribution for all $k$ and $i$; $M = 2$; $c = 0.0001$. In each method, we estimate the Hessian matrix in two different approaches: using the gradient of the log-likelihood function or using the log-likelihood function values when the gradient is not available. Results based on 50 independent replications are summarized in Table B.1 ($P$-values correspond to $t$-tests of the comparison between the relative norms of the deviation matrices from two approaches).

Table B.1 indicates that there is statistical evidence for the advantage of the feedback-based Monte Carlo method over the basic Monte Carlo resampling method. The difference between the two methods is more significant when the gradient information of the log-likelihood function is available (row 2) or the number of iterations increases when only likelihood function is available (rows 4).

Keeping all other settings and parameters the same, we now test on the independent perturbation per measurement idea in section B2.3. Table B.2 summarizes the simulation results based on 50 independent realizations ($P$-values correspond to $t$-tests of the comparison between the relative norms of the deviation matrices from two approaches: independent perturbation alone and feedback and independent perturbation combined).

**Table B.1**: Sample mean value of $\|\boldsymbol{F}_{\text{est}}(\boldsymbol{\theta}) - \boldsymbol{F}_n(\boldsymbol{\theta})\| / \|\boldsymbol{F}_n(\boldsymbol{\theta})\|$ with approximate 95% confidence intervals (CIs) shown in brackets. *P*-values based on one-sided *t*-test using 50 independent runs.

| Input Information | Basic Approach | Feedback-based Approach | *P*-value |
|---|---|---|---|
| Gradient Function $N = 40{,}000$ | 0.0104 [0.0096, 0.0111] | 0.0063 [0.0058, 0.0067] | $<10^{-10}$ |
| Log-likelihood Function Only $N = 40{,}000$ | 0.0272 [0.026, 0.0283] | 0.0261 [0.0251, 0.0271] | 0.0016 |
| Log-likelihood Function Only $N = 80{,}000$ | 0.0204 [0.0194, 0.0213] | 0.0191 [0.0184, 0.0198] | $2.52\times10^{-5}$ |

Table B.2 demonstrates the improvement in estimation accuracy when the sample is independent and separate perturbation is applied to each independent measurement. Specifically, the estimation accuracy is improved by independent perturbation alone (column 2) and is improved even more by the combination of independent perturbation and feedback approach (column 3).

**Table B.2**: Sample mean value of $\|\boldsymbol{F}_{\text{est}}(\boldsymbol{\theta})-\boldsymbol{F}_n(\boldsymbol{\theta})\| / \|\boldsymbol{F}_n(\boldsymbol{\theta})\|$ when using independent perturbation per measurement. Approximate 95% CIs shown in brackets. *P*-value based on one-sided *t*-test using 50 independent runs.

| Input Information | Indep. Perturbation Alone | Feedback and Indep. Perturbation | *P*-value |
|---|---|---|---|
| Gradient Function $N = 40{,}000$ | 0.0066 [0.0043, 0.0103] | 0.0062 [0.0044, 0.0097] | $7.622\times10^{-9}$ |

B4.2 *Example 2—Mixture Gaussian distribution*

Mixture Gaussian distribution is of great interest and is popularly used in practical applications (see Wang, 2001; Stein et al., 2002). In this study, we consider a mixture of two scale normal distributions. Specifically, let $\boldsymbol{Z} = [z_1, z_2, \ldots, z_n]^T$ be an independent and identically distributed sequence with probability density function:

$$f(z,\boldsymbol{\theta}) = \lambda \exp\left(-(z-\mu_1)^2/(2\sigma_1^2)\right)\Big/\sqrt{2\pi\sigma_1^2} + (1-\lambda)\exp\left(-(z-\mu_2)^2/(2\sigma_2^2)\right)\Big/\sqrt{2\pi\sigma_2^2}\,,$$

where $\boldsymbol{\theta} = [\lambda,\mu_1,\sigma_1,\mu_2,\sigma_2]^T$. There are 5(5+1)/2 = 15 unique terms in $\boldsymbol{F}_n(\boldsymbol{\theta})$ that are to be estimated. The analytical form of the true Fisher information matrix is not attainable in this case. But the closed form of the Hessian matrix is computable (see Boldea and Magnus 2009). We thus approximate the true Fisher information using the sample average of the Hessian matrix over a large number ($10^6$) of independent replications. This should be a fairly good approximation since the first three decimal digits do not vary as the amount of averaging increases.

In this numerical study, we consider the case where $\boldsymbol{\theta}$ = $[0.2, 0, 1, 4, 9]^T$. As in Example 1, elements in perturbation $\boldsymbol{\Delta}_{k|i}$ have symmetric Bernoulli ± 1 distribution for all $k$ and $i$; $M$ = 2; $c = 0.0001$. In each method, we estimate the Hessian matrix in two different approaches: using the gradient of the log-likelihood function or using the log-likelihood function values only. Results based on 50 independent replications are summarized in Table B.3 ($P$-values correspond to $t$-tests of the comparison between the relative norms of the deviation matrices from two approaches).

**Table B.3**: Sample mean value of $\|\boldsymbol{F}_{\text{est}}(\boldsymbol{\theta})-\boldsymbol{F}_n(\boldsymbol{\theta})\| / \|\boldsymbol{F}_n(\boldsymbol{\theta})\|$ with approximate 95% CIs shown in brackets. $P$-values based on one-sided $t$-test using 50 independent runs.

| Input Information | Basic Approach | Feedback-based Approach | $P$-value |
|---|---|---|---|
| Gradient Function $N = 40{,}000$ | 0.0038 [0.0035, 0.0042] | 0.0013 [0.0011, 0.0015] | $<10^{-10}$ |
| Log-likelihood Function Only $N = 40{,}000$ | 0.0094 [0.0088, 0.01] | 0.0088 [0.0083, 0.0094] | $2.39\times10^{-4}$ |
| Log-likelihood Function Only $N = 80{,}000$ | 0.0065 [0.006, 0.0069] | 0.0059 [0.0054, 0.0063] | $3.6\times10^{-7}$ |

Table B.3 indicates statistical evidence for the advantage of the feedback-based Monte Carlo method over the basic Monte Carlo resampling method. The difference between the performances of the two methods is more significant when gradient information of the log-likelihood function is available (row 2) or the number of iterations increases when only likelihood function is available (row 4).

## B5. Conclusions

This appendix demonstrates two enhanced Monte Carlo methods for estimating the Fisher information matrix: feedback-based approach and independent perturbation approach. Numerical examples show that both of these two methods improve the estimation accuracy as compared to the basic Monte Carlo approach.